\documentclass[12pt]{amsart}
\newcommand{\filename}{homog.tex}   

\usepackage{latexsym, amssymb, amscd, amsfonts}

\usepackage{pstricks}
\usepackage{pst-poly}
\usepackage{multido}
\usepackage[all]{xypic}
\usepackage{ifthen}
\usepackage{longtable}

\DeclareMathAlphabet{\mathpzc}{OT1}{pzc}{m}{it}

\newboolean{showpicture}
\setboolean{showpicture}{true}

\newboolean{bourbaki}
\setboolean{bourbaki}{true}

\newboolean{draft}
\setboolean{draft}{true}

\newboolean{shproofs}
\setboolean{shproofs}{false}

\newcommand{\np}{\medskip\noindent}

\newcommand{\point}{\vspace{3mm}\par\refstepcounter{subsection}\noindent{\bf \thesubsection.} }
\newcommand{\tpoint}[1]{\vspace{3mm}\par\refstepcounter{subsubsection}\noindent{\em #1 {\rm(}{\em \thesubsubsection}{\rm)} ---} }

\newcommand{\bpoint}[1]{\vspace{3mm}\par\refstepcounter{subsection}\noindent{\bf \thesubsection.} {\bf #1.} }

\renewenvironment{equation}{\medskip\noindent\refstepcounter{subsubsection}\makebox[0pt][l]{({\bf\thesubsubsection})}\begin{minipage}[b]{\textwidth}$$}{$$\end{minipage}\medskip\noindent}

\renewcommand{\labelenumi}{{({\em \alph{enumi}})}}

\ifthenelse{\boolean{draft}}
{
}
{
}

\newcommand{\bpf}{\noindent {\em Proof.  }}
\newcommand{\epf}{\qed \vspace{+10pt}}

\newcommand{\Sl}{\operatorname{SL}}
\newcommand{\So}{\operatorname{SO}}
\newcommand{\Gl}{\operatorname{GL}}
\newcommand{\Tr}{\operatorname{Tr}}

\renewcommand{\geq}{\geqslant}
\renewcommand{\leq}{\leqslant}
\newcommand{\Stab}{\operatorname{Stab}}

\newcommand{\st}{\,\,|\,\,}                            
\newcommand{\Osh}{{\mathcal O}}                        
\newcommand{\Fsh}{{\mathcal F}}                        
\newcommand{\Gsh}{{\mathcal G}}                        
\newcommand{\Ksh}{{\mathrm{K}}}                        
\newcommand{\Msh}{{\mathcal M}}                        
\newcommand{\Nsh}{{\mathcal N}}                        
\renewcommand{\Rsh}{{\mathcal R}}                      
\newcommand{\Ssh}{{\mathcal S}}                        

\newcommand{\uu}{\underline{u}}                        
\newcommand{\vv}{\underline{v}}                        
\renewcommand{\ll}{\underline{\lambda}}                
\newcommand{\mm}{\underline{\mu}}                      
\newcommand{\nn}{\underline{\nu}}                      
\newcommand{\useq}{{\underline{\mathbf{u}}}}           
\newcommand{\vseq}{{\underline{\mathbf{v}}}}           
\newcommand{\suseq}{{\!\underline{\mathbf{u}}}}        
\newcommand{\svseq}{{\!\underline{\mathbf{v}}}}        
\newcommand{\smu}{{\!\mu}}                             
\newcommand{\snu}{{\!\nu}}                             
\newcommand{\emptyseq}{\underline{\emptyset}}          

\renewcommand{\H}{\mathrm{H}}                          
\newcommand{\A}{\mathrm{A}}                            
\newcommand{\B}{\mathrm{B}}                            
\newcommand{\C}{\mathrm{C}}                            
\newcommand{\F}{\mathrm{F}}                            
\newcommand{\G}{\mathrm{G}}                            
\renewcommand{\P}{\mathrm{P}}                          
\newcommand{\X}{\mathrm{X}}                            
\newcommand{\M}{\mathrm{M}}                            
\newcommand{\N}{\mathrm{N}}                            
\newcommand{\V}{\mathrm{V}}                            
\newcommand{\U}{\mathrm{U}}                            
\renewcommand{\O}{\mathrm{O}}                          
\newcommand{\Q}{\mathrm{Q}}                            
\newcommand{\Y}{\mathrm{Y}}                            
\newcommand{\Z}{\mathrm{Z}}                            
\newcommand{\bt}{\,\mbox{\tiny$\boxtimes$}\,}          
\newcommand{\W}{\mathcal{W}}                           
\newcommand{\WW}{\mathrm{W}}                           
\newcommand{\id}{\operatorname{id}}                    
\newcommand{\pr}{\operatorname{pr}}                    
\newcommand{\m}{\mathfrak{m}}                          
\newcommand{\delneg}{\Delta^{\!-}}                     
\newcommand{\delpos}{\Delta^{\!+}}                     
\newcommand{\del}{\Delta}                              
\newcommand{\invset}{\Phi}                             
\newcommand{\wo}{w_{0}}                                
\newcommand{\fo}{f_{\circ}}                            
\newcommand{\so}{s_{\circ}}                            
\renewcommand{\c}{\operatorname{c}}                    
\newcommand{\T}{\mathrm{T}}                            
\renewcommand{\L}{\mathrm{L}}                          
\newcommand{\R}{\mathrm{R}}                            
\newcommand{\Serre}{\mathrm{S}}                        
\newcommand{\D}{\mathrm{D}}
\renewcommand{\SS}{\mathrm{S}}                         
\newcommand{\Sym}{\operatorname{Sym}}                  
\newcommand{\spec}{\operatorname{Spec}}                

\newcommand{\gb}{{\mathfrak{b}}}
\newcommand{\gc}{{\mathfrak{c}}}

\newcommand{\gh}{\mathfrak{t}}
\newcommand{\gm}{\mathfrak{m}}
\newcommand{\gn}{\mathfrak{n}}
\newcommand{\gp}{{\mathfrak{p}}}
\newcommand{\gs}{\mathfrak{s}}
\newcommand{\gt}{\mathfrak{a}}

\newcommand{\kostP}{\mbox{{\large$\mathpzc{p}$}}}      
\newcommand{\LRC}{{\mathcal{C}}}                       
\newcommand{\LRP}{{\mathcal{P}}}                       
\newcommand{\die}{\partial}                            
\newcommand{\mult}{\operatorname{mult}}                

\newcommand{\Umult}{\stackrel{\mbox{\tiny mult}}{\bigcup}}

\newcommand{\Ch}{\operatorname{Ch}}                    
\newcommand{\etale}{\'etale }

\newcommand{\Ad}{\operatorname{Ad}}
\newcommand{\tr}{\operatorname{tr}}
\newcommand{\supp}{\operatorname{supp}}
\newcommand{\cP}{\mathcal{P}}
\newcommand{\cQ}{\mathcal{Q}}
\def\vep{\varepsilon}

\newcommand{\Lie}{\operatorname{Lie}}                  
\renewcommand{\ggg}{\mathfrak{g}}                      

\renewcommand{\tt}{\mathfrak{t}}                       

\newcommand{\quot}{/\!\!/}

\newcommand{\E}{\mathrm{E}}                            
\newcommand{\K}{\mathrm{K}}                            

\newcommand{\Span}{\operatorname{span}}                
\newcommand{\rank}{\operatorname{rank}}                

\newcommand{\subsmash}[1]{{\makebox[0cm][l]{\scriptsize$#1$}}} 

\newcommand{\I}{\mathrm{I}}

\newcommand{\opp}{\operatornamewithlimits{\oplus}}
\newcommand{\ott}{\operatornamewithlimits{\otimes}}

\ifthenelse{\boolean{bourbaki}}
{
\renewcommand{\AA}{\mathbf{A}} 
\newcommand{\CC}{\mathbf{C}} 
\newcommand{\PP}{\mathbf{P}} 
\newcommand{\QQ}{\mathbf{Q}} 
\newcommand{\ZZ}{\mathbf{Z}} 
}
{
\renewcommand{\AA}{\mathbb{A}} 
\newcommand{\CC}{\mathbb{C}} 
\newcommand{\PP}{\mathbb{P}} 
\newcommand{\QQ}{\mathbb{Q}} 
\newcommand{\ZZ}{\mathbb{Z}} 
}

\ifthenelse{\boolean{draft}}
{

\newcommand{\remind}[1]{{\bf[#1]}}
\newcommand{\lremind}[1]{{\bf[label:  #1]}}
\newcommand{\bremind}[1]{{\bf[label:  #1]}}

\newcommand{\comment}[1]{{\bf [#1]}}
}
{

\newcommand{\remind}[1]{{}}
\newcommand{\lremind}[1]{{}}
\newcommand{\bremind}[1]{{}}

\newcommand{\comment}[1]{{}}
}

\newcommand{\hiddenproof}[1]{
\ifthenelse{\boolean{shproofs}}{
\medskip
\begin{centering}
\begin{minipage}{0.9\textwidth}
\hrule
\vspace{-0.75\baselineskip}
\small
#1
\vspace{0.25\baselineskip}
\hrule
\end{minipage}\\
\end{centering}
\medskip
}
{
}
}

\newcommand{\Agrid}[4]{%
\pscustom[linecolor=gray]{
\code{8 dict begin /lx}
\dim{#1}
\code{def /ly}
\dim{#2}
\code{def /ux}
\dim{#3}
\code{def /uy}
\dim{#4}
\code{def /sc}
\dim{1}
\code{def}
\code{lx ly ux lx sub uy ly sub rectclip}
\code{/horgrid {3 2 roll 3 sqrt div sc div 2 mul floor 2 div sc mul 3 sqrt mul exch 3 sqrt sc mul 2 div exch {3 copy dup 3 1 roll 0 0 [7 1 roll] [1 0 0 1 0 0] 6 array concatmatrix aload pop pop pop moveto lineto pop} for } def }
\code{/lgrid {3 2 roll 3 sqrt div sc div 2 mul floor 2 div sc mul 3 sqrt mul exch 3 sqrt sc mul 2 div exch {3 copy dup 3 1 roll 0 0 [7 1 roll] [60 cos -60 sin  60 sin 60 cos 0 0] 6 array concatmatrix aload pop pop pop moveto lineto pop} for } def }
\code{/rgrid {3 2 roll 3 sqrt div sc div 2 mul floor 2 div sc mul 3 sqrt mul exch 3 sqrt sc mul 2 div exch {3 copy dup 3 1 roll 0 0 [7 1 roll] [-60 cos 60 sin -60 sin -60 cos 0 0] 6 array concatmatrix aload pop pop pop moveto lineto pop} for } def }
\code{lx ly ux uy horgrid}
\code{60 cos ux mul -60 sin ly mul add 60 sin lx mul 60 cos ly mul add 60 cos lx mul -60 sin uy mul add 60 sin ux mul 60 cos uy mul add lgrid}
\code{60 cos lx mul 60 sin ly mul add -60 sin ux mul 60 cos ly mul add 60 cos ux mul 60 sin uy mul add -60 sin lx mul 60 cos uy mul add rgrid} 
\code{end}
}
}

\newcommand{\wtcirc}[2]{\pscircle[fillcolor=gray,fillstyle=solid](!#1 2 div #1 #2 2 mul add 3 0.5 exp mul 6 div){0.1}}

\newcommand{\bigwtcirc}[2]{\pscircle[linestyle=solid](!#1 2 div #1 #2 2 mul add 3 0.5 exp mul 6 div){0.2}} 
\newcommand{\bigwtsq}[2]{
\pscustom{
\msave
\translate(!#1 2 div #1 #2 2 mul add 3 0.5 exp mul 6 div)
\pspolygon(-0.16,-0.16)(-0.16,0.16)(0.16,0.16)(0.16,-0.16)
\translate(!#1 -1 mul 2 div #1 #2 2 mul add 3 0.5 exp mul 6 div -1 mul)
\mrestore
}
}
\newcommand{\wtmvto}[2]{\moveto(!#1 2 div #1 #2 2 mul add 3 0.5 exp mul 6 div)} 
\newcommand{\wtlnto}[2]{\lineto(!#1 2 div #1 #2 2 mul add 3 0.5 exp mul 6 div)}

\newcommand{\Bgrid}[4]{%
\pscustom[linecolor=gray]{
\code{9 dict begin /lx}
\dim{#1}
\code{def /ly}
\dim{#2}
\code{def /ux}
\dim{#3}
\code{def /uy}
\dim{#4}
\code{def /sc}
\dim{1}
\code{def}
\code{lx ly ux lx sub uy ly sub rectclip}
\code{/diaggrid {lx uy sub 2 sqrt div ux ly sub 2 sqrt div lx ly add 2 sqrt div 1 2 sqrt div sc mul ux uy add 2 sqrt div {3 copy dup 3 1 roll 0 0 [7 1 roll] [45 cos -45 sin  45 sin 45 cos 0 0] 6 array concatmatrix aload pop pop pop moveto lineto pop} for } def }  
\code{/odiaggrid {lx ly add 2 sqrt div ux uy add 2 sqrt div lx uy sub 2 sqrt div 1 2 sqrt div sc mul ux ly sub 2 sqrt div {3 copy dup 4 1 roll exch 0 0 [7 1 roll] [45 cos -45 sin  45 sin 45 cos 0 0] 6 array concatmatrix aload pop pop pop moveto lineto pop} for } def }  
\code{ly sc div floor sc mul sc uy sc div ceiling sc mul {dup lx exch moveto ux exch lineto} for}  
\code{lx sc div floor sc mul sc ux sc div ceiling sc mul {dup ly moveto uy lineto} for}            
\code{diaggrid}                                                                        
\code{odiaggrid}                                                                        
\code{end}
}
}

\newcommand{\bwtcirc}[2]{\pscircle[fillcolor=gray,fillstyle=solid](!#1 #2 add #2){0.1}} 
\newcommand{\bigbwtcirc}[2]{\pscircle[linestyle=solid](!#1 #2 add #2){0.2}} 
\newcommand{\bigbwtsq}[2]{
\pscustom{
\msave
\translate(!#1 #2 add #2)
\pspolygon(-0.16,-0.16)(-0.16,0.16)(0.16,0.16)(0.16,-0.16)
\translate(!#1 #2 add -1 mul #2 -1 mul)
\mrestore
}
}
\newcommand{\bwtmvto}[2]{\moveto(!#1 #2 add #2)} 
\newcommand{\bwtlnto}[2]{\lineto(!#1 #2 add #2)} 
\newcommand{\bigbwttriang}[2]{
\pscustom{
\msave
\translate(!#1 #2 add #2)
\pspolygon(0,.3)(0.2590,-0.15)(-0.25980,-0.15)
\translate(!#1 -1 mul 2 div #1 #2 2 mul add 3 0.5 exp mul 6 div -1 mul)
\mrestore
}
}

\newcommand{\widewtlnto}[2]{%
\lineto(!#1 2 div #1 #2 2 mul add 3 0.5 exp mul 6 div) 
\stroke[linewidth=1.4\pslinewidth,linecolor=white]
\stroke[linewidth=1.3\pslinewidth,linecolor=black]
}
\newcommand{\bigwhitewtcirc}[2]{\pscircle[linestyle=solid,fillcolor=white,fillstyle=solid](!#1 2 div #1 #2 2 mul add 3 0.5 exp mul 6 div){0.25}} 
\newcommand{\littlewhitewtcirc}[2]{\pscircle[linestyle=none,fillcolor=white,fillstyle=solid](!#1 2 div #1 #2 2 mul add 3 0.5 exp mul 6 div){0.2}} 
\newcommand{\biggerwtsq}[2]{
\pscustom{
\msave
\translate(!#1 2 div #1 #2 2 mul add 3 0.5 exp mul 6 div)
\pspolygon(-0.20,-0.20)(-0.20,0.20)(0.20,0.20)(0.20,-0.20)
\translate(!#1 -1 mul 2 div #1 #2 2 mul add 3 0.5 exp mul 6 div -1 mul)
\mrestore
}
}

\newgray{superlightgray}{0.92}

\begin{document}
\pagestyle{plain} \title{{ \large{Cup Products of Line Bundles on Homogeneous Varieties and Generalized PRV 
Components of Multiplicity One}}
}
\author{Ivan Dimitrov}
\address{Department of Mathematics and Statistics, Queen's University, Kingston,
Ontario,  K7L 3N6, Canada} 
\email{dimitrov@mast.queensu.ca} 
\thanks{Research partially supported by NSERC grants}
\author{Mike Roth}
\email{mikeroth@mast.queensu.ca}


\begin{abstract}
Let $\X = \G/\B$ and let $\L_1$ and $\L_2$ be two line bundles on $\X$. Consider the
cup product map 
$$
\H^{d_1} (\X, \L_1) \otimes \H^{d_2} (\X, \L_2) \stackrel{\cup}{\longrightarrow} \H^d(\X, \L),
$$
where $\L = \L_1 \otimes \L_2$ and $d = d_1 + d_2$. 
We answer two natural questions about the map above: When is it a nonzero homomorphism of representations of $\G$?
Conversely, given generic irreducible representations $\V_1$ and $\V_2$, which irreducible components 
of $\V_1 \otimes \V_2$ may appear in the right hand side of the equation above? 
For the first question we find a combinatorial condition expressed in terms of inversion sets of Weyl group elements.
The answer to the second question is especially elegant - 
the representations $\V$ appearing in the right 
hand side of the equation above are exactly the generalized PRV components of $\V_1 \otimes \V_2$ of stable 
multiplicity one. Furthermore, the highest weights $(\lambda_1, \lambda_2, \lambda)$ corresponding to
the representations $(\V_1, \V_2, \V)$ 
fill up the generic faces of the Littlewood-Richardson cone of $\G$ of codimension
equal to the rank of $\G$. In particular, we conclude that the corresponding Littlewood-Richardson coefficients
equal one.

\np
\end{abstract}

\maketitle
\tableofcontents

\date{\today.\hspace{0.5cm}  {\em \filename}}

\maketitle
\tableofcontents

\vspace{0.5cm}

\section{Introduction}
\bpoint{Main problems}
The main object of study of this paper is the cup product map

\begin{equation} \label{eq:cupprod}
\H^{d_1}(\X,\L_1)\otimes\cdots \otimes \H^{d_k}(\X,\L_k) \stackrel{\cup}{\longrightarrow} \H^{d}(\X,\L),
\end{equation}
where $\X=\G/\B$; $\G$ is a semisimple algebraic group over an algebraically closed field of characteristic zero, 
$\B$ is a Borel subgroup of $\G$; 
$\L_1, \ldots, \L_k$ are arbitrary line bundles on $\X$, 
$\L = \L_1 \otimes \cdots \otimes \L_k$;  $d_1, \ldots, d_k$ are non-negative integers, and 
$d = d_1 + \cdots + d_k$.

\np
We assume that both sides of \eqref{eq:cupprod} are nonzero for otherwise the cup product map is the zero map.
Without loss of generality we may also assume that the line bundles $\L_1$,\ldots, $\L_k$, and $\L$ 
are $\G$-equivariant;  then both sides of \eqref{eq:cupprod} carry a natural $\G$-module structure and the 
cup product map is $\G$-equivariant.  
Furthermore by the Borel-Weil-Bott theorem there are irreducible representations $\V_{\smu_1}$, \ldots, $\V_{\smu_k}$,
and $\V_{\smu}$ so that $\H^{d_i}(\X,\L_i)=\V_{\smu_i}^{*}$ for $i=1$,\ldots, $k$, and $\H^{d}(\X,\L)=\V_{\smu}^{*}$
as representations of $\G$.  The dual of \eqref{eq:cupprod} is thus a $\G$-homomorphism

\begin{equation} \label{eq1.1.2}
\V_{\smu} \longrightarrow \V_{\smu_1} \otimes \cdots \otimes \V_{\smu_k}.
\end{equation}

\np
Since $\V_{\smu_1}$, \ldots, $\V_{\smu_k}$, and $\V_{\smu}$ are irreducible representations, 
(\ref{eq:cupprod}) is either surjective or zero; respectively, (\ref{eq1.1.2}) is either injective or zero.  
This leads us naturally to the two main problems of this paper.

\np
{\bf Problem I.} When is (\ref{eq:cupprod}) a surjection of nontrivial representations? 

\np
{\bf Problem II.} For which $(k+1)$-tuples $(\V_{\smu_1}, \ldots, \V_{\smu_k}, \V_{\smu})$ 
of irreducible representations of $\G$ can $\V_{\smu}$ be realized as a component of 
$\V_{\smu_1} \otimes \cdots \otimes \V_{\smu_k}$ via (\ref{eq1.1.2}) 
for appropriate line bundles $\L_1, \ldots, \L_k$ on $\X$? 

\np
We call 
an irreducible representation $\V_{\smu}$ which can be embedded into $\V_{\smu_1} \otimes \cdots \otimes \V_{\smu_k}$ 
via (\ref{eq1.1.2}) a {\em cohomological component of 
$\V_{\smu_1} \otimes \cdots \otimes \V_{\smu_k}$}.  Fixing $\V_{\smu_1}$,\ldots, $\V_{\smu_k}$, a 
variation of Problem II is to determine the cohomological components of 
$\V_{\smu_1} \otimes \cdots \otimes \V_{\smu_k}$.  

\np
With the exception of some quite degenerate cases for Problem II, we provide a complete solution to both problems.

\bpoint{Solution of Problem I} 
\label{sec:setup} 
Fix a maximal torus $\T\subseteq\B$. 
The $\G$-equivariant line bundles on $\X$ are in one-to-one correspondence with the characters of $\T$.
For a character $\lambda$ of $\T$, we denote by $\L_{\lambda}$ the line bundle on $\X$ corresponding to 
the one dimensional representation of $\B$ on which $\T$ acts via $-\lambda$.

\np
The affine action of the Weyl group $\W$ of $\G$ on the lattice of $\T$-characters $\Lambda$ is defined as
$$
w \cdot \lambda = w(\lambda + \rho) - \rho,
$$
where $\rho$, as usual, denotes the half-sum of the roots of $\B$. A character $\lambda \in \Lambda$ is {\em regular}
if there exists a (necessarily unique) element $w \in \W$ such that $w \cdot \lambda$ is a dominant character. 
Following Kostant, \cite[Definition 5.10]{K1}, we define {\em the 
inversion  set $\Phi_w$ of $w \in \W$} as the set $\Phi_w = w^{-1} \Delta^- \cap \Delta^+$, where $\Delta^- = - \Delta^+$
is the set of negative roots of $\G$. 

\np
Let $\lambda_1, \ldots, \lambda_k \in \Lambda$ be the (regular) characters such that $\L_i = \L_{\lambda_i}$ for 
$1 \leq i \leq k$. Then $\L = \L_\lambda$, where $\lambda = \sum_{i = 1}^k \lambda_i$. Assume that $\lambda$ is also 
regular and denote by $w_1$, \ldots, $w_k$, and $w$ the Weyl group elements for which 
$w_i \cdot \lambda_i$ for $1 \leq i \leq k$
and $w \cdot \lambda$ are dominant. With this notation we prove the following criterion for surjectivity of (\ref{eq:cupprod}).

\np
{\em Theorem I ---} For any semisimple $\G$, if $\H^{d}(\X,\L_{\lambda})\neq 0$, then 
the cup product map \eqref{eq:cupprod} is surjective if and only if

\begin{equation}\label{eqn:liningup}
\invset_{w} = \bigsqcup_{i=1}^{k} \invset_{w_i}.
\end{equation}

\np
Studying the structure of $(k+1)$-tuples $(w_1,\ldots, w_k,w)$ satisfying \eqref{eqn:liningup} is an interesting
combinatorial problem which we do not address here.  
For some open questions concerning \eqref{eqn:liningup} see the expository article \cite{dr}.

\bpoint{Solution of Problem II} 
We say that a component $\V_{\smu}$ of $\V_{\smu_1}\otimes\cdots\otimes\V_{\smu_k}$ 
has {\em stable multiplicity one} if 
the multiplicity of $\V_{\!m\mu}$ in $\V_{\!m\mu_1}\otimes\cdots\otimes\V_{\!m\mu_k}$ is one for all $m\gg 0$.
We say that $\V_{\smu}$ is a {\em generalized PRV component of $\V_{\smu_1}\otimes\cdots\otimes\V_{\smu_k}$} if
there exist $w_1$,\ldots, $w_k$, and $w\in\W$ such that 
$w^{-1}  \mu = w_1^{-1}  \mu_1 + \cdots + w_k^{-1} \mu_k$.
(See \S\ref{sec:PRVdiscuss} and \S\ref{sec:componentdefs} for further discussion of these conditions.)

\np
{\em Theorem II ---}\label{sec:thmIII}

\begin{enumerate}
\item Let $\V_{\smu}$ be a cohomological component of $\V_{\smu_1}\otimes\cdots\otimes \V_{\smu_k}$.  Then 
$\V_{\smu}$ is a generalized PRV component of $\V_{\smu_1}\otimes\cdots\otimes \V_{\smu_k}$ of stable multiplicity one.

\medskip
\item Conversely, assume that 
$\V_{\smu}$ is a generalized PRV component of $\V_{\smu_1}\otimes\cdots\otimes \V_{\smu_k}$ of stable multiplicity one.
If, in addition, one of the following holds:

\medskip
\begin{itemize}
\item[({\em i})] at least one of $\mu_1$,\ldots, $\mu_k$ or $\mu$ is strictly dominant,
\item[({\em ii})] $\G$ is a simple classical group or a product of simple classical groups, 
\end{itemize}

\np
then $\V_{\smu}$ is a cohomological component of $\V_{\smu_1}\otimes\cdots\otimes \V_{\smu_k}$.
\end{enumerate}

\np
It is unfortunate that in part ({\em b}) above we require condition ({\em i}) or ({\em ii}). 
Indeed, we believe that we do not need these conditions 
but we impose them due to our inability to overcome a combinatorial problem. 

\bpoint{Representation-theoretic implications of Theorem II} \label{sec:thmIIdiscuss}
The representation-theoretic significance of
Theorem II is twofold: it provides both a geometric construction of special components of a tensor product
via the Bott theorem and a new way of generalizing the classical PRV component.

\np
The Borel-Weil-Bott theorem provides a geometric realization of every irreducible representation of $\G$ as the
cohomology (in any degree) of an appropriate line bundle on $\X$. In particular, every irreducible representation
equals the space of global sections of a unique line bundle on $\X$.
In this sense the Borel-Weil theorem (the statement about cohomology in degree zero) 
suffices since the Bott theorem 
(the statement about higher cohomology) yields the same representations. 
However, in addition to being representations, the cohomology groups carry a ring structure induced from
the cup product.
Theorem II employs this structure to give a geometric realization of certain components 
of a tensor product of representations.
As far as we know this is the first use of the Bott theorem
for a geometric construction of representations in the case when $\G$ is a semisimple algebraic 
group over a field of characteristic zero. 

\np
We are borrowing the term "generalized PRV component" from the
case when $k=2$. In \cite{prv} Parthasarathy, Ranga Rao, and Varadarajan established that 
if $\mu$ is in the $\W$-orbit of $\mu_1+\wo\mu_2$ (where $\wo$ denotes the longest element of $\W$), 
then $\V_{\smu}$ is a component
of $\V_{\smu_1} \otimes \V_{\smu_2}$. Moreover, they proved that $\V_\mu$ has multiplicity one in,
and is the smallest component of, $\V_{\smu_1} \otimes \V_{\smu_2}$. 
It is true more generally that if $\mu$ is in the $\W$-orbit of $\mu_1+ v\mu_2$ 
(where $v$ is now an arbitrary element of $\W$) then $\V_{\smu}$ is again a component of 
$\V_{\smu_1}\otimes\V_{\smu_2}$;  this was established independently 
by Kumar \cite{ku2} and Mathieu \cite{m}.  

\np
Unlike the
original PRV component, a generalized PRV component $\V_{\smu}$  may have multiplicity greater than one in 
$\V_{\mu_1} \otimes \V_{\mu_2}$.
However, 
by Theorem II, every cohomological component is a 
generalized PRV component of stable multiplicity one.
The cohomological components also retain an aspect of the minimality  of the original PRV component:
every cohomological component of $\V_{\smu_1}\otimes\cdots\otimes \V_{\smu_k}$
is extreme among all components of $\V_{\smu_1}\otimes\cdots\otimes \V_{\smu_k}$
(see Proposition \ref{prop:cohomologicalvertex}). 
These properties of cohomological components suggest that they may be viewed as
 the ``true'' analog of the original PRV component.

\np
The following examples illustrate Theorem II and Proposition \ref{prop:cohomologicalvertex} when $k=2$.


\bigskip
\bigskip


\ifthenelse{\boolean{showpicture}}
{
\noindent
\begin{tabular}{ccc}
\begin{pspicture}(-0.8,1.5)(3.8,6.5)
\rput(1.5,1.1){\scriptsize {$\V_{3,5}\otimes \V_{1,2}$}}
\rput(1.5,6.9){\scriptsize {$\Sl_3$}}
\SpecialCoor
\pscustom[fillstyle=solid,fillcolor=superlightgray,linecolor=gray]{
\wtmvto{1}{4}
\wtlnto{5}{2}
\wtlnto{6}{3}
\wtlnto{4}{7}
\wtlnto{2}{8}
\wtlnto{0}{6}
\wtlnto{1}{4}
}
\Agrid{-0.8}{1.5}{3.8}{6.5}
\wtcirc{0}{6}
\wtcirc{1}{4}
\wtcirc{1}{7}
\wtcirc{2}{5}
\wtcirc{2}{8}
\wtcirc{3}{3}
\wtcirc{3}{6}
\wtcirc{4}{4}
\wtcirc{4}{7}
\wtcirc{5}{2}
\wtcirc{5}{5}
\wtcirc{6}{3}
\bigwtcirc{4}{7}
\bigwtcirc{6}{3}
\bigwtcirc{2}{8}
\bigwtcirc{5}{2}
\bigwtcirc{0}{6}
\bigwtcirc{1}{4}
\end{pspicture}
&
\begin{pspicture}(-1,-0.5)(5,4.5)
\rput(2,4.9){\scriptsize {$\So_5$}}
\rput(2,-0.9){\scriptsize {$\V_{\!\rho}\otimes \V_{\!\rho}$}}
\SpecialCoor
\pscustom[fillstyle=solid,fillcolor=superlightgray,linecolor=gray]{
\bwtmvto{0}{0}
\bwtlnto{3}{0}
\bwtlnto{2}{2}
\bwtlnto{0}{4}
\bwtlnto{0}{0}
}
\Bgrid{-0.5}{-0.5}{4.5}{4.5}
\bwtcirc{0}{0}
\bwtcirc{1}{0}
\bwtcirc{2}{0}
\bwtcirc{3}{0}
\bwtcirc{2}{2}
\bwtcirc{0}{4}
\bwtcirc{0}{2}
\bwtcirc{1}{2}
\bigbwtcirc{0}{0}
\bigbwtcirc{3}{0}
\bigbwtcirc{2}{2}
\bigbwtcirc{0}{4}
\bigbwtsq{0}{2}
\bigbwtsq{1}{2}
\bigbwttriang{1}{0}  
\end{pspicture}
&
\begin{pspicture}(0.7,1)(5.8,6)
\rput(3.25,6.4){\scriptsize {$\Sl_3$}}
\rput(3.25,0.6){\scriptsize {$\V_{7,2}\otimes \V_{1,3}$}}
\SpecialCoor
\pscustom[fillstyle=solid,fillcolor=superlightgray,linecolor=gray]{
\wtmvto{4}{1}
\wtlnto{6}{0}
\wtlnto{9}{0}
\wtlnto{10}{1}
\wtlnto{8}{5}
\wtlnto{6}{6}
\wtlnto{3}{3}
\wtlnto{4}{1}
}
\Agrid{0.7}{1}{5.8}{6}
\wtcirc{3}{3}
\wtcirc{4}{1}
\wtcirc{4}{4}
\wtcirc{5}{2}
\wtcirc{5}{5}
\wtcirc{6}{0}
\wtcirc{6}{3}
\wtcirc{6}{6}
\wtcirc{7}{1}
\wtcirc{7}{4}
\wtcirc{8}{2}
\wtcirc{8}{5}
\wtcirc{9}{0}
\wtcirc{9}{3}
\wtcirc{10}{1}
\bigwtcirc{8}{5}
\bigwtcirc{10}{1}
\bigwtcirc{6}{6}
\bigwtcirc{3}{3}
\bigwtcirc{4}{1}
\bigwtsq{8}{2}
\end{pspicture}
\\
\rule{0cm}{0.6cm}\\
\multicolumn{3}{c}{
\begin{pspicture}(0,-0.2)(14,1.9)
\pspolygon[linecolor=gray](0,-0.15)(0,1.9)(12,1.9)(12,-0.15)
\pscircle[fillcolor=gray,fillstyle=solid](0.5,1.7){0.1}
\rput(3.05,1.7){\tiny -- component of the tensor product}
\pscircle[fillcolor=gray,fillstyle=solid](0.5,1.2){0.1}
\pscircle(0.5,1.2){0.2}
\rput(6.1,1.2){\tiny -- cohomological component/generalized PRV component of stable multiplicity one}
\pscircle[fillcolor=gray,fillstyle=solid](0.5,0.7){0.1}
\pscustom{
\msave
\translate(0.5,0.7)
\pspolygon(-0.16,-0.16)(-0.16,0.16)(0.16,0.16)(0.16,-0.16)
\translate(-0.5,-0.7)
\mrestore
}
\rput(4.77,0.7){\tiny -- generalized PRV component of multiplicity greater than one}
\pscircle[fillcolor=gray,fillstyle=solid](0.5,0.1){0.1}
\pscustom{
\msave
\translate(0.5,0.1)
\pspolygon(0,.3)(0.2590,-0.15)(-0.25980,-0.15)
\translate(-0.5,-0.1)
\mrestore
}
\rput(5.92,0.1){\tiny -- generalized PRV component of multiplicity one, but not stable multiplicity one}
\end{pspicture}
}\\
\multicolumn{3}{c}{Figure 1} \\
\end{tabular}

}
{}

\bpoint{The combinatorics of the Littlewood-Richardson cone} 
The {\em Littlewood-Richardson cone} $\LRC(k)$ is defined as the the rational cone 
 generated by $(\mu_1,\ldots, \mu_k,\mu)$ such that $\V_{\smu}$ is a component of 
$\V_{\smu_1}\otimes\cdots\otimes\V_{\smu_k}$.  Let $n$ denote the rank of $\G$. 
Combining Theorem I with results of Ressayre we obtain:

\newpage
\np
{\em Theorem IV  ---}

\begin{enumerate}
\item 
Assume that $(w_1,\ldots, w_k,w)\in\W^{k+1}$ satisfies the conditions 
$\invset_{w}=\sqcup_{i=1}^{k} \invset_{w_i}$ and $\cap_{i=1}^{k}[\Omega_{w_i}]\cdot[\X_{w}]=1$. Then the  
integral weights $(\mu_1,\ldots, \mu_k,\mu)\in\LRC(k)$ such that $\V_{\smu}$ is a cohomological component  
of $\V_{\smu_1}\otimes\cdots\otimes\V_{\smu_k}$ with respect to $(w_1,\ldots, w_k, w)$ are exactly the lattice 
points of a codimension $n$ face of $\LRC(k)$.

\medskip
\item
Conversely, given a codimension $n$ face $\F$ of $\LRC(k)$ which intersects the locus of strictly dominant
weights, there exists $(w_1,\ldots, w_k,w)\in\W^{k+1}$ such that for every lattice point 
$(\mu_1,\ldots, \mu_k,\mu)$ on $\F$, $\V_{\smu}$ is a cohomological component of $\V_{\smu_1}\otimes\cdots\otimes
\V_{\smu_k}$ with respect to $(w_1,\ldots, w_k, w)$.
\end{enumerate}

\np
We also establish multiplicity bounds for the Littlewood-Richardson coefficients. 
Let $\kostP_{\!\delpos}$ be the Kostant partition function and,
for any integer $r\geq 1$, let $\kostP_{\!r\delpos}$ be the partition function whose generating function is 
the $r$-th power
of the generating function for $\kostP_{\!\delpos}$.    For any $w_1$,\ldots, $w_k$, and $w$ such that 
\eqref{eqn:liningup} holds  we prove the inequality (see \eqref{eqn:cohbound})

\begin{equation} \label{eq:introbound}
\mult(\V_{\smu},\V_{\smu_1}\otimes\cdots\otimes\V_{\smu_k})\leq\kostP_{\!(k-1)\delpos}
\left({ \sum_{i=1}^{k}w_i^{-1}\mu_i - w^{-1}\mu }\right),
\end{equation}

\np
where $\mult(\V_{\smu},\V_{\smu_1}\otimes\cdots\otimes\V_{\smu_k})$ denotes the multiplicity of $\V_{\smu}$ in 
$\V_{\smu_1}\otimes\cdots\otimes\V_{\smu_k}$.  Each of these bounds has the same asymptotic order of growth
along rays as the multiplicity function.  
Applied to faces of $\LRC(2)$, \eqref{eq:introbound} yields
a strengthening of the statement of Theorem IV({\em b}). Namely, if
 $\F$ is a face of $\LRC(2)$ of codimension $n-1$ which intersects the locus of strictly 
dominant weights,
then for every triple $(\mu_1,\mu_2,\mu)$ of integral weights in $\F$, the multiplicity 
of $\V_{\smu}$ in $\V_{\smu_1}\otimes\V_{\smu_2}$ is at most one.
After scaling (i.e., modulo the saturation problem) $\V_{\smu}$ has multiplicity exactly one.
Even this statement can be strengthened, see Corollary \ref{cor:nearvertex}.

\bpoint{Other results} In conclusion we mention several other results which may be of independent interest.

\np
{\bf The cup product and Schubert calculus.}
Recall that a basis for the cohomology ring 
$\H^{*}(\X,\ZZ)$ of $\X=\G/\B$ is given by the classes of the Schubert cycles 
$\{[\X_{w}]\}_{w\in \W}$ indexed 
by the elements of the Weyl group $\W$. 
The dual basis $\{[\Omega_{w}]\}_{w\in \W}$, is given by $\Omega_{w}:=\X_{\wo w}$.  
With the notation of \S\ref{sec:setup} we prove the following:

\np
{\em Theorem III ---} For any semisimple algebraic group $\G$,
\begin{enumerate}
\item[({\em a})] if $\displaystyle{\bigcap_{i=1}^{k}[\Omega_{w_i}]\cdot[\X_{w}]=1}$ then the cup product map 
\eqref{eq:cupprod} is surjective;
\item[({\em b})] if $\displaystyle{\bigcap_{i=1}^{k}[\Omega_{w_i}]\cdot[\X_{w}]=0}$ then the cup product map 
\eqref{eq:cupprod} is zero.
\end{enumerate}

\np
We use Theorem III as stated above and a variation of its proof to prove Theorem I.  
If $\G$ is a product of  simple classical groups, then one can show that
condition \eqref{eqn:liningup} implies
that ${\cap_{i=1}^{k}[\Omega_{w_i}] \cdot [\X_{w}]=1}$.  
We do not know if condition \eqref{eqn:liningup} implies that the intersection number is one for the 
exceptional groups of rank greater than two.

\np
{\bf Diagonal Bott-Samelson-Demazure-Hansen varieties.} We construct a class of varieties which
generalize the Bott-Samelson-Demazure-Hansen varieties. These varieties may find applications beyond 
this paper. In particular we express them as resolutions of singularities of the total space of 
intersections of translates of Schubert varieties, cf. Theorem \ref{thm:fibres}.  Other notable results related
to this construction include
Lemma \ref{lem:key} which controls the multiplicity of 
cohomological components, and Theorem \ref{thm:hammer} which provides a new proof of the necessity 
of the inequalities determining the Littlewood-Richardson cone.

\np
{\bf On invariants of $\Sym^.(\ggg)$.}
The proof of Theorem II({\em b}) requires the solution of a combinatorial problem which we reduce to a problem
of finding a parabolic subalgebra of $\ggg$ with particular properties.  The proof of the existence
of such a subalgebra  uses methods very different from the rest of the paper and is
provided in a self-contained Appendix. The spirit of the result, 
Theorem \ref{thm:classical}, is very classical. Notably, this is 
the only place in the paper where a case-by-case analysis is carried out. This is not an accident --
Theorem \ref{thm:classical} does not hold when $\ggg$ is an exceptional Lie algebra of rank greater than two.

\bpoint{Acknowledgments} 
We thank 
P.\ Belkale, W.\ Fulton, B.\ Kostant, S.\ Kumar, O.\ Mathieu, K.\ Purbhoo, and D.\ Wehlau  for numerous useful
discussions.
Ivan Dimitrov acknowledges excellent working conditions at the Max-Planck Institute.
Mike Roth acknowledges the hospitality of the University of Roma III.

\section{Notation and background results}
\label{sec:background}

\bpoint{Notation and conventions}
\label{sec:notation}
The ground field is algebraically closed of characteristic zero. 
Throughout the paper we fix  a semisimple 
connected algebraic group $\G$, a Borel subgroup $\B \subset \G$, and a maximal torus $\T \subset \B$.
All parabolic subgroups we consider contain $\T$.
The Lie algebras of algebraic groups are denoted by fraktur letters, e.g. $\ggg$, $\gb$, $\tt$, etc.
We use the term "$\G$-module" instead of "representation of $\G$" to avoid differentiating between representations 
of algebraic groups and modules over the respective Lie algebras; likewise, since $\T$ is fixed, we use 
the term "weight" both for characters of $\T$ and weights of $\tt$; in particular we only consider
integral weights of $\tt$.

\np
The point $w\B/\B \in \X_w \subseteq \X = \G/\B$, where $w \in \W$ and $\X_w$ is the corresponding Schubert variety, 
is denoted by $w$ for short. If $\M=\G/\P$ for some parabolic $\P$ we similarly use $w$ to indicate the point 
$w\P/\P\in\M$.

\np
If $\Lambda$ is the lattice of weights of $\T$ we denote the group ring of $\Lambda$ by $\ZZ[\Lambda]$, i.e. 
$$
\ZZ[\Lambda] = \left\{ \sum_{i=1}^k c_i e^{\lambda_i} \, | \, c_i \in \ZZ, \lambda_i \in \Lambda\right\}.
$$
For a $\T$--module $\mathcal{M}$, the {\it formal character of $\mathcal{M}$} is 
$$
\Ch \mathcal{M} = \sum_{\lambda \in \Lambda} \dim \mathcal{M}^\lambda e^\lambda \in \ZZ[\Lambda],
$$
where $\mathcal{M}^\lambda = \{ x \in \mathcal{M} \, | \, t \cdot x = \lambda(t) x {\text { for every }} t \in \tt\}$.
All formal characters discussed in this paper are contained in
$\ZZ[\Delta]$. For a subset $\Phi \subseteq \Delta$, the formal character of 
$\oplus_{\alpha \in \Phi} \ggg^\alpha$ is denoted by $\langle \Phi \rangle$, i.e.
$$
\langle \Phi \rangle = \sum_{\alpha \in \Phi} e^\alpha.
$$

\np
If $w$ is an element of the Weyl group $\W$, then $\ell(w)$ means the length of 
any minimal expression giving $w$ as a product of simple reflections.  
 If $\vv$ is a word in the simple reflections, then $\ell(\vv)$ is the 
number of reflections in the word.  
Note that, if $\vv$ is a word in simple reflections, and $v \in \W$ is the corresponding element of the Weyl group,
then $\ell(\vv) = \ell(v)$ if and only if $\vv$ is a reduced word.  If $\vv = s_{i_1} \ldots s_{i_m}$ is a non-empty
word, we denote by $\vv_R$ the word $s_{i_1} \ldots s_{i_{m-1}}$ obtained from $\vv$ by dropping the rightmost
reflection in $\vv$.
If $\vseq=(\vv_1,\ldots, \vv_k)$ is a sequence of words then
we set $\ell(\vseq)=\sum_{i=1}^{k} \ell(\vv_i)$.  

\np
The following notation is used consistently throughout the paper. 

\medskip

\begin{longtable}{llll}
$\sqcup_{i=1}^k$ & -- && disjoint union\\
$\kappa(\cdot,\cdot)$  & -- && the Killing form  of $\G$ \\
$\Lambda$, $\Lambda^+$ & -- && weight lattice and cone of dominant weights\\
$\V_{\smu}$  & -- &&irreducible $\G$-module of highest weight $\mu$ \\
$\mult(\V_{\smu}, \V)$ & -- && the multiplicity of $\V_{\smu}$ in $\V$\\ 
$\LRC(k)$ & -- && the Littlewood-Richardson cone in $\Lambda^{k+1}$, see \S\ref{sec:LR-cone}\\
$\mathcal{P}(\mu_1, \ldots, \mu_k)$ & -- && the Littlewood-Richardson polytope, see \S\ref{sec:LR-cone}\\
$\{\alpha_1, \ldots, \alpha_n\}$ & -- && base of simple roots of $\B$\\
$\W$ & -- && Weyl group of $\ggg$\\
$w \cdot \lambda$ & -- && $w(\lambda + \rho) - \rho$, the result of the affine action of $w \in \W$ on $\lambda \in \Lambda$\\
$s_i$  & -- && simple reflection along $\alpha_i$\\
$\P_{\!\alpha_i} $  & -- && the minimal parabolic subgroup of $\G$ associated to $\alpha_i$ \\
$\P_{\I}$  & -- && the minimal parabolic subgroup of $\G$ associated to a set $\I$ of simple roots\\
$\W_\P$ & -- && the Weyl group of a parabolic subgroup $\P \subseteq \G$\\
$\Span_{\ZZ_{\geq 0}} \Phi$  & -- && the set of non-negative integer combinations 
of elements of $\Phi \subseteq \Delta$\\
$\uu$ or $\vv$   & -- &&a word $s_{i_{1}} \ldots s_{i_{m}}$ in the simple reflections 
of the Weyl group\\
$u$, $v$ & -- && the element of $\W$ corresponding to $\uu$ or $\vv$\\
$\vv_{R}$  & -- &&the word obtained by dropping the rightmost reflection of $v$\\
$\useq$ or $\vseq$  & -- &&a sequence $(\uu_1,\ldots, \uu_k)$ or $(\vv_1,\ldots, \vv_k)$ of words\\
$\invset_{w}$ & -- && $w^{-1} \Delta^- \cap \Delta^+$, the inversion set of $w \in \W$ \\
$\langle \invset \rangle$ & -- && $\sum_{\alpha \in \invset} e^{\alpha}$, the formal character of $\oplus_{\alpha \in \invset} \ggg^\alpha$,
where $\invset \subset \Delta$\\
$\ell(w)$  & -- &&the length of $w \in \W$\\
$\L_{\lambda}$  & --&& the line bundle on $\X$ corresponding to $\B$-module
on which $\T$ acts via $-\lambda$ \\
$\N$  & --&& the dimension of $\X$ \\
$\pi_i$  & --&& the projection $\pi_i\colon \X\longrightarrow \G/\P_{\alpha_i}$ (a $\PP^1$-fibration) \\
 \end{longtable}

\bpoint{Inversion sets}\label{sec:inversionset}
Let $\delpos$ be the set of positive roots of $\ggg$ (with respect to $\B$). Following Kostant
\cite[Definition 5.10]{K1}, for any element $w$ of the Weyl group $\W$ 
we define $\invset_{w}$, the {\em inversion set} of $w$, to be the set of positive roots  sent to 
negative roots by $w$, i.e., 

\begin{equation}
\invset_{w} := w^{-1}\delneg\cap \delpos.
\end{equation}

\np
For a subset $\Phi$ of $\delpos$, we set $\Phi^{\c}:=\delpos\setminus\Phi$. 
We will need the following formulas, which follow easily from the definition:

\begin{equation}\label{eqn:weightcomplement}
\invset_{\wo w} = \invset_{w}^{\c};
\end{equation}

\begin{equation}\label{eqn:tgt-at-winv}
w^{-1}\Delta^{+} = \invset_{w}^{\c} \sqcup -\invset_{w};
\end{equation}

\begin{equation}\label{eqn:rootsum}
w^{-1}\cdot 0 
=
w^{-1}\rho-\rho
= -\sum_{\alpha\in\invset_{w}} \alpha.
\end{equation}

\bpoint{Generalized PRV components} \label{sec:PRVdiscuss}
For fixed dominant weights $\mu_1$, $\mu_2$, and $\mu$   
it is clear that the two conditions

\begin{enumerate}
\renewcommand{\labelenumi}{{({\em \roman{enumi}})}}
\item[({\em a})] there exist $w_1$, $w_2$, and $w$ in $\W$ such that $w^{-1}\mu=w_1^{-1}\mu_1+w_2^{-1}\mu_2$,
\item[({\em b})] there exists $v$ in $\W$ such that $\mu$ is in the $\W$-orbit of $\mu_1+v\mu_2$
\end{enumerate}

\np
are equivalent. If these conditions are satisfied we call $\V_{\smu}$ a {\em generalized PRV component
of $\V_{\smu_1}\otimes\V_{\smu_2}$}.

\np
As is suggested by the name, but is far from obvious from the definition, every generalized PRV component
of $\V_{\smu_1}\otimes\V_{\smu_2}$ is in fact a component of the tensor product 
$\V_{\smu_1}\otimes\V_{\smu_2}$ of $\G$-modules.  This 
was first proved when $v=\wo$ (i.e., when $\mu$ is in
the $\W$-orbit of $\mu_1+\wo \mu_2$) in \cite{prv}.  
In the literature this component is referred to simply as the {\em PRV component}.
The general case, that $\V_{\smu}$ is a 
component of $\V_{\smu_1}\otimes\V_{\smu_2}$ for an arbitrary $v$,
became known as the PRV conjecture, and was established independently by Kumar \cite{ku2} and Mathieu \cite{m}.  

\np
In the present paper we extend the notion of generalized PRV component to components of the tensor product of
$k$ irreducible $\G$-modules for $k\geq 2$.
We call $\V_{\smu}$ a {\em generalized PRV component of $\V_{\smu_1}\otimes\cdots\otimes\V_{\smu_k}$} if there
exist $w_1$,\ldots, $w_k$, and $w$ in $\W$ such that $w^{-1}\mu=\sum_{i=1}^{k} w_i^{-1}\mu_i$.
A straightforward induction from the case $k=2$ implies that every generalized PRV component of 
$\V_{\smu_1}\otimes\cdots\otimes\V_{\smu_k}$ is a component of the tensor product 
$\V_{\smu_1}\otimes\cdots\otimes\V_{\smu_k}$  of $\G$-modules.
We record the special case when $\mu=0$ for use in the proof of Theorem I.

\tpoint{Lemma} \label{lem:PRVsymmetric}
For any dominant weights $\mu_1$,\ldots, $\mu_{k}$, and Weyl group elements $w_1$,\ldots, $w_{k}$,
if $\sum_{i=1}^{k}w_i^{-1}\mu_i = 0$ then $(\V_{\smu_1}\otimes\cdots\otimes\V_{\smu_{k}})^{\G}\neq 0$.

\bpoint{Borel-Weil-Bott theorem} \label{sec:bwb} 
Suppose that $\lambda$ is a regular weight, so there is a unique $w\in \W$ with $w\cdot\lambda\in\Lambda^{+}$.
The Borel-Weil-Bott theorem identifies the cohomology of the line bundle $\L_{\lambda}$ on $\X$ as $\G$-modules:

$$
\H^{d}(\X,\L_{\lambda}) = \left\{{
\begin{array}{cl}
\V_{w\cdot\lambda}^{*} & \mbox{if $d=\ell(w)$} \\
0 & \mbox{otherwise.} \rule{0cm}{0.6cm}\\
\end{array}
}\right.
$$

\np
If $\lambda$ is not a regular weight then the cohomology of $\L_{\lambda}$ is zero in all degrees.

\bpoint{Serre Duality on $\X$} \label{sec:serredual}
For any weight $\lambda$ set $\Serre(\lambda)=-\lambda-2\rho$.  Since the canonical bundle
$\Ksh_{\X}$ of $\X$ is equal to $\L_{-2\rho}$ 
we see that $\L_{\Serre(\lambda)}=\Ksh_{\X}\otimes \L_{\lambda}^{*}$.
In other words, $\Serre$ is the function that for each weight $\lambda$ returns the weight $\Serre(\lambda)$ 
of the line bundle Serre dual to $\L_{\lambda}$;  the map $\Serre$ is clearly an involution. 
Let $w$ be any element of the Weyl group and $\lambda$ any weight.  
A straightforward computation shows that $\Serre$ commutes with the affine action of the Weyl group,
i.e.\ that $w\cdot\Serre(\lambda)=\Serre(w\cdot\lambda)$.

\tpoint{Lemma} \label{lem:serredual}
If $\lambda$ is a regular weight and $w$ the unique element of the Weyl group with $w\cdot\lambda\in\Lambda^{+}$ 
then $(\wo w)\cdot\Serre(\lambda)\in\Lambda^{+}$.

\np
\bpf
If $\mu$ is a dominant weight then $\V_{\smu}^{*}=\V_{-\wo \mu}$.
Therefore if $w\cdot\lambda=\mu\in\Lambda^{+}$ then 

\begin{equation}\label{eq:weightserreduality}
(\wo w)\cdot\Serre(\lambda)=\wo\cdot\Serre(w\cdot\lambda)=\wo\cdot\Serre(\mu)=-\wo\mu \in\Lambda^{+}. \qed
\end{equation}

\np
Since $\ell(\wo w)=\N-\ell(w)$,
the calculation above fits in neatly
with the Borel-Weil-Bott theorem and Serre duality.  If $\lambda$ is a regular weight and $w$ an
element of the Weyl group with $w\cdot\lambda=\mu\in\Lambda^{+}$ then we have

$$
\xymatrix{
\V_{\smu}
\ar@{=}[r]^(0.3){\mbox{\tiny BWB}} &
\left({\H^{\ell(w)}(\X,\L_{\lambda}) }\right)^{*}
\ar@{=}[r]^(0.45){\mbox{\tiny Serre}} &
\H^{\N-\ell(w)}(\X,\Ksh_{\X}\otimes\L_{\lambda}^{*}) 
\ar@{=}[r]^(0.55){\mbox{\tiny \ref{sec:serredual}}}
&
\H^{\ell(\wo w)}(\X,\L_{\Serre(\lambda)} )
\ar@{=}[r]^(0.65){\mbox{\tiny BWB }}_(0.65){\mbox{\tiny +
\eqref{eq:weightserreduality}}} &
\V_{-\wo\mu}^{*}.
}
$$

\bpoint{Schubert varieties} \label{sec:schubertvarieties}
For an element $w\in\W$ of the Weyl group the {\em Schubert variety} $\X_{w}$ is defined by

$$\X_{w} := \overline{\B w\B/\B}\subseteq \G/\B=\X.$$

\np
Recall that the classes of the Schubert cycles $\{[\X_{w}]\}_{w\in \W}$ give 
a basis for the cohomology ring $\H^{*}(\X,\ZZ)$ of $\X$.
Each $[\X_{w}]$ is a cycle of complex dimension $\ell(w)$.
The dual Schubert cycles $\{[\Omega_{w}]\}_{w\in \W}$, given by $\Omega_{w}:=\X_{\wo w}$, also form a basis. 
Each $[\Omega_{w}]$ is a cycle of complex codimension $\ell(w)$.
The work of Demazure \cite{dem1}, Kempf \cite{Ke}, Ramanathan \cite{R}, and Seshadri \cite{Sh} 
shows
that each Schubert variety $\X_{w}$ is normal with rational singularities.

\np
{\bf Remark.} If $w_1$,\ldots, $w_k$, and $w\in\W$ are such that 
$\ell(w)=\sum \ell(w_i)$, then the intersection $\cap_{i=1}^{k}[\Omega_{w_i}]\cdot[\X_{w}]$ is a number.
The number is the coefficient
of $[\Omega_{w}]$ when writing the product $\cap_{i=1}^{k}[\Omega_{w_i}]$ in terms of the basis 
$\{[\Omega_{v}]\}_{v\in \W}$.

\np
To reduce notation we use $w$ to also refer to the point $w\B/\B\in \X_{w} \subseteq\X$. 
In particular for the identity $e\in\W$,  $\X_{e}=\{e\}$.  
Note that $e\in\X$ is also the image of $1_{\G}$ under the projection from $\G$ onto $\X$.

\np
{\bf Bruhat order.}\label{sec:bruhat}
The {\em Bruhat order} on the Weyl group $\W$ is the partial order given by the relation
$v\leq w$ if and only if $\X_{v}\subseteq \X_{w}$.
The minimum element in this order is $e$ and the maximum element is $\wo$,
corresponding to the subvarieties $\X_{e}=\{e\}$ and $\X_{\wo}=\X$ respectively.

\np
The following result will be used several times throughout the paper. 

\tpoint{Lemma}\label{lem:nonzero}
Suppose that $w_1$,\ldots, $w_k$ are elements of the Weyl group 
such that $\delpos=\sqcup_{i=1}^{k} \invset_{w_i}$.  Then $\cap_{i=1}^{k}[\Omega_{w_i}]\neq 0$.  

\np
\bpf
Each class $[\Omega_{w_i}]$ is represented by any translation of the cycle $\Omega_{w_i}$, so to understand
$\cap_{i=1}^{k}[\Omega_{w_i}]$ we can study the intersection  of schemes

\begin{equation}\label{eqn:cycleintersection}
\bigcap_{i=1}^{k} (\wo w_i)^{-1}\Omega_{w_i}.
\end{equation}

\np
Each of the schemes $(\wo w_i)^{-1}\Omega_{w_i}$ passes through $e\in \X$.
The tangent space to $(\wo w)^{-1}\Omega_{w}$ at $e$ is 
$$
\Lie\left({(\wo w)^{-1}\B (\wo w)\rule{0cm}{0.5cm}}\right)/\Lie(\B) = 
\opp_{\alpha\in-\invset_{\wo w}} \ggg^{\alpha}\stackrel{\eqref{eqn:weightcomplement}}{=}
\opp_{\alpha\in-\invset_{w}^{\c}}\ggg^{\alpha} \subseteq 
\gb^{-} = \T_{e}\X,
$$

\np
where we have identified $\T_{e}\X$ with $\gb^{-}$ via the projection $\G\longrightarrow\X$.
Noting that

$$
\bigcap_{i=1}^k \Phi_w^{\c} = (\bigcup_{i=1}^k \Phi_w)^{\c} = (\Delta^+)^{\c} = \emptyset,
$$

\np 
we conclude that
the intersection of the tangent spaces of the varieties $(\wo w_i)^{-1}\Omega_{w_i}$ at $e\in \X$ is $0$.
Hence the intersection \eqref{eqn:cycleintersection} is transverse at the identity.    
By Kleiman's transversality theorem \cite[Corollary 4(ii)]{kl}, 
small translations of each of the varieties $(\wo w_i)^{-1}\Omega_{w_i}$ will intersect properly 
and compute the intersection number.  Small translations of varieties cannot remove transverse 
points of intersection and thus $\cap_{i=1}^{k}[\Omega_{w_i}]\neq 0$.
\epf

\bpoint{The Littlewood-Richardson Cone}
\label{sec:LR-cone}
For any $k\geq 1$, let $\LRC(k)$ be the {\em Littlewood-Richardson cone}, i.e., 
the rational cone generated by the tuples $(\mu_1,\ldots, \mu_k,\mu)$ of dominant weights
such that $\V_{\mu}$ is a component of $\V_{\mu_1}\otimes\cdots\otimes\V_{\mu_k}$.  
It is known that $\LRC(k)$ is polyhedral. 
A (nonminimal) set of equations describing $\LRC(k)$ is given by the conditions that
$(\mu_1,\ldots, \mu_k,\mu)\in\LRC(k)$ if and only if, 
for every $w_1$,\ldots, $w_k$, $w\in \W$ such that $\sum_{i=1}^{k}\ell(w_i)=\ell(w)$ and
$\cap_{i=1}^{k}[\Omega_{w_i}]\cdot \X_{w}\neq 0$, the weight $\sum_{i=1}^{k} w_i^{-1}\mu_i-w^{-1}\mu$ belongs to
$\Span_{\QQ_{\geq 0}} \Delta^+$ (see \cite[Theorem 5.2.1]{BeSj}). 

\np
For any dominant weights $\mu_1$,\ldots, $\mu_k$, let $\LRP(\mu_1,\ldots, \mu_k)$ be the slice of $\LRC(k)$ 
obtained by fixing the first $k$ coordinates to be $(\mu_1,\ldots, \mu_k)$.  
We call $\LRP(\mu_1,\ldots, \mu_k)$
the {\em Littlewood-Richardson polytope} of $\mu_1$,\ldots, $\mu_k$.

\np
{\bf Results of Ressayre.}
For any set $\I$ of simple roots, we define $\P_{\I}$ to be the parabolic subgroup associated to $\I$. 
For any parabolic $\P\supseteq\B$ we denote the Weyl group of $\P$ by $\W_{\P}$.

\np
For a set $\I$ of simple roots we wish to consider elements $w_1$,\ldots, $w_k$, and $w$ 
of $\W$ satisfying the following conditions with respect to $\I$:

\begin{equation}\label{eqn:Iconditions}
\left\{
\mbox{
\begin{minipage}{0.8\textwidth}
\begin{itemize}
\item[({\em i})] $\ell(w)=\sum_{i=1}^{k} \ell(w_i)$ and $\cap_{i=1}^{k}[\Omega_{w_i}]\cdot[\X_{w}]=1$.
\item[({\em ii})] Each $w_i$ is of minimal length in the coset $w_i\W_{\P_{\I}}$ and $w$ is of minimal length
in the coset $w\W_{\P_{\I}}$.
\item[({\em iii})] 
The weight $\sum_{i=1}^{k} w_i^{-1}\cdot 0 - w^{-1}\cdot 0$ belongs to $\Span_{\ZZ \geq 0} \I$.
\end{itemize}
\end{minipage}
}
\right.
\end{equation}

\np
The work of Ressayre gives an explicit description of some of the boundary components of $\LRC(k)$.
The following is a translation of \cite[Theorem D]{Re} into our notation:

\tpoint{Theorem} \label{thm:Ressayre}

\begin{enumerate}
\item
Let $\I$ be a set of simple roots and 
$w_1$,\ldots, $w_k$, and $w$ elements of $\W$ satisfying conditions \eqref{eqn:Iconditions} with respect to $\I$.
Then the set 

$$
\left\{(\mu_1,\ldots,\mu_k,\mu)\in\LRC(k) \, | \, \sum_{i=1}^{k} w_i^{-1}\mu_i-w^{-1}\mu \in \Span_{\QQ_\geq 0} \I \right\}
$$

\np
is a face of codimension $(n-\#\I)$ of $\LRC(k)$.

\medskip
\item Any face of $\LRC(k)$ which intersects the locus of strictly dominant weights is of the type in part ({\em a}).
\end{enumerate}

\bpoint{Symmetric and nonsymmetric forms}  \label{sec:nonsymmetric-symmetric}
Most questions we consider, including Problem I and Problem II, 
can be stated in nonsymmetric and symmetric
forms and it is frequently convenient to switch from one to the other.
We illustrate this procedure by showing how to switch from the nonsymmetric to the 
symmetric form of Problem I.

\np
In the nonsymmetric form 
we are given $w_1$,\ldots, $w_k$, and $w$, 
such that $\ell(w)=\sum \ell(w_i)$, 
and $\lambda_1$,\ldots, $\lambda_k$, and $\lambda$, such that $\lambda=\sum \lambda_i$, 
satisfying the additional conditions that $w_i\cdot\lambda_i\in \Lambda^{+}$ for $i=1$,\ldots, $k$, and 
$w\cdot\lambda\in\Lambda^{+}$.
This corresponds to the data of a cup product problem:

\begin{equation}\label{eq:nonsymmetric-sample}
\H^{\ell(w_1)}(\X,\L_{\lambda_1})\otimes\cdots\otimes\H^{\ell(w_k)}(\X,\L_{\lambda_k})
\stackrel{\cup}{\longrightarrow} 
\H^{\ell(w)}(\X,\L_{\lambda}).
\end{equation}

\np
Set $\mu_i=w_i\cdot\lambda_i$ for $i=1$,\ldots, $k$, and $\mu=w\cdot\lambda$ to keep track of the 
modules which appear as cohomology groups. By the Borel-Weil-Bott theorem the map 
\eqref{eq:nonsymmetric-sample} corresponds to a $\G$-equivariant map 

$$
\V_{\mu_1}^{*}
\otimes 
\V_{\mu_2}^{*}
\otimes\cdots\otimes\V_{\mu_k}^{*}\longrightarrow\V_{\mu}^{*}.
$$

\label{sec:serre-reduction}
\np
By Serre duality 
$\H^{\N-\ell(w)}(\X,\Ksh_{\X}\otimes\L_{\lambda}^{*})\neq 0$ and the cup product map
$$\H^{\ell(w)}(\X,\L_{\lambda})\otimes
\H^{\N-\ell(w)}(\X,\Ksh_{\X}\otimes\L_{\lambda}^{*})
\stackrel{\cup}{\longrightarrow}
\H^{\N}(\X,\Ksh_{\X})$$ 
is a perfect pairing.  Since $\H^{\ell(w)}(\X,\L_{\lambda})$ is an irreducible $\G$-module,
the surjectivity of \eqref{eq:nonsymmetric-sample} is equivalent to the surjectivity of 
the cup product map

\begin{equation}\label{eq:symmetric-sample}
\H^{\ell(w_1)}(\X,\L_{\lambda_1})\otimes\cdots\otimes \H^{\ell(w_k)}(\X,\L_{\lambda_k})
\otimes\H^{\N-\ell(w)}(\X,\Ksh_{\X}\otimes\L_{\lambda}^{*})
\stackrel{\cup}{\longrightarrow}
\H^{\N}(\X,\Ksh_{\X}).
\end{equation}

\np
To get the symmetric form of this problem, we set 
$w_{k+1}=\wo w$, $\lambda_{k+1}=\SS(\lambda)=-\lambda-2\rho$,
and $\mu_{k+1}=-\wo\mu=w_{k+1}\cdot\lambda_{k+1}$. 
Then 
$\L_{\lambda_{k+1}}=\Ksh_{\X}\otimes\L_{\lambda}^{*}$ by \S\ref{sec:serredual},
$w_{k+1}\cdot\lambda_{k+1}\in\Lambda^{+}$ by Lemma \ref{lem:serredual},
and $\ell(w_{k+1})=\N-\ell(w)$,
so that \eqref{eq:symmetric-sample} becomes 

$$
\H^{\ell(w_1)}(\X,\L_{\lambda_1})\otimes\cdots\otimes \H^{\ell(w_k)}(\X,\L_{\lambda_k})
\otimes\H^{\ell(w_{k+1})}(\X,\L_{\lambda_{k+1}})
\stackrel{\cup}{\longrightarrow}
\H^{\N}(\X,\Ksh_{\X}).
$$

\np
Since $\sum_{i=1}^{k+1}\lambda_i = \lambda + (-\lambda-2\rho)=-2\rho$ and $\L_{-2\rho}=\Ksh_{\X}$, this is again
a cup product problem of the type we consider, but now all weights $\lambda_1$, \ldots, $\lambda_{k+1}$ and 
Weyl group elements $w_1$,\ldots, $w_{k+1}$ play equal roles.

\np
By \eqref{eqn:weightcomplement} $\invset_{w_{k+1}}=\invset_{w}^{\c}$ 
and therefore the condition that 
$\invset_{w}=\sqcup_{i=1}^{k} \invset_{w_i}$ is equivalent to the condition $\delpos=\sqcup_{i=1}^{k+1}\invset_{w_i}$.
Since $[\Omega_{w_{k+1}}] = [\X_{\wo w_{k+1}}] = [\X_{w}]$, the intersection numbers
$\bigcap_{i=1}^{k} [\Omega_{w_i}]\cdot [\X_{w}]$  and $\bigcap_{i=1}^{k+1} [\Omega_{w_i}]$ are the same.
Finally,
the multiplicity of $\V_{\mu}$ in
$\V_{\mu_1}\otimes\cdots\otimes\V_{\mu_k}$ is the same as the multiplicity of the trivial module in 
$\V_{\mu_1}\otimes\cdots\otimes\V_{\mu_k}\otimes\V_{\mu_{k+1}}$
because $\V_{\mu_{k+1}}=\V_{-\wo\mu}=\V_{\mu}^{*}$.

\np
To go from the symmetric form to the nonsymmetric form we simply reverse the above procedure, although of course
we are free to desymmetrize with respect to any of the indices $i=1$,\ldots, $k+1$, and not just the last one.

\np
For convenience we list below the symmetric and nonsymmetric forms of some formulas and expressions 
we are interested in.

\medskip

{
\renewcommand{\arraystretch}{2.7}

\begin{centering}
\begin{longtable}{|c|c|}
\hline
Nonsymmetric & Symmetric \\
\hline
$\displaystyle \ott_{i=1}^{k} \H^{\ell(w_i)}(\X,\L_{\lambda_i}) \longrightarrow \H^{\ell(w)}(\X,\L_{\lambda})$
&
$\displaystyle \ott_{i=1}^{k+1} \H^{\ell(w_i)}(\X,\L_{\lambda_i}) \longrightarrow \H^{\N}(\X,\K_{\X})$
\\
\hline
$\displaystyle \sum_{i=1}^{k} \ell(w_i) = \ell(w) $ & $\displaystyle \sum_{i=1}^{k+1} \ell(w_i) = \N$ \\
\hline
$\displaystyle \sum_{i=1}^{k} \lambda_i = \lambda$ & $\displaystyle \sum_{i=1}^{k+1} \lambda_i = -2\rho $ \\
\hline
$\displaystyle \sum_{i=1}^{k} w_i^{-1}\mu_i - w^{-1}\mu$ & $\displaystyle \sum_{i=1}^{k+1} w_i^{-1}\mu_i$ \\
\hline
$\displaystyle \sum_{i=1}^{k} w_i^{-1}\cdot 0 - w^{-1}\cdot 0$ & $\displaystyle \sum_{i=1}^{k+1} w_i^{-1}\cdot0 + 2\rho$ \\
\hline
$\displaystyle \invset_{w} = \bigsqcup_{i=1}^{k} \invset_{w_i}$ & $\displaystyle \delpos = \bigsqcup_{i=1}^{k+1} \invset_{w_i}$ \\
\hline
$\displaystyle \bigcap_{i=1}^{k} [\Omega_{w_i}]\cdot [\X_{w}]$ & $\displaystyle \bigcap_{i=1}^{k+1} [\Omega_{w_i}]$ \\
\hline
\end{longtable} 
\end{centering}
\renewcommand{\arraystretch}{1}
}

\medskip
\np
Since $k$ is an arbitrary positive integer, after switching to the symmetric form we often use $k$ in place of $k+1$
to reduce notation.

\bpoint{Demazure reflections} \label{sec:demazure}
Suppose that $\WW$ and $\M$  are varieties and $\pi\colon\WW\longrightarrow\M$ is a $\PP^1$-fibration, 
i.e., a smooth morphism with fibres isomorphic to $\PP^1$.  Let $\L$ be a line bundle on $\WW$
and $b$ be the degree of $\L$ on the fibres of $\pi$.  
Demazure \cite[Theorem 1]{dem2} proves the following isomorphism of vector bundles on $\M$:

\begin{equation}\label{eq:demazure-isom}
\R^{i}\pi_{*}\L \cong \R^{1-i}\pi_{*}\left({\L\otimes\omega_{\pi}^{b+1}}\right)\quad\quad\mbox{for $i=0,1$.}
\end{equation}

\np
where $\omega_{\pi}$ is the relative cotangent bundle of $\pi$.
The line bundle $\L\otimes\omega_{\pi}^{b+1}$ is called the {\em Demazure reflection of $\L$ with respect to
$\pi$}.

\np
Note that there is at most one value of $i$ for which the resulting vector bundles are nonzero:  $i=0$ if $b\geq 0$, 
$i=1$ if $b\leq -2$, and neither if $b=-1$.  Equation \eqref{eq:demazure-isom} and the corresponding
Leray spectral sequence
give the isomorphisms 
$$
\H^{j}(\WW,\L) \cong \left\{
\begin{array}{cl}
\H^{j+1}(\WW,\L\otimes\omega_{\pi}^{b+1}) & \mbox{if $b\geq 0$} \\
\H^{j-1}(\WW,\L\otimes\omega_{\pi}^{b+1}) & \mbox{if $b\leq -2$} \\
\end{array}\right.
\mbox{\,\,for all $j$.}
$$

\np
{\bf Link between Demazure reflections and the affine action.}
Let $\alpha_i$ be any simple root, $\P_{\!\alpha_i}$ the parabolic associated to $\alpha_i$, 
and $\pi_i\colon\X\longrightarrow \M_i:=\G/\P_{\!\alpha_i}$ the corresponding $\PP^1$-fibration. 
The relative cotangent bundle $\omega_{\pi_i}$ of $\pi_i$ is the line bundle $\L_{-\alpha_i}$.
Given any $\lambda\in\Lambda$, the degree of 
the line bundle $\L_{\lambda}$ on the fibres of $\pi_i$ is $\lambda(\alpha_i^{\vee})$ where 
$\alpha_i^{\vee}$ is the coroot corresponding to $\alpha_i$. 
We thus obtain that the Demazure reflection of $\L_{\lambda}$ with respect to the fibration $\pi_i$ is 
the line bundle:
$$\L_{\lambda}\otimes \L_{-\alpha_i}^{\lambda(\alpha_i^{\vee})+1} = 
\L_{\lambda-(\lambda(\alpha_i^{\vee})+1)\alpha_i}=
\L_{s_i\lambda-\alpha_i}=
\L_{s_i\cdot\lambda}
$$

\np
where $s_i$ is the simple reflection corresponding to $\alpha_i$.   
The combinatorics of performing 
Demazure reflections with respect to the various $\PP^1$-fibrations of $\X$ is therefore kept track of by
the affine action of the Weyl group on $\Lambda$.  In particular, if $v=s_{i_1}\cdots s_{i_m}\in\W$
and $\lambda\in\Lambda$ the result of applying the Demazure reflections 
with respect to the fibrations $\pi_{i_m}$, $\pi_{i_{m-1}}$, \ldots, $\pi_{i_1}$ in that order to $\L_{\lambda}$
is $\L_{v\cdot\lambda}$.

\np
{\bf Demazure reflections and base change.} Given any morphism $h\colon\Y_2\longrightarrow\M$ we can form the fibre product diagram
$$
\begin{array}{c}
\xymatrix{
\Y_1\ar[r]^{f}\ar[d]^{\pi_1} & \WW\ar[d]^{\pi} \\
\Y_2\ar[r]^{h}\ar@{}[ur]|{\Box} & \M \\
}
\end{array}.
$$

\np
If $\pi$ is a $\PP^1$-fibration then so is $\pi_1$, and $\omega_{\pi_1}=f^{*}\omega_{\pi}$.  Therefore, for any
line bundle $\L$ on $\V$, 
we have
$$f^{*}(\L\otimes\omega_{\pi}^{b+1}) = (f^{*}\L)\otimes\omega_{\pi_1}^{b+1}$$
where $b$ is the degree of $\L$ on the fibres of $\pi$.
The degree of $f^{*}\L$ on $\pi_1$ is also $b$ and therefore the formula above shows that 
the pullback of the Demazure reflection of $\L$ with respect to $\pi$  is the 
Demazure reflection of the pullback of $\L$ with respect to $\pi_1$.
Furthermore, by the theorem on cohomology and base change, the natural morphisms

$$
\R^{i}\pi_{1*}\left({f^{*}\L}\right) \stackrel{\sim}{\longleftarrow}h^{*}\left({\R^{i}\pi_{*}\L}\right)
\,\,\mbox{and}\,\,
\R^{1-i}\pi_{1*}\left(({f^{*}\L)\otimes\omega_{\pi_1}^{b+1}}\right) \stackrel{\sim}{\longleftarrow}
h^{*}\left({\R^{1-i}\pi_{*}(\L\otimes\omega_{\pi}^{b+1})}\right)
$$

\np
are isomorphisms for $i=0,1$.

\bpoint{$\E_2$-terms and computation of maps on cohomology} \label{sec:E2compute-begin}
Suppose that we have a commutative diagram of varieties

$$
\xymatrix{
\WW'\,\ar@{^{(}->}[r]^{\gamma}\ar[d]^{\pi'} & \WW\ar[d]^{\pi} \\
\M'\,\ar@{^{(}->}[r]& \M \\
}
$$

\np
where the vertical maps are proper and the horizontal maps are closed immersions.  
Suppose further that we have 
coherent sheaves $\Fsh$ on $\WW$ and $\Fsh'$ on $\WW'$, 
and a map $\varphi\colon\gamma^{*}\Fsh\longrightarrow \Fsh'$ of sheaves on $\WW'$.   The map $\varphi$ induces
maps $\varphi_{d}\colon\H^{d}(\WW,\Fsh)\longrightarrow\H^{d}(\WW',\Fsh')$ on cohomology and maps 
$\varphi_{d,k}\colon \H^{d-k}(\M,\R^{k}_{\pi*}\Fsh)\longrightarrow \H^{d-k}(\M',\R^{k}_{\pi*}\Fsh')$
on the $\E_2$-terms of the Leray spectral sequences for $\Fsh$ and $\Fsh'$ with respect to $\pi$ and $\pi'$.
Assume that both spectral sequences degenerate at the $\E_2$-term.
In \S\ref{sec:ThmIcor} we will need to know when we can compute $\varphi_{d}$ by knowing the maps $\varphi_{d,k}$.

\np
By the definition of convergence of a spectral sequence there are increasing filtrations

$$
0=\U_{-1}\subseteq
\U_{0} \subseteq
\cdots
\subseteq
\U_{d} 
=\H^{d}(\WW,\Fsh)
\,\,\mbox{and}\,\,
0=\U_{-1}'\subseteq
\U_{0}' \subseteq
\cdots
\subseteq
\U_{d}'
=\H^{d}(\WW',\Fsh')
$$

\np
such that $\U_{k}/\U_{k-1}=\H^{d-k}(\M,\R^{k}_{\pi*}\Fsh)$ and 
$\U_{k}'/\U_{k-1}'=\H^{d-k}(\M',\R^{k}_{\pi*}\Fsh')$ for $k=0$,\ldots, $d$.
Since the map $\varphi_{d}$ on the cohomology groups is compatible with the filtrations (in the sense that 
$\varphi_{d}(\U_{k})\subseteq \U_{k}'$ for $k=-1$,\ldots, $d$), $\varphi_d$ induces maps between the associated graded
pieces of the filtrations; these maps are exactly the maps $\varphi_{d,k}$.

\np
We will need to know that $\varphi_{d}$ can be computed from the maps $\varphi_{d,k}$ 
in an elementary case. Suppose there is a unique $k$ such that $\U_{k}/\U_{k-1}$ is nonzero (and so 
$\U_{k}/\U_{k-1}=\H^{d}(\WW,\Fsh)$), and a unique $k'$ such that
$\U_{k'}'/\U_{k'-1}'$ is nonzero (and so $\U_{k'}'/\U_{k-1}'=\H^{d}(\WW',\Fsh')$).  
Then we can compute $\varphi_{d}$ from the maps $\varphi_{d,k}$ if and only if $k=k'$;  
if this occurs then $\varphi_{d}=\varphi_{d,k}$.  

\np
In order to show that we must check the condition $k=k'$ above, 
i.e., that the map on $\E_2$-terms does not always determine the map $\varphi_d$, 
we give the following example of a nonzero map between cohomology groups of sheaves 
where the induced map on $\E_2$-terms is zero.  This example is also a cup product map.

\tpoint{Example} Let $\WW=\PP^{m}\times\PP^{m}$ for some $m\geq 1$, $\Fsh=\Osh_{\PP^{m}}(1)\bt\Osh_{\PP^{m}}(-r)$ with 
$r\geq m+2$, and let $\Gsh=\Osh_{\Delta}(1-r)$ be the restriction of $\Fsh$ to the diagonal of $\WW$.  We have
$\H^{m}(\WW,\Fsh)=\H^{0}(\PP^{m},\Osh_{\PP^m}(1))\otimes\H^{m}(\PP^{m},\Osh_{\PP^m}(-r))$ and 
$\H^{m}(\WW,\Gsh)=\H^{m}(\PP^{m},\Osh_{\PP^{m}}(1-r))$.  The natural restriction map 
$\varphi\colon\Fsh\longrightarrow\Gsh$
induces the cup product map

$$\varphi_{m}\colon\H^{0}(\PP^{m},\Osh_{\PP^m}(1))\otimes\H^{m}(\PP^{m},\Osh_{\PP^m}(-r))\stackrel{\cup}{\longrightarrow}
\H^{m}(\PP^{m},\Osh_{\PP^{m}}(1-r))$$

\np
which is a surjective map of nonzero groups.

\np
If $\pi\colon\WW\longrightarrow\M=\PP^{m}$ is the projection onto the first factor then both of 
the Leray spectral sequences degenerate at the $\E_2$ term with only one nonzero entry in each sequence.  We have
$\H^{m}(\WW,\Fsh)=\H^{0}(\M,\R^{m}_{\pi*}\Fsh)$ (i.e., $k=m$) and 
$\H^{m}(\WW,\Gsh)=\H^{m}(\M,{\pi_*}\Gsh)$ (i.e., $k'=0$).  
The maps $\varphi_{m,k}$ on the $\E_2$-terms are clearly zero, 
even though $\varphi_{m}$ is nonzero.

\bpoint{Bott-Samelson-Demazure-Hansen Varieties}
\label{sec:Zmaxpoint}
\label{sec:BSDH}
Let $\vv=s_{i_1}\cdots s_{i_m}$ be a word, not necessarily reduced, of simple reflections.  
Associated to $\vv$ is a variety $\Z_{\vv}$, a left action of $\B$ on $\Z_{\vv}$, and a
$\B$-equivariant map $f_{\vv}\colon\Z_{\vv}\longrightarrow \X$.  
If $\vv$ is nonempty there is also a
$\B$-equivariant map $\pi_{\vv}\colon \Z_{\vv}\longrightarrow \Z_{\vv_R}$ expressing $\Z_{\vv}$ as a $\PP^1$-bundle 
over $\Z_{\vv_R}$ together with a $\B$-equivariant $\sigma_{\vv}\colon \Z_{\vv_R}\longrightarrow \Z_{\vv}$ section
such that $f_{\vv_R}=f_{\vv}\circ \sigma_{\vv}$.

\np
These varieties were originally constructed by Demazure \cite{dem1} and Hansen \cite{han} following an analogous
construction by Bott and Samelson \cite{botsam} in the compact case.
In this subsection we recall their construction and several related facts.  We give two different descriptions of
the construction; both will be used in the constructions in section \ref{sec:Ys}.

\np
{\bf Recursive Construction.}\label{firstconst}
Recall that $e$ is unique point of $\X$ fixed by $\B$. 
If the word $\vv$ is empty we define $\Z_{\vv}$ to be $e$, the map $f_{\vv}$ to be the inclusion
$e\hookrightarrow \X$, and the $\B$-action on $\Z_{\vv}$ to be trivial.

\np
If $\vv=s_{i_1}\cdots s_{i_m}$ is nonempty, let $\uu=\vv_{R}=s_{i_1}\cdots s_{i_{m-1}}$ be the word 
obtained by dropping the rightmost reflection of $\vv$.  By induction we have already constructed $\Z_{\uu}$ 
and the map $f_{\uu}\colon\Z_{\uu}\longrightarrow \X$.
Set $h=\pi_{i_m}\circ f_{\uu}$, where 
$\pi_{i_m}$ is the $\G$-equivariant projection (and $\PP^1$-fibration) 
$\X\longrightarrow \M_{i_m}=\G/\P_{\alpha_{i_m}}$.
We then define $\Z_{\vv}$ to be the fibre product $\Z_{\uu} \times_{\M_{i_m}}\X$, and $f_{\vv}$ and $\pi_{\vv}$ 
to be the maps from the fibre product to $\X$ and to $\Z_{\uu}$ respectively.  
Since $h=\pi_{i_m}\circ f_{\uu}$
by the universal property of the fibre product there exists a unique map 
$\sigma_{\vv}\colon\Z_{\uu}\longrightarrow \Z_{\vv}$ such that  
$f_{\uu}=f_{\vv}\circ\sigma_{\vv}$ and $\id_{\Z_{\uu}}=\pi_{\vv}\circ\sigma_{\vv}$.
These maps are summarized in the following diagram, where the square is a fibre product:

\begin{equation}\label{eqn:fibredef}
\begin{array}{c}
\xymatrix{
\Z_{\vv}\ar[rr]^{f_{\vv}}\ar[dd]^{\pi_{\vv}} & & \X\ar[dd]^{\pi_{i_m}} \\
& \Box & \\
\Z_{\uu}\ar[rr]^{h}\ar'[ur][uurr]^(.2){f_{\uu}} \ar@/^20pt/[uu]^{\sigma_{\vv}} & & 
\M_{\subsmash{i_m}}  \\
}
\end{array}.
\end{equation}

\np
Since $\B$ acts on $\Z_{\uu}$ and on $\X$, and the maps $f_{\uu}$, $\pi_{i_m}$, and $h$ are $\B$-equivariant, 
by the universal property of the fibre product, the diagram \eqref{eqn:fibredef}
induces  a $\B$-action on $\Z_{\vv}$ such that $f_{\vv}$ and $\sigma_{\vv}$ are $\B$-equivariant maps.
Since each morphism $\sigma_{\vv}$ is a $\PP^1$-fibration it follows immediately that each $\Z_{\vv}$ is a smooth
proper variety of dimension $\ell(\vv)$.

\np
{\bf Direct Construction. } \label{secondconst}
For any word $\vv$ set

$$\P_{\vv}:= \left\{{
\begin{array}{cl}
e & \mbox{if $\vv$ is empty} \\
\P_{\alpha_{i_1}}\times\cdots\times\P_{\alpha_{i_m}} & \mbox{if $\vv=s_{i_1}\cdots s_{i_m}$ is nonempty.} \\
\end{array}
}\right.
$$

\np
If $\vv$ is empty we define $\Z_{\vv}$, $f_{\vv}$, and the $\B$-action as in the direct construction.

\np
If $\vv=s_{i_1}\cdots s_{i_m}$ is nonempty then 
$\Z_{\vv}$ is the quotient of $\P_{\vv}$ by $\B^{m}$, where an element
$(b_1,\ldots, b_m)$ of $\B^m$ acts on the right on $(p_1,\ldots, p_m)$ by 

$$(p_1,\, \ldots,\,  p_m)\cdot(b_1,\, \ldots,\,  b_m) = 
(p_1b_1,\, b_1^{-1}p_2b_2,\, b_2^{-1}p_3b_3,\,\ldots,\, b_{m-1}^{-1}p_mb_m).$$

\np
The left action of $\B$ on $\P_{\vv}$ given by 
$$b\cdot(p_1,\,p_2,\,\cdots,\,p_m)=(bp_1,\, p_2,\,\cdots,\,p_m)$$
commutes with the right action of $\B^{m}$ and therefore descends to a left action of $\B$ on $\Z_{\vv}$.
We denote the corresponding $\B$-equivariant quotient map by $\psi_{\vv}\colon \P_{\vv}\longrightarrow\Z_{\vv}$.

\np
The product map $\P_{\vv}\stackrel{\phi_{\vv}}{\longrightarrow} \G$ given by 
$(p_1,\ldots, p_m)\mapsto p_1\cdots p_m$ is equivariant for the left $\B$-action described above and left 
multiplication of $\G$ by $\B$. Under the homomorphism of groups
$\B^m\longrightarrow \B$ given by the projection $(b_1,\ldots,b_m)\mapsto b_m$ the product 
map $\phi_{\vv}$ is also equivariant for the right action
of $\B^m$ on $\P_{\vv}$ and the right multiplication of $\G$ by  $\B$.  The product map therefore
descends to a left $\B$-equivariant morphism $f_{\vv}\colon \Z_{\vv}\longrightarrow \X$.

\np
Let 
$\uu=\vv_{R}=s_{i_1}\cdots s_{i_{m-1}}$ the word obtained by dropping the rightmost reflection in $\vv$.
The projection map $\pr_{\vv}\colon\P_{\vv}{\longrightarrow} \P_{\uu}$ 
sending $(p_1,\ldots, p_m)$ to $(p_1,\ldots, p_{m-1})$ 
is equivariant with respect to the projection $\B^{m}\longrightarrow \B^{m-1}$
sending $(b_1,\ldots, b_m)$ to $(b_1,\ldots, b_{m-1})$. 
Similarly 
the inclusion map $j_{\vv}\colon\P_{\uu}{\hookrightarrow} \P_{\vv}$ 
sending $(p_1,\ldots, p_{m-1})$ to $(p_1,\ldots, p_{m-1},1_{\G})$
is equivariant with respect to the inclusion $\B^{m-1}\hookrightarrow\B^{m}$
sending $(b_1,\ldots, b_{m-1})$ to $(b_1,\ldots, b_{m-1},b_{m-1})$.  
The maps $\pr_{\vv}$ and $j_{\vv}$ respect the left $\B$-action
on $\P_{\vv}$ and $\P_{\uu}$, and
therefore descend to $\B$-equivariant maps $\pi_{\vv}\colon\Z_{\vv}\longrightarrow\Z_{\uu}$ and
$\sigma_{\vv}\colon\Z_{\uu}\longrightarrow\Z_{\vv}$.  
Since $\pr_{\vv}\circ j_{\vv}=\id_{\P_{\uu}}$ and
$\phi_{\vv}\circ j_{\vv}=\phi_{\uu}$, taking quotients we obtain
$\pi_{\vv}\circ\sigma_{\vv}=\id_{\Z_{\uu}}$ and $f_{\vv}\circ\sigma_{\vv}=f_{\uu}$.
Finally, the fibres of $\pi_{\vv}$ are isomorphic to $\P_{\alpha_{i_m}}/\B\cong \PP^1$.

\np
We record the following well-known facts about the construction above.

\tpoint{Proposition} \label{prop:ZvandXv}

\begin{enumerate}
\item The varieties $\Z_{\vv}$ produced by the recursive and direct constructions above are isomorphic over $\X$.
\item If $\vv=s_{i_1}\cdots s_{i_m}$ is a reduced word with product $v$ 
then the image of $f_{\vv}\colon\Z_{\vv}\longrightarrow \X$ is $\X_{v}$ and $f_{\vv}$ 
is a resolution of singularities of $\X_{v}$.

\end{enumerate}

\bpf 
Part ({\em b}) is proved in \cite{dem1} and \cite{han}.
To show ({\em a}) it is enough to show that the varieties produced by the direct construction 
satisfy the fibre product diagram \eqref{eqn:fibredef}. This is most easily checked after pulling 
back \eqref{eqn:fibredef} via the maps $\G\longrightarrow\X$ 
and $\P_{\uu}=\P_{i_1}\times\cdots\times \P_{i_{m-1}}\longrightarrow\Z_{\uu}\,$; the details are omitted here.
\epf

\np
{\bf Maximum points.}  
Let $\vv=s_{i_1}\cdots s_{i_m}$ be a reduced word with product $v$.
By 
Proposition \ref{prop:ZvandXv}({\em b}) the image of $\Z_{\vv}$ under $f_{\vv}$ is $\X_{v}$ and one can check
that there is a unique point $p_{\vv}$ of $\Z_{\vv}$ which maps to ${v}\in\X_{v}$.  
More specifically, from the point of view of the direct construction the point 
$(s_{i_1},\ldots,s_{i_m})$ is a point of $\P_{\vv}$ and its image under the quotient map 
$\P_{\vv}\longrightarrow\Z_{\vv}$
is $p_{\vv}$.  From the point of view of the recursive construction one starts with $p_{\emptyset}=e$, 
and recursively defines $p_{\vv}$ to be unique torus fixed point in the $\PP^1$-fibre of 
$\pi_{\vv}\colon\Z_{\vv}\longrightarrow \Z_{\uu}$ over $p_{\uu}$ which is not equal to $\sigma_{\vv}(p_{\uu})$, 
where $\uu=\vv_{R}=s_{i_1}\cdots s_{i_{m-1}}$.
Note that $p_{\vv}$ is the unique torus fixed point of $\Z_{\vv}$ 
whose image in $\X_{v}$ is the largest in the Bruhat order among torus-fixed points of $\X_{v}$. 
We call $p_{\vv}$ the {\em maximum point} of $\Z_{\vv}$.

\np
Since $p_{\vv}$ is a torus fixed point, the torus acts on the tangent space $\T_{p_{\vv}}\Z_{\vv}$ and it will be 
important for us to know the formal character of $\T_{p_{\vv}}\Z_{\vv}$.  
It follows inductively from the recursive construction that

\begin{equation}\label{eqn:Z-tgtweights}
\Ch(\T_{p_{\vv}}\Z_{\vv})  = \langle \invset_{v^{-1}} \rangle. 
\end{equation}

\bpoint{Semi-stability of torus fixed points}\label{sec:PRVpoints} 
The following lemma is due to Kostant.

\np
\tpoint{Lemma}\label{lem:PRV-point}
Let $\WW$ be a projective variety with a $\G$-action and $\L$ a $\G$-equivariant ample line bundle 
on $\WW$.  A torus fixed point $q\in\WW$ is semi-stable with respect to $\L$ if and only if the weight of 
 $\L_{q}$ is zero.  In this case the orbit of $q$ is closed in the semi-stable locus.

\np
\bpf
If the action of the torus on the fibre $\L_{q}$ is non-trivial then it is easy to see (for instance using the
Hilbert-Mumford criterion for semi-stability, \cite[Theorem 2.1, p.\ 49]{git}) that $q$ is not a stable point.

\np
Conversely, suppose that the weight of $\L_{q}$ is zero.  Replacing $\L$ by a multiple we may assume
that $\L$ is very ample and gives an embedding $\WW\hookrightarrow \PP^{r}$ for some $r$.  Let 
$\AA^{\!r+1}$ be the affine space corresponding to $\PP^{r}$ and 
$\AA^{\!r+1}\setminus\{0\}\longrightarrow\PP^{r}$ be the quotient map.
Then $\G$ acts linearly on $\AA^{\!r+1}$ inducing an action on $\PP^{r}$ compatible with the action on $\WW$.  Let
$\tilde{q}$ be any lift to $\AA^{\!r+1}$ of the image of $q$ in $\PP^{r}$.  The condition that the torus
act trivially on $\L_{q}$ is equivalent to the condition that $\tilde{q}$ be fixed by $\T$ under
the $\G$-action on $\AA^{\!r+1}$.
Kostant (\cite[p. 354, Remark 11]{kos}) proves that for any finite dimensional 
module of a reductive group $\G$ and any point $\tilde{q}$ fixed by $\T$, the $\G$-orbit 
of $\tilde{q}$ is closed; this result was also later generalized by Luna \cite[Theorem (**)]{lu2}.
Since $\G$ is reductive and the orbit of $\tilde{q}$ does not meet zero, 
there is a $\G$-invariant homogeneous form of some degree $m$ 
which is nonzero on $\tilde{q}$.  This corresponds to a $\G$-invariant section $s\in\H^{0}(\WW,\L^{m})^{\G}$ 
such that $s(q)\neq 0$.  We thus see that if the weight of  $\L_{q}$ is zero then $q$ is a semi-stable
point, and the orbit of $q$ is closed in the semi-stable locus.   \epf

\section{Diagonal Bott-Samelson-Demazure-Hansen-Kumar Varieties}\label{sec:Ys}

\np
In this section
we give a generalization of the varieties from \S\ref{sec:BSDH}.  
The construction is a variation of a construction of Kumar in \cite{ku2};
see \S\ref{sec:compareKumar} for a comparison.
These varieties are obtained by applying the idea of the Bott-Samelson resolution to the diagonal inclusion 
$\X\hookrightarrow \X^{k}$.  They can also be thought of as a desingularization of the total space of the variety
of intersections of translates of Schubert cycles.    This alternate description is established in Theorem 
\ref{thm:fibres}. 

\np
More specifically,
for each sequence $\vseq=(\vv_1,\ldots, \vv_k)$ of words we construct a smooth variety
$\Y_{\svseq}$ of dimension $\N+\ell(\vseq)$ with a $\G$-action together with a proper 
map $f_{\vseq}\colon\Y_{\svseq}\longrightarrow \X^{k}$ which is 
$\G$-equivariant for the diagonal action of $\G$ on $\X^{k}$.
If $\useq$ is the sequence obtained by dropping a single simple reflection from the right of one of the $\vv_j$'s
then $\Y_{\svseq}$ is a $\PP^1$-fibration over $\Y_{\suseq}$, and there is a section 
$\Y_{\suseq}\hookrightarrow \Y_{\svseq}$ compatible with the maps $f_{\vseq}$ and $f_{\useq}$ to $\X^{k}$.
The fibration and section maps are $\G$-equivariant; moreover 
they are compatible with the $\PP^1$-fibrations on factors of $\X^{k}$. 
These relationships are summarized in diagram \eqref{eqn:bigfibredef}.

\bpoint{Recursive Construction} \label{sec:1constY}
Let $\vseq=(\vv_1,\ldots, \vv_k)$ be a sequence of words.  If
all $\vv_j$ are empty, i.e., if $\vseq=(\emptyset,\ldots, \emptyset)$, we set $\Y_{\svseq}=\X$ and
let $f_{\vseq}\colon \Y_{\svseq}\longrightarrow \X^{k}$ be the diagonal embedding.

\np
Otherwise suppose that $\vv_j$ is nonempty. Let 

\begin{equation}\label{eqn:uudef}
\uu_l :=
\left\{{
\begin{array}{lc}
\hphantom{(}\vv_l & \mbox{if $l\neq j$} \\
(\vv_j)_R & \mbox {if $l=j$} \rule{0cm}{0.7cm}\\
\end{array}
}\right.,\quad\quad\mbox{$l=1,\ldots,k$}
\end{equation}

\np
and set 
$\useq=(\uu_1,\ldots,\uu_{k})$.
By induction on $\ell(\vseq)$ we may assume that $\Y_{\suseq}$ and the map
$f_{\useq}\colon\Y_{\useq}\longrightarrow\X^{k}$ have been constructed.
If $\vv_j=s_{i_1}\cdots s_{i_m}$, so that $\uu_j=s_{i_1}\cdots s_{i_{m-1}}$
then
we define $\Y_{\svseq}$, the map $f_{\vseq}$, the projection $\pi_{\vseq,\useq}$, and the section 
$\sigma_{\vseq,\useq}$ by the following fibre product square:

\begin{equation}\label{eqn:bigfibredef}
\begin{array}{c}
\xymatrix@C=20pt{
\Y_{\svseq}\ar[rr]^{f_{\vseq}}\ar[dd]^{\pi_{\vseq,\useq}} & & 
\X^{\subsmash{k}}
\ar[dd]^{(\id_{\X})^{j-1}\times\pi_{i_m}\times(\id_{\X})^{k-j}} 
\\
& & \\
\Y_{\suseq}\ar[rr]\ar[uurr]|*+<4pt>{\rule{0.3cm}{0cm}\mbox{\large$\,\Box\,$}\rule{0.3cm}{0cm}}^(0.55){\rule[-0.3cm]{0.3cm}{0cm}f_{\useq}} 
\ar@/^20pt/[uu]^{\sigma_{\vseq,\useq}} & & 
\X^{j-1}\times\M_{i_m}\times\X^{\subsmash{k\!-\!j}}  \\
}
\end{array}.
\end{equation}

\np
Here $\pi_{i_m}\colon \X\longrightarrow \M_{i_m}:=\G/\P_{\alpha_{i_m}}$ is the natural projection, and
$\X^{k}\longrightarrow \X^{j-1}\times\M_{i_m}\times \X^{k-j}$ is the projection $\pi_{i_m}$ on the $j$-th factor
and the identity on all others.  The bottom map $\Y_{\suseq}\longrightarrow \X^{j-1}\times\M_{i_m}\times\X^{k-j}$ is
the map $f_{\useq}$ to $\X^{k}$ followed by the map 
$\X^{k}\longrightarrow \X^{j-1}\times\M_{i_m}\times \X^{k-j}$ above.

\np
Since $\X^{k}\longrightarrow \X^{j-1}\times \pi_{i_m}\times \X^{k-j}$ is a $\PP^1$-fibration the same is true of
$\pi_{\vseq,\useq}$.  We conclude by induction that the variety $\Y_{\svseq}$ is smooth, proper, and irreducible of dimension
$\N+\ell(\vseq)$.  
The maps $f_{\useq}$ and $\id_{\Y_{\suseq}}$ from $\Y_{\suseq}$ to $\X^{k}$ and $\Y_{\suseq}$ respectively
give rise to the section $\sigma_{\vseq,\useq}$.  
By construction we have $f_{\useq}=f_{\vseq}\circ\sigma_{\vseq,\useq}$ and 
$\id_{\Y_{\useq}}=\pi_{\vseq,\useq}\circ\sigma_{\vseq,\useq}$.

\np
This construction is well-defined.  Indeed, assume that we had dropped a simple reflection from the 
right of $\vv_{j'}$, $j'\neq j$ to obtain a sequence of words $\useq'$ and used $\Y_{\suseq'}$ instead
of $\Y_{\suseq}$ to construct $\Y_{\svseq}$.  We claim that the resulting variety $\Y_{\svseq}$ is the same.
This follows easily by induction on $\ell(\vseq)$ 
and the fact that the diagram expressing the commutativity of the projections on the different factors is
a fibre square:

$$
\begin{array}{c}
\xymatrix{
\X^{k} \ar[r]^{(\id_{\X})^{j'-1}\times\pi_{i_m}\times(\id_{\X})^{k-j'}}   
\ar[dd]_{(\id_{\X})^{j-1}\times\pi_{i_m}\times(\id_{\X})^{k-j}}   &
\X^{j'-1}\times\M_{i_m'}\times\X^{k-j'} \ar[dd] \\
\\
\X^{j-1}\times \M_{i_m} \times \X^{k-j} \ar[r] \ar@{}[uur]|{\mbox{$\Box$}}
& \X^{j'-1}\times\M_{i_{m'}}\times\X^{j-j'-1}\times
\M_{i_m}\times \X^{k-j} 
}
\end{array}.
$$

\np
Here, by symmetry, we have assumed that $j'<j$.

\bpoint{Direct Construction}\label{sec:2constY}\label{sec:2constYmap}
Let $\vseq=(\vv_1,\ldots, \vv_k)$ be a sequence of words.  The group
$\B$ acts diagonally on $\Z_{\vv_1}\times\cdots\times\Z_{\vv_k}$ on the left.   
We define $\Y_{\svseq}$
to be the quotient of $\G\times(\Z_{\vv_1}\times\cdots\times\Z_{\vv_k})$ by the left $\B$-action

\begin{equation}\label{eqn:B-action-for-quotient}
b\cdot (g,z_1,\ldots, z_k)= (gb^{-1},b\cdot z_1,\ldots,b\cdot z_k).
\end{equation}

\np
Since $\G\times(\Z_{\vv_1}\times\cdots\times\Z_{\vv_k})$ is smooth and $\B$ acts without 
fixed points, the quotient $\Y_{\svseq}$ is smooth.  

\np
The group $\G$ acts on $\G\times(\Z_{\vv_1}\times\cdots\times\Z_{\vv_k})$ by left multiplication on the first factor.
Since this action commutes with the action of $\B$, it descends to an action of $\G$ on $\Y_{\svseq}$.
The map from $\G\times(\Z_{\vv_1}\times\cdots\times\Z_{\vv_k})$ to $\X^{k}$ given by 

\begin{equation}\label{eqn:fvseq-def}
(g,z_1,\ldots, z_k) \mapsto
\left({g\cdot f_{\vv_1}(z_1),\, g\cdot f_{\vv_2}(z_2)\,\ldots,\, g\cdot f_{\vv_k}(z_k)}\right)
\end{equation}

\noindent
is invariant under the $\B$-action.  If we let $\G$ act on $\X^k$ diagonally then \eqref{eqn:fvseq-def}
is also $\G$-equivariant
and hence 
descends to a $\G$-equivariant morphism $f_{\vseq}\colon\Y_{\svseq}\longrightarrow \X^{k}.$

\label{sec:Yconst:proj}
\np
As in the direct construction, we suppose that $\vv_j$ is nonempty, define $\uu_{l}$ by
\eqref{eqn:uudef} 
and set 
$\useq=(\uu_1,\ldots,\uu_{k})$.
The $\B$-equivariant morphisms $\pi_{\vv}\colon\Z_{\vv_j}\longrightarrow \Z_{\uu_{j}}$ and 
$\sigma_{\vv_j}\colon \Z_{\uu_{j}}\longrightarrow\Z_{\vv_j}$ from \S\ref{sec:BSDH}
give rise to $\B$-equivariant morphisms between 
$\G\times(\Z_{\vv_1}\times\cdots\times\Z_{\vv_k})$ and $\G\times(\Z_{\uu_1}\times\cdots\times\Z_{\uu_k})$
and hence to a $\G$-equivariant $\PP^1$-fibration $\pi_{\vseq,\useq}\colon\Y_{\svseq}\longrightarrow\Y_{\suseq}$
and a $\G$-equivariant section $\sigma_{\vseq,\useq}\colon\Y_{\suseq}\longrightarrow\Y_{\svseq}$.
These maps fit together to give diagram \eqref{eqn:bigfibredef}.

\bpoint{Expanded Version of the Direct Construction}
\label{sec:3constY}
Combining the formulas for $\P_{\vv}$ from \S\ref{secondconst} with the direct construction above
we obtain a more explicit expression for $\Y_{\svseq}$.
If $\vseq=(\vv_1,\ldots, \vv_k)$ with $\vv_{j}=s_{i_{1,j}}\cdots s_{i_{m_j,j}}$ for $j=1,\ldots,k$ then we define
$\Y_{\svseq}$ to be the quotient of 

$$ \G\times \P_{\vv_1}\times\cdots\times \P_{\vv_k} =
\G\times(\P_{i_{1,1}}\times\cdots\times \P_{i_{m_1,1}})\times\cdots\times(\P_{i_1,k}\times\cdots\times\P_{i_{m_k,k}})$$

\np
by the right action of $\B\times\B^{m_1}\times\cdots\times\B^{m_k}$, where an element

$$(b_0 \st b_{1,1},\ldots, b_{m_1,1} \st b_{1,2}, \ldots, b_{m_2,2} \st \cdots \st b_{1,k},\ldots, b_{m_k,k})$$

\np
acts from the right on

$$(g \st p_{i_{1,1}}, p_{i_{2,1}},\ldots, p_{i_{m_1,1}} \st p_{i_{1,2}}, p_{i_{2,2}}, \ldots, p_{i_{m_2,2}} \st 
\cdots \st p_{i_{1,k}},\ldots, p_{i_{m_k,k}})$$

\np
to give

$$(gb_0 \st b_0^{-1}p_{i_{1,1}}b_{1,1},\, b_{1,1}^{-1}p_{i_{2,1}}b_{2,1},\,\ldots,\, b_{m_1-1,1}^{-1}p_{i_{m_1,1}}b_{m_1,1} 
\st 
\cdots \st b_0^{-1}p_{i_{1,k}}b_{1,k},\,\ldots, \,b_{m_k-1,k}^{-1}p_{i_{m_k,k}}b_{m_k,k}).$$

\np
(In the expressions above the vertical lines ``$|$'' are used to indicate logical groupings, but otherwise
have no significance.)
The group $\G$ acts on $\G\times \P_{\vv_1}\times\cdots\times\P_{\vv_k}$ by left multiplication 
on the $\G$ factor, this action descends to a left action on $\Y_{\svseq}$.

\np
The map $f_{\vseq}$ is induced by the map sending an element 

$$({g \st p_{i_{1,1}}, p_{i_{2,1}},\ldots, p_{i_{m_1,1}} \st p_{i_{1,2}}, p_{i_{2,2}}, \ldots, p_{i_{m_2,2}} \st 
\cdots \st p_{i_{1,k}},\ldots, p_{i_{m_k,k}}})$$

\noindent
of $\G\times\P_{\vv_1}\times\cdots\times\P_{\vv_k}$ to

\begin{equation}\label{eqn:f-vseq:concrete}
({g p_{i_{1,1}} p_{i_{2,1}}\cdots p_{i_{m_1,1}} \st g p_{i_{1,2}} p_{i_{2,2}} \cdots p_{i_{m_2,2}} \st 
\cdots \st g p_{i_{1,k}}\cdots p_{i_{m_k,k}}})
\end{equation}

\noindent
in $\X^{k}$.  From the explicit formulas this is clearly a $\G$-equivariant map.

\np
Finally, if $\vseq$ is a sequence of words, and $\useq$ is a sequence obtained by dropping the rightmost 
reflection of a single word in $\vseq$ (as in \S\ref{sec:Yconst:proj}) then 
the $\G$-equivariant $\PP^1$-fibration $\pi_{\vseq,\useq}\colon \Y_{\svseq}\longrightarrow\Y_{\suseq}$ 
and 
the $\G$-equivariant section $\sigma_{\vseq,\useq}\colon\Y_{\useq}\longrightarrow\Y_{\vseq}$ 
are constructed using the obvious formulas analogous
to those in \S\ref{secondconst}. 
It again follows easily from these formulas that $f_{\useq}=f_{\vseq}\circ\sigma_{\vseq,\useq}$.

\np
{\bf Remark.} Note that the variety $\Y_{\svseq}$ depends on the sequence of words $\vseq=(\vv_1,\ldots, \vv_k)$
and not just on the corresponding sequence $(v_1,\ldots, v_k)$ of Weyl group elements.  If we choose a different
reduced factorization of each $v_i$ the resulting variety is birational to $\Y_{\svseq}$ over $\X^{k}$.
The proof is omitted because we do not need this fact.

\bpoint{The map $\fo$} \label{sec:fo}
As before, let $\vseq=(\vv_1,\ldots, \vv_k)$ be a sequence of words.
Besides the map $f_{\vseq}$ to $\X^{k}$, each $\Y_{\svseq}$ comes with a $\G$-equivariant map $\fo$ to $\X$
expressing $\Y_{\svseq}$ as a $\Z_{\vv_1}\times\cdots\times\Z_{\vv_k}$-bundle over $\X$.

\np
From the point of view of the construction in \S\ref{sec:1constY} $\fo$ is the composite map 

$$\Y_{\svseq}\stackrel{\sigma_{\vseq,\useq}}{\longrightarrow} \Y_{\useq} \longrightarrow \cdots \longrightarrow 
\Y_{\underline{\emptyset}}=\X$$

\np
obtained by dropping the elements in the entries of $\vseq$ one at a time.  The fibre over $e$ in $\X$ is 
then the result of applying the recursive construction in \S\ref{firstconst} separately for each $\vv_i$, $i=1$,\ldots,
$k$, and so the fibre is $\Z_{\vv_1}\times\cdots\times\Z_{\vv_k}$.

\np
From the point of view of the construction in \S\ref{sec:2constY} one starts with the projection 
$\G\times(\Z_{\vv_1}\times\cdots\times\Z_{\vv_k})\longrightarrow\G$ onto the first factor. This is $\B$-equivariant
for the right action of $\B$ on $\G$ and hence descends to a morphism $\fo\colon\Y_{\svseq}\longrightarrow \X$ 
expressing $\Y_{\svseq}$ as a $\Z_{\vv_1}\times\cdots\times\Z_{\vv_k}$-bundle over $\X$.

\np
Let $\useq=(\emptyseq,\vv_1,\ldots,\vv_k)$.
Since the action of $\B$ on the point $\Z_{\emptyseq}=e$ is trivial, we have 
an isomorphism

$$
\G\times \Z_{\vv_1}\times\cdots\times\Z_{\vv_k} 
\simeq
\G\times \Z_{\emptyseq}\times\Z_{\vv_1}\times\cdots\times\Z_{\vv_k} 
$$

\np
of $\B$-varieties and hence a $\G$-isomorphism $\phi\colon\Y_{\svseq}\longrightarrow\Y_{\suseq}$. 
From the explicit description in 
\eqref{eqn:fvseq-def}
we see that the composite map 
$f_{\useq}\circ\phi\colon\Y_{\svseq}\longrightarrow\X^{k+1}$ followed by projection onto the first factor
is $\fo$, and that $f_{\useq}\circ\phi$ followed by projection onto the last $k$ factors is $f_{\vseq}$.

\np
Thus the map $\fo\times f_{\vseq}\colon\Y_{\svseq}\longrightarrow \X\times\X^{k}$ is equal to the map
$f_{(\emptyseq,\vv_1,\ldots,\vv_k)}\colon\Y_{(\emptyseq,\vv_1,\ldots,\vv_k)}\longrightarrow\X^{k+1}$ under
the isomorphism $\phi$.  This will be used in the proof of Theorem \ref{thm:fibres}.

\bpoint{Maximum point}\label{sec:maxpt}
Let $\vseq=(\vv_1,\ldots,\vv_k)$ be a sequence of words.  We define the {\em maximum point}
$p_{\vseq}$ of $\Y_{\svseq}$ to be the product maximum point (\S\ref{sec:Zmaxpoint}) 
$p_{\vv_1}\times\cdots\times p_{\vv_k}$ in the fibre
$\Z_{\vv_1}\times\cdots\times\Z_{\vv_k}$ of $\fo$ over $e$ in $\X$.
Alternatively,  if $\vv_j=(s_{i_{1,j}},\ldots,s_{i_{m_j,j}})$ for $j=1$,\ldots, $k$ then (in the notation of
\S\ref{sec:3constY}) the point 
$$(e \st s_{i_{1,1}}, s_{i_{2,1}},\ldots, s_{i_{m_1,1}} \st \cdots \st s_{i_{1,k}},\ldots, s_{i_{m_k,k}})$$

\np
is a point of 

$$
\G\times(\P_{i_{1,1}}\times\cdots\times \P_{i_{m_1,1}})\times\cdots\times(\P_{i_1,k}\times\cdots\times\P_{i_{m_k,k}})$$

\np
and its image in $\Y_{\svseq}$ under the quotient map by $\B\times\B^{m_1}\times\cdots\times\B^{m_k}$ 
is the maximum
point $p_{\vseq}$.  If each $\vv_j$ is a factorization of some $v_j\in \W$, then the image $f_{\vseq}(p_{\vseq})$
of the maximum point in $\X^{k}$ is the point $q_{\vseq}:=(v_1,\ldots,v_k)$.

\bpoint{Tangent space formulas}
We will need to know the formal character (see \S\ref{sec:notation}) of the tangent space of $\Y_{\svseq}$ at the 
maximum point $p_{\vseq}$.  If each $\vv_j$ is a reduced word with product $v_j$,
then the formal character of the tangent space to $\Z_{\vv_j}$ 
at $p_{\vv_j}$ is $\langle \invset_{v_j^{-1}} \rangle$ and the  formal character of the 
tangent space of $\X$ at $e$ is $\langle \delneg \rangle$.

\np
Since the fibration $\fo$ is smooth, the formal character of $\T_{p_{\vseq}}\Y_{\svseq}$ is the
sum of these formal characters, i.e.,\ 

$$
\Ch(\T_{p_{\vseq}}\Y_{\svseq}) =  \langle \delneg \rangle +
 \sum_{i=1}^{k} 
\langle \invset_{v_i^{-1}} \rangle.
$$

\np
If  $v_j=w^{-1}_j\wo$ for $j=1$,\ldots, $k$, then by \eqref{eqn:weightcomplement} this is the same as

\begin{equation}\label{eqn:Yvtgt}
\Ch(\T_{p_{\vseq}}\Y_{\svseq}) = \delneg + \sum_{i=1}^k
\langle \invset_{w_i}^{\c}  \rangle.
\end{equation}

\bigskip
\np
\bpoint{Fibres and images of $f_{\vseq}$}\label{sec:fibreinfo}

\tpoint{Lemma} \label{lem:Xv-properties2}
Let $\vseq=(\emptyset,\vv_2,\ldots, \vv_k)$ be a sequence of words, 
with each $\vv_i$ a reduced factorization of $v_i$, 
and let $\X_{\vseq}$ be the (reduced) image of $f_{\vseq}$ in $\X^{k}$.  Then:

\begin{enumerate}

\item Projection onto the first factor of $\X^{k}$ endows $\X_{\vseq}$ with the structure of a 
fibre bundle over $\X$ with fibre isomorphic to $\X_{v_2}\times\cdots\times\X_{v_k}$.
\item The variety $\X_{\vseq}$ is normal with rational singularities of dimension $\N+\ell(\vseq)$, and 
the induced map $\Y_{\svseq}\longrightarrow\X_{\vseq}$ is birational with connected fibres.
\end{enumerate}

\medskip
\bpf
Projection on the first factor of $\X^{k}$ gives a $\G$-equivariant morphism 
$\X_{\vseq}\stackrel{\eta}{\longrightarrow} \X$.  Since
$\G$ acts transitively on $\X$ this morphism is surjective and all fibres are isomorphic, i.e., this expresses
$\X_{\vseq}$ as a fibre bundle over $\X$.  To study the fibres we look at the fibre $\eta^{-1}(e)$ 
over the $\B$-fixed point $e$ of $\X$.

\np
Consider the diagram

\begin{equation}\label{eqn:fv-diagram}
\begin{array}{c}
\xymatrix{
\G\times e\times\Z_{\vv_2}\times\cdots\times\Z_{\vv_k}\ar[rr]^(.7){\psi_{\vseq}}\ar[dd]_{\id_{\G}\times \id_{x_o}\times f_{\vv_1}\times\cdots\times f_{\vv_k}} & & \Y_{\subsmash{\svseq}} \ar[dd]^{f_{\vseq}} \\
\\
\G\times e\times\X_{v_2}\times\cdots\times\X_{v_k}\ar[rr]^(.7){\phi}  & & \X^{\subsmash{k}}  \\
}
\end{array}
\end{equation}

\np
where $\phi$ is given by $\phi(g,e,x_2,\ldots, x_k)=(g\cdot e, g\cdot x_2,\ldots, g\cdot x_k)\in \X^{k}$.
Since $\psi_{\vseq}$ and the leftmost vertical map are
surjective, the image of $f_{\vseq}$ is the same as the image of $\phi$. 
Since $\B$ is the stabilizer of $e$, the fibre $\eta^{-1}(e)$ is the image of 
$\B\times e\times\X_{v_2}\times\cdots\times\X_{v_k}$ under $\phi$.  But each Schubert
variety $\X_{w}$ is stable under the action of $\B$ and therefore the image above is just 
$e \times \X_{v_2}\times\cdots\times\X_{v_k}$, proving ({\em a}).

\np
From the fibration $\eta$ it is clear that 
$$\dim(\X_{\vseq})=\dim(\X)+\sum_{i=2}^{k} \dim(\X_{v_i}) = \N + \sum_{i=2}^{k} \ell(v_i) = \N + \ell(\vseq)$$
because each $\vv_i$ is reduced and hence $\ell(v_i)=\ell(\vv_i)$ for $i\geq 2$.

\np
The product of normal varieties is again normal, and the product of varieties with rational singularities also has
rational singularities.  Since each $\X_{w}$ is normal with rational singularities 
(\S\ref{sec:schubertvarieties}), the fibres also have this property, 
and therefore so does $\X_{\vseq}$ (since the properties of being normal or having rational singularities are local,
and $\X_{\vseq}$ is locally the product of the fibre and a smooth variety).  

\np
Since each map $f_{\vv_i}\colon\Z_{\vv_i}\longrightarrow\X_{v_i}$ is a resolution of singularities of a normal
variety, each $f_{\vv_i}$ is birational with connected fibres. It follows that 
the map $\Y_{\svseq}\longrightarrow\X_{\vseq}$,
which is the 
quotient of the leftmost vertical map in \eqref{eqn:fv-diagram} by the action of $\B$ 
is also birational with connected fibres.  This proves ({\em b}). \epf

\tpoint{Definition} If $\X_{w}$ is any Schubert subvariety of $\X$ 
and $q$ is any point of $\X$, we define the subvariety $q\X_{w}$ of $\X$ to be the 
result of translating $\X_{w}$ by any element in the $\B$-coset corresponding to $q$.
Since $\X_{w}$ is $\B$-stable 
the result is independent of the choice of representative for $q$.

\np
The following theorem gives more precise information about the image and fibres of $f_{\svseq}$. 

\tpoint{Theorem} 
\label{thm:fibres}
Let $\vseq=(\vv_1,\ldots,\vv_k)$ be a sequence of reduced words with corresponding Weyl group elements 
$(v_1,\ldots, v_k)$.
Then there exists a factorization 
$f_{\vseq}\colon \Y_{\svseq}\stackrel{\tau}{\longrightarrow}\Q_{\vseq}\stackrel{h}{\longrightarrow} \X^{k}$ such that

\begin{enumerate}
\item $\Q_{\vseq}$ is normal with rational singularities; 
\item the map $\tau\colon\Y_{\svseq}\longrightarrow\Q_{\vseq}$ is proper and birational with connected fibres; 
\item for each point $(q_1,\ldots, q_k)$ of $\X^{k}$ there is a natural inclusion 
of the scheme-theoretic fibre $h^{-1}(q_1,\ldots,q_k)$ 
into the scheme-theoretic intersection
$\bigcap_{i=1}^{k} q_i\X_{v_i^{-1}}$;
\item the inclusion of schemes in ({\em c}) induces an isomorphism at the level of reduced schemes,
or in other words,
the set theoretic fibre $h^{-1}(q_1,\ldots,q_k)$ is equal to the set theoretic intersection 
$\bigcap_{i=1}^{k} q_i\X_{v_i^{-1}}$. 
\end{enumerate}

\np
\bpf
Let 
$\fo\times f_{\vseq}\colon\Y_{\svseq}\longrightarrow \X\times\X^{k}$ be the product of $f_{\vseq}$ and the map
$\fo\colon\Y_{\svseq}\longrightarrow \X$ from \S\ref{sec:fo} expressing $\Y_{\svseq}$ as a 
$\Z_{\vseq_1}\times\cdots\times \Z_{\vseq_k}$-bundle over $\X$.  We define $\Q_{\svseq}$ to be the image of 
$\fo\times f_{\vseq}$ with the reduced scheme structure, $\tau$ to be the map from $\Y_{\svseq}$ onto $\Q_{\vseq}$, 
and $h$ to be the map from $\Q_{\vseq}$ to $\X^{k}$ induced by the projection $\X\times\X^{k}\longrightarrow \X^{k}$.
By construction $f_{\vseq}=h\circ \tau$.

\np
Letting $\psi_{\vseq}$  be the map (from \S\ref{sec:2constY})
defining $\Y_{\svseq}$ as a quotient of $\B\times\Z_{\vv_1}\times\cdots\times \Z_{\vv_k}$  and
$\phi\colon \G\times \X_{v_1} \times \cdots \times \X_{v_k}\longrightarrow \Q_{\vseq}\subseteq \X\times\X^{k}$ 
as the map sending $(g,x_1,\ldots, x_k)$ to $(g\B/\B,g\cdot x_1,\ldots, g\cdot x_k)$ in $\X\times\X^{k}$,
we obtain a refinement of diagram \eqref{eqn:fv-diagram}:

$$
\xymatrix{
\G\times \Z_{\vv_1}\times \cdots \times \Z_{\vv_k}\ar[rr]^(0.6){\psi_{\vseq}}\ar[dd]\ar@{}[ddrr]|\Box & &
\Y_{\svseq}\ar[dd]_{\tau} \ar[ddr]^{\fo\times f_{\vseq}}
\\
\\
\G\times \X_{v_1} \times \cdots \times \X_{v_k}\ar[rr]^(0.6){\phi} & & 
\Q_{\vseq}\ar@{^{(}->}[r]\ar[d]_{h} & \X\times \X^{k}\ar[d] \\
 & & \X^{k} \ar@{=}[r] & \X^{k} \\
}
$$

\np
Since $\Y_{\vseq}\simeq\Y_{(\emptyset,\vv_1,\ldots,\vv_k)}$ (see \S\ref{sec:fo}) and 
under this isomorphism the map $\fo\times f_{\vseq}$ is the map $f_{(\emptyset, \vv_1,\ldots, \vv_k)}$,
it follows from Lemma
\ref{lem:Xv-properties2}({\em b}) that $\Q_{\vseq}$ 
is normal with rational singularities and that $\tau\colon\Y_{\svseq}\longrightarrow \Q_{\vseq}$  
is birational with connected fibres, proving ({\em a}) and ({\em b}).

\np
The composite map 
$$ \G\times \P_{\vv_1}\times\cdots\times \P_{\vv_k} \longrightarrow 
\G \times \Z_{\vv_1}\times\cdots\times \Z_{\vv_k} \stackrel{\psi_{\vseq}}{\longrightarrow}
\Y_{\svseq} \stackrel{\tau}{\longrightarrow} \Q_{\vseq} $$

\np
is given (in the notation of \S\ref{sec:3constY}) by sending 
$$({g \st p_{i_{1,1}}, p_{i_{2,1}},\ldots, p_{i_{m_1,1}} \st p_{i_{1,2}}, p_{i_{2,2}}, \ldots, p_{i_{m_2,2}} \st 
\cdots \st p_{i_{1,k}},\ldots, p_{i_{m_k,k}}})$$

\np
to

$$(g\st {g p_{i_{1,1}} p_{i_{2,1}}\cdots p_{i_{m_1,1}} \st g p_{i_{1,2}} p_{i_{2,2}} \cdots p_{i_{m_2,2}} \st 
\cdots \st g p_{i_{1,k}}\cdots p_{i_{m_k,k}}})$$

\np
in $\X\times\X^{k}$.  A point $q$ of
$\X$ is therefore in the fibre $h^{-1}(q_1,\ldots, q_k)\subseteq \X\times q_1\times\cdots\times q_k=\X$ 
if for any $\B$-coset representatives 
$g$, $g_1$, \ldots, $g_k$ of $q$, $q_1$,\ldots, $q_k$, there exist elements $\left\{p_{i,j}\right\}$ 
in the respective parabolic subgroups such that we can solve the equations 

\begin{eqnarray*}
g p_{i_{1,1}} p_{i_{2,1} }\cdots p_{i_{m_1,1}} & = & g_1 \\
\vdots \rule{1.3cm}{0cm}& \vdots & \vdots \\
g p_{i_{1,k}} p_{i_{2,k}}\cdots p_{i_{m_k,k}} & = & g_k. \\
\end{eqnarray*}

\np
Moving the $p_{i,j}$'s to the right hand side, the system above becomes

\begin{eqnarray*}
g & = & g_1 p_{i_{m_1,1}}^{-1} \cdots p_{i_{2,1}}^{-1}  p_{i_{1,1}}^{-1}  \\
\vdots & \vdots & \rule{1.3cm}{0cm}\vdots \\
g & = & g_k p_{i_{m_k,k}}^{-1} \cdots p_{i_{2,k}}^{-1}p_{i_{1,k}}^{-1} \\
\end{eqnarray*}

\np
which is equivalent to $q$ belonging in the intersection 
$\displaystyle{\bigcap_{i=1}^{k} q_i\X_{v_i^{-1}}},$ proving ({\em d}).

\np
Let $\vv$ be a reduced word with product $v$.
By part ({\em d}) the set
$$\Q_{\vv}':=\left\{{ (q,p)\in\X\times \X \st p\in q\X_{v}\rule{0cm}{0.4cm}}\right\}$$
is the image of 
$f_{(\emptyset,\vv)}\colon\Y_{(\emptyset,\vv)}\longrightarrow \X\times\X$ and is therefore a closed subvariety of
$\X\times\X$.
Alternatively $\Q'_{\vv}$ is the Zariski closure of the 
set $\left\{{ (g,g\cdot v) \st g \in \G}\right\}\subseteq \X\times\X$.

\np
For $i=1$, \ldots, $k$, let $p_i\colon \X\times\X^{k}\longrightarrow \X\times \X$ be the map which is the product of
$\id_{\X}$ with projection $\X^{k}\longrightarrow \X$ onto the $i$-th factor.  The intersection
$$\Q'_{\vseq}:=\bigcap_{i=1}^{k} p_i^{-1}(\Q'_{\vv_i})$$

\np
is a closed subscheme of $\X\times\X^{k}$ which, by ({\em d}), agrees set theoretically with $\Q_{\vseq}$.  Since
$\Q_{\vseq}$ is reduced, we have the inclusion of schemes $\Q_{\vseq}\subseteq \Q'_{\svseq}$.  
If $h'$ is the map $h'\colon\Q'_{\vseq}\longrightarrow \X^{k}$ induced by projection, then the scheme-theoretic fibres
of $h$ are naturally a subscheme of the scheme-theoretic fibres of $h'$ (and both are naturally subschemes of $\X$).
The scheme-theoretic fibre of $h'$ is the scheme-theoretic intersection 
$\bigcap_{i=1}^{k} q_i\X_{v_i^{-1}}$, proving ({\em c}). \epf

\np
The image $\X_{\vseq}$ is therefore the set of translations $(q_1,\ldots, q_k)$ in $\X^{k}$ 
for which the intersection $\bigcap_{i=1}^{k} q_i\X_{v_i^{-1}}$
of translated Schubert varieties is non-empty, and the 
set-theoretic fibres of $h$ are the intersections themselves.  
Moreover, $\Q_{\vseq}$ is the incidence correspondence of intersections of translates of Schubert varieties
(the first coordinate in $\X\times\X^{k}$ is the intersection, the remaining $k$ coordinates are the parameters 
$(q_1,\ldots, q_k)$ controlling the translates).  
Theorem \ref{thm:fibres} shows that $\Y_{\svseq}$ is a resolution of singularities of $\Q_{\vseq}$.

\tpoint{Corollary} \label{cor:map-degree}
Let $\vseq = (\vv_1,\ldots, \vv_k)$ be a sequence of reduced words with corresponding Weyl group elements
$(v_1,\ldots, v_k)$  such that $\sum_{i=1}^{k} \ell(\vv_i) = (k-1)\N$.
Then the degree of the map 
$f_{\vseq}\colon \Y_{\svseq}\longrightarrow \X^{k}$ is given by the intersection number 
$\bigcap_{i=1}^{k} [\Omega_{\wo v_i^{-1}}]=\bigcap_{i=1}^{k}[\X_{v_i^{-1}}]$.

\np
{\bf Remark.} The dimension of $\Y_{\svseq}$ in this case is $\N+\sum \ell(\vv_i) = k\N = \dim(\X^{k})$ so 
it is reasonable to ask for the degree of the map.

\medskip
\bpf
Since we are working in characteristic zero, the degree of $f_{\vseq}$ is given by the number of points 
in a generic fibre.
By Theorem \ref{thm:fibres} the map $p\colon\Y_{\svseq}\longrightarrow \Q_{\vseq}$ is birational, 
and so the generic fibre of $f_{\vseq}$ is the same as the generic fibre of 
$h\colon \Q_{\vseq}\longrightarrow \X^{k}$. 
By the Kleiman transversality theorem, if $q_1$,\ldots, $q_k$ are generic, the scheme-theoretic intersection
$\displaystyle{\cap_{i=1}^{k} q_i\X_{v_i^{-1}}}$ is reduced and finite, and 
the number of points is equal to the intersection number 
$\displaystyle{\cap_{i=1}^{k} [\X_{v_i^{-1}}}] 
=
\cap_{i=1}^{k} [\Omega_{\wo v_i^{-1}}]$
in $\H^{*}(\X,\ZZ)$.  
By Theorem \ref{thm:fibres}({\em c--d}) if the scheme-theoretic intersection 
$\displaystyle{\cap_{i=1}^{k} q_i\X_{v_i^{-1}}}$ is reduced it is equal to the scheme-theoretic
fibre $h^{-1}(q_1,\ldots, q_k)$, proving the corollary. \epf

\bpoint{Key Lemma}
We now prove an important lemma which will allow us to derive several results necessary for the proofs of
Theorems I and II.  The lemma itself will also be used in the proof of Theorem I.  

\tpoint{Lemma}\label{lem:key}
Let $\vseq$ be a sequence of reduced words, $\L$ be a $\G$-equivariant 
line bundle on
$\Y_{\svseq}$ and $s\in\H^{0}(\Y_{\svseq},\L)^{\G}$ be a nonzero $\G$-invariant section.  Then 

\begin{enumerate}
\item the weight of $\L$ at the $\T$-fixed maximum point (\S\ref{sec:maxpt}) 
$p=p_{\vseq}\in \Y_{\svseq}$ belongs to $\Span_{\ZZ_{\geq 0}} \Delta^+$; 
\item the weight of $\L$ at $p$ is zero if and only if $s$ does not vanish at $p$;
\item without supposing that $\L$ has a $\G$-invariant section, if $\L$ is an equivariant bundle on $\Y_{\svseq}$
and the weight of $\L$ at $p$ is zero, then $\dim \H^{0}(\Y_{\svseq},\L)^{\G}\leq 1$.
\end{enumerate}

\np
{\bf Remark.} 
Part ({\em c}) will be used repeatedly to control the size of the $\G$-invariant sections.

\bpf
Let $\fo\colon \Y_{\svseq}\longrightarrow \X$ be the map 
from \S\ref{sec:fo}
expressing $\Y_{\svseq}$ as a $\Z_{\vv_1}\times \cdots\times \Z_{\vv_k}$-bundle over $\X$.
The section $s$ cannot vanish on any fibre of $\fo$ 
since (by $\G$-invariance and transitivity of $\G$-action on $\X$) $s$ would vanish on all of $\Y_{\svseq}$.  
We can thus restrict $s$ to get a nonzero section on the fibre $\Z_{\vv_1}\times\cdots\times \Z_{\vv_k}$ of 
$\fo$ over $e\in \X$; this fibre contains the maximum point $p$.

\np
The formal character of the tangent space at the maximum point $p_i$ of $\Z_{\vv_i}$ is
$\langle \invset_{v_i^{-1}}\rangle$;
i.e., all the weights of this space are positive roots.  Since the maximum point 
$p=p_1\times\cdots \times p_k \in \Z:=\Z_{\vv_1}\times\cdots\times\Z_{\vv_k}$ is the product of the maximum points of
the factors, each of the weights on the tangent space of $p$ in $\Z$ is also a positive root.

\np
Let $\m_p$ be the maximal ideal of $p$ in $\Osh_{\Z,p}$.  For every $r\geq 0$ we get a $\T$-equivariant 
restriction map 

$$\H^{0}(\Z,\L|_{\Z}) \longrightarrow \L\otimes_{\Osh_{\Z}} \left({\Osh_{\Z,p}/\m_p^{r+1}}\right) = \L\otimes 
\left({\Osh_{\Z}/\m_p \oplus \m_p/\m_p^{2}}\oplus\cdots\oplus \m_p^{r}/\m_p^{r+1}\right)
$$

\np
which is an injection for $r$ sufficiently large.  In particular, for sufficiently large $r$, the section $s$ restricts
to a nonzero element of $\L\otimes_{\Osh_{\Z}} \Osh_{\Z,p}/\m_p^{r+1}$.  Since $s$ is an invariant section, this means
that the zero weight is a weight of  $\L\otimes_{\Osh_{\Z}} (\Osh_{\Z,p}/\m_p^{r+1})$, and
so must appear in one of the factors 
$\L\otimes_{\Osh_{\Z}}( \m_p^{i}/\m_p^{i+1})=\L\otimes_{\Osh_{\Z}}\Sym^{i}(\m_p/\m_p^2)$ for $i=0,\ldots,r$.  

\np
Since $\m_p/\m_p^2$ is dual to the tangent space at $p$, all weights of $\m_p/\m_p^2$ are negative roots,
and therefore the weights of $\Sym^{i}(\m_p/\m_p^2)$ belong to $\Span_{\ZZ\leq0} \delpos$.  Tensoring with $\L$ 
multiplies the formal character of $\Sym^{i}(\m_p/\m_p^2)$ by the weight of 
 $\L$ at $p$. Thus the zero weight is a weight of  $\L\otimes_{\Osh_{\Z}}\Sym^{i}(\m_p/\m_p^2)$
 only if the weight of  $\L$ at $p$ belongs to $\Span_{\ZZ_{\geq 0}} \Delta^+$. This proves ({\em a}).

\np 
The value of $s$ at $p$ is the restriction of $s$ to the factor
$\L\otimes_{\Osh_{\Z}}(\Osh_{\Z,p}/\m_p )= \L_p$.    If $s$ does not vanish at $p$ the weight of  $\L_p$ is therefore 
zero.  Conversely, if the weight of $\L_p$ is zero then the weights of
$\L\otimes_{\Osh_{\Z}} \Sym^{i}(\m_p/\m_p^2)$ 
are non-zero for $i\geq 1$. Hence the only possibility for the invariant section 
$s$ under the restriction map is to have nonzero restriction to $\L\otimes_{\Osh_{\Z}}(\Osh_{\Z,p}/\m_p )= \L_p$,
proving ({\em b}).

\np
Suppose that the weight of $\L$ at $p$ is zero.
If there were two linearly independent sections $s_1,s_2\in \H^{0}(\Y_{\svseq},\L)^{\G}$ 
then some nonzero linear combination would vanish at $p$ contradicting ({\em b}). 
Hence if the weight is zero we must have $\dim \H^{0}(\Y_{\svseq},\L)\leq 1$, giving ({\em c}).
\epf

\bpoint{Applications of Lemma \ref{lem:key}}
\label{sec:keyconsequences}

\tpoint{Theorem}\label{thm:hammer}
Suppose that $w_1$, \ldots, $w_k$, and $w$ are elements of the Weyl group such that $\ell(w)=\sum_{i=1}^{k} \ell(w_i)$
and $\displaystyle{\bigcap_{i=1}^{k}[\Omega_{w_i}]\cdot[\X_{w}]\neq 0}$ in $\H^{*}(\X,\ZZ)$.  Then:

\begin{enumerate}
\item \rule{0cm}{0.6cm}For any dominant weights $\mu_1$, \ldots, $\mu_k$, and $\mu$ such that 
the irreducible module $\V_{\mu}$ is a component of $\V_{\mu_1}\otimes\cdots\otimes \V_{\mu_k}$, the 
weight $\sum_{i=1}^{k} w_i^{-1}\mu_i - w^{-1}\mu$ belongs to $\Span_{\ZZ_{\geq 0}} \Delta^+$.
\item If $\sum_{i=1}^{k} w_i^{-1}\mu_i - w^{-1}\mu=0$ then 
$\mult(\V_{\smu},\V_{\smu_1}\otimes\cdots\otimes\V_{\smu_k})=1$.
\rule{0cm}{0.6cm}
\item $\sum_{i=1}^{k} w_i^{-1}\cdot 0 - w^{-1}\cdot 0 = 
\sum_{i=1}^{k} (w_k^{-1}\rho-\rho) - (w^{-1}\rho-\rho)$ belongs to $\Span_{\ZZ_{\geq 0}} \Delta^+$. \rule{0cm}{0.6cm}
\item If $\sum_{i=1}^{k} w_i^{-1}\cdot 0 = w^{-1}\cdot 0$ then 
$\invset_{w} = \bigsqcup_{i=1}^{k} \invset_{w_i}$.  \rule{0cm}{0.6cm}
\end{enumerate}

\np
Note that the action of the Weyl group in parts ({\em a}) and ({\em b}) is the homogeneous
action, while the action in parts ({\em c}) and ({\em d}) is the affine action.

\np
\bpf Let $v_i=w_i^{-1}\wo$ for $i=1,\ldots, k$, $v_{k+1}=w^{-1}$, let $\vv_i$ be a reduced word with product $v_i$,
for $i=1$,\ldots, $k+1$, and set $\vseq=(\vv_1,\ldots, \vv_{k+1})$.
Then $\sum \ell(\vv_i)=(k+1-1)\N$ and so, by Corollary \ref{cor:map-degree}, the degree of 
$f_{\vseq}\colon\Y_{\svseq}\longrightarrow \X^{k+1}$ is given by the intersection number 

$$\bigcap_{i=1}^{k+1} [\Omega_{\wo v_i^{-1}}]= \bigcap_{i=1}^{k}[\Omega_{w_i}]\cdot[\X_{w}].$$

\np
By hypothesis this intersection number is nonzero and therefore $f_{\vseq}$ is surjective.

\np
Given dominant weights $\mu_1$, \ldots, $\mu_k$, and $\mu$ let $\lambda_i=-\wo\mu_i$ for 
$i=1\ldots, k$ and $\lambda_{k+1}=\mu$.  
Set $\L$ to be the line bundle $\L_{\lambda_1}\bt\cdots\bt\L_{\lambda_{k+1}}$ on $\X^{k+1}$,
so that $\H^{0}(\X^{k+1},\L)= \V_{\mu_1}\otimes\cdots\otimes\V_{\mu_k}\otimes\V_{\mu}^{*}$ and 
$\dim \H^{0}(\X^{k+1},\L)^{\G}$ is the multiplicity of $\V_{\mu}$ in the tensor product.
$\V_{\mu_1}\otimes\cdots\otimes\V_{\mu_k}$.

\np
Since $f_{\vseq}$ is surjective, pullback induces an inclusion 

$$\H^{0}(\Y_{\svseq},f_{\vseq}^{*}\L) \stackrel{f_{\vseq}^{*}}{\longleftarrow} \H^{0}(\X^{k+1},\L)$$

\np
and, in particular, $\dim \H^{0}(\Y_{\svseq},f_{\vseq}^{*}\L)^{\G}\geq \dim\H^{0}(\X^{k+1},\L)^{\G}$.  Applying
Lemma \ref{lem:key}({\em a}), we know that if $f_{\vseq}^{*}\L$ has a nonzero $\G$-invariant section then the weight
of  $f_{\vseq}^{*}\L$ at the maximum point $p_{\vseq}$ belongs to $\Span_{\ZZ_{\geq 0}} \Delta^+$.    
This weight is 

\begin{equation}\label{eqn:lemqweight}
\sum_{i=1}^{k+1} v_i(-\lambda_i) = \sum_{i=1}^{k} (w_i^{-1} \wo)(\wo\mu_i) + w^{-1}(-\mu) = 
\sum_{i=1}^{k} w_i^{-1}\mu_i - w^{-1}\mu,
\end{equation}

\np
proving ({\em a}).

\np
If the weight in \eqref{eqn:lemqweight} is zero then 
$\dim\H^{0}(\X^{k+1},\L)^{\G} \leq \dim \H^{0}(\Y_{\svseq},f_{\vseq}^{*}\L)^{\G} \leq 1$
by Lemma \ref{lem:key}({\em c}), and so if $\V_{\mu}$ is
a component of 
$\V_{\mu_1}\otimes\cdots\otimes\V_{\mu_k}$ then it is of multiplicity at most one.  The fact that $\V_{\mu}$ actually
is a component of the tensor product is a consequence of the solution of the PRV conjecture --
see \S\ref{sec:PRVdiscuss} for a discussion. This proves ({\em b}).

\np
The map $f_{\vseq}\colon\Y_{\svseq}\longrightarrow \X^{k+1}$ induces a natural map 
$f_{\vseq}^{*}\K_{\X^{k+1}}\longrightarrow \K_{\Y_{\svseq}}$ which is given by a global section $s$ of 
$\H^0(\Y_{\svseq},(f_{\vseq}^{*}\K_{\X^{k+1}})^{*}\otimes \K_{\Y_{\svseq}})$.  Since $\Y_{\svseq}$ and $\X^{k+1}$ have
the same dimension and since $f_{\svseq}$ is surjective this section is nonzero.  
Because the pullback morphism is natural, the section $s$ is $\G$-invariant.
By Lemma \ref{lem:key}({\em a}) the weight of the line bundle 
$\K_{\Y_{\svseq}/\X^{k+1}}:=(f_{\vseq}^{*}\K_{\X^{k+1}})^{*}\otimes\K_{\Y_{\svseq}}$ at the maximum point $p_{\vseq}$ 
belongs to $\Span_{\ZZ\geq0}\delpos$.

\np
By \eqref{eqn:Yvtgt}, \eqref{eqn:weightcomplement}, and \eqref{eqn:tgt-at-winv} the formal characters of
 the tangent spaces at $p_{\vseq}$ in $\Y_{\svseq}$ and $q_{\vseq}:=f_{\vseq}(p_{\vseq})$ 
in $\X^{k+1}$ are, respectively:

\begin{equation}\label{eqn:ChTp}
\Ch(\T_{p_{\vseq}}\Y_{\svseq}) =  \langle \invset_{w} \rangle + \langle \delneg \rangle +
\sum_{i=1}^{k} \langle \invset_{w_i}^{\c} \rangle
\end{equation}

\np
and

\begin{equation}\label{eqn:ChTq}
\Ch(\T_{q_{\vseq}}\X^{k+1})= 
\left({\langle \invset_{w} \rangle + \langle -\invset_{w}^{\c} \rangle}\right)
+ \sum_{i=1}^{k}  \left({\langle \invset_{w_i}^{\c} \rangle  + \langle -\invset_{w_i} \rangle}\right).
\end{equation}

\np
A short calculation using formula \eqref{eqn:rootsum}
shows that the  weight of  $\K_{\Y_{\svseq}/\X^{k+1}}$ at $p_{\vseq}$ is 
$\sum_{i=1}^{k} (w_k^{-1}\rho-\rho) - (w^{-1}\rho-\rho)$, proving ({\em c}).

\np
If the weight
$\sum_{i=1}^{k} (w_k^{-1}\rho-\rho) - (w^{-1}\rho-\rho)$ is zero then, by Lemma \ref{lem:key}({\em b}),
the section $s$ is nonzero at $p_{\vseq}$.  This means that 
$f_{\vseq}$ is unramified at $p_{\vseq}$ and therefore the tangent space map 
$\T_{\Y_{\svseq},p}\stackrel{df_{\vseq}}{\longrightarrow} \T_{\X^{k+1},q}$ is an isomorphism.
Hence both spaces must have the same formal characters.
Comparing the negative roots and their multiplicities in 
\eqref{eqn:ChTp}
and
\eqref{eqn:ChTq}
gives 
$\delneg= (\sqcup_{i=1}^{k} -\invset_{w_i}) \bigsqcup -\invset_{w}^{\c}$
which is equivalent to $\invset_{w} = \bigsqcup_{i=1}^{k} \invset_{w_i}$, proving ({\em d}). \epf

\bpoint{Relation with existing results}\label{sec:attribution}
Part ({\em a}) of  Theorem \eqref{thm:hammer}
is due to Berenstein and Sjamaar.
A theorem of this type was first proved by Klyacho \cite{kly} for $\Gl_n$.
This was later extended to all semisimple groups by Berenstein and Sjamaar in \cite{BeSj} and 
by Kapovich, Leeb, and Milson in \cite{klm}.  
Parts ({\em c}) and ({\em d}) are due to Belkale and Kumar:
part ({\em c}) is \cite[Theorem 29]{BK1} and ({\em d}) is \cite[Theorem 15]{BK1}, both in the case when the 
parabolic group $\P$ is the Borel group $\B$.  

\np
Part ({\em b}) is new and crucial for controlling the multiplicities of cohomological components.
The remaining statements have been included because 
Lemma \ref{lem:key} allows us to give a new, short, and unified proof of these results.
In particular, we obtain a new proof of the necessity of the inequalities determining the Littlewood-Richardson cone.
Namely, these inequalities are obtained by requiring that the weights in Theorem \ref{thm:hammer}({\em a})
(for all $w_1$,\ldots, $w_k$, $w$ satisfying the conditions of the theorem) belong to $\Span_{\ZZ\geq 0}\delpos$.
(The proof that these inequalities are sufficient requires a separate GIT argument.)

\np
{\bf Relation with a construction of Kumar.}
\label{sec:compareKumar}
Given a sequence $\uu$ of simple reflections, Kumar \cite[\S1.1]{ku2} defined a variety $\widetilde{\Z}_{\uu}$
along with a map $\theta_{\uu}$ from $\widetilde{\Z}_{\uu}$ to $\X^{2}$.    For any pair of words $\vseq=(\vv_1,\vv_2)$
let $\uu=\vv_1^{-1}\vv_2$ be the word obtained by reversing $\vv_1$ and concatenating it onto the left of $\vv_2$.
By comparing the construction of $\Y_{\svseq}$ and $\widetilde{\Z}_{\uu}$ 
it is not hard to find an isomorphism $\widetilde{\Z}_{\uu} = \Y_{\svseq}$ over $\X^2$ (i.e., such 
that $\theta_{\uu}=f_{\vseq}$ under the isomorphism).  Therefore when $k=2$ the varieties produced by our construction
are the same the ones constructed in \cite[\S1.1]{ku2}.

\section{Proof of Theorem III}

\np
\point
We will prove Theorem III in its symmetric form.  After applying the symmetrization procedure from 
\S\ref{sec:nonsymmetric-symmetric}  (and replacing $k+1$ by $k$) we obtain:

\tpoint{Theorem} {\em (Symmetric form of Theorem III)} \label{thm:topsurjectivity}
Let $w_1$, \ldots, $w_k$ be elements of $\W$ group such that $\sum_i \ell(w_i)=\N$; 
$\lambda_1$, \ldots, $\lambda_k$ be weights such that $w_i\cdot\lambda_i$ are dominant weights for $i=1,\ldots, k$;
and $\sum_{i=1}^{k} \lambda_i=-2\rho$. 

\begin{enumerate}
\item[({\em a})] If $\displaystyle{\bigcap_{i=1}^{k}[\Omega_{w_i}]=1}$ then the cup product map 

\begin{equation}\label{eqn:symmetriccup}
\H^{\ell(w_1)}(\X,\L_{\lambda_1})\otimes\cdots\otimes \H^{\ell(w_k)}(\X,\L_{\lambda_k})
\stackrel{\cup}{\longrightarrow}
\H^{\N}(\X,\K_{\X})
\end{equation}
is surjective.
\vskip 0.1cm
\item[({\em b})]  If $\displaystyle{\bigcap_{i=1}^{k}[\Omega_{w_i}]=0}$ then 
\eqref{eqn:symmetriccup} is zero.
\end{enumerate}

\np
The proof of Theorem \ref{thm:topsurjectivity} is given in \S\ref{proof-of-top-dim}.  
We will use the following common notation.
For any sequence $\ll=(\lambda_1,\ldots, \lambda_k)$ of weights let $\L_{\ll}$ be the line bundle 

$$\L_{\ll} := 
\L_{\lambda_1}\bt \cdots \bt \L_{\lambda_k} = 
pr_{1}^{*}\L_{\lambda_1}\otimes\cdots\otimes pr_{k}^{*}\L_{\lambda_k}
$$

\np
on $\X^{k}$, where $pr_i\colon\X^{k}\longrightarrow\X$ denotes projection onto the $i$-th factor.

\bpoint{Inductive Lemma}
Let $\ll=(\lambda_1,\ldots, \lambda_k)$ be a sequence of weights and $\vseq=(\vv_1,\ldots,\vv_k)$ a sequence
of words.  
Let $\useq$ be a sequence of words as in \eqref{eqn:uudef}, i.e., $\useq$ is a sequence of words
obtained by dropping a simple reflection from the right of a single member of $\vseq$.  
The following lemma lets us propagate information about the pullback map 

\begin{equation}
\label{eqn:top-map-hypoth}
\H^{\N+\ell(\vseq)}(\Y_{\svseq},f^{*}_{\vseq}\L_{\ll}) \stackrel{f_{\vseq}^{*}}{\longleftarrow} 
\H^{\N+\ell(\vseq)}(\X^{k},\L_{\ll})
\end{equation}

\np
on the top degree cohomology of $\Y_{\svseq}$ to information about an analogous pullback map to the 
top degree cohomology of $\Y_{\suseq}$.  
If $\vv_j=s_{i_1}\cdots s_{i_m}$, so that we are dropping $s_{i_m}$ from $\vv_j$ to get $\uu_j$, we denote by 
$\mm$ the sequence $\mm:=(\lambda_1,\ldots,\lambda_{j-1},s_{i_m}\cdot\lambda_j,\lambda_{j+1},\ldots, \lambda_k)$.
Finally, we assume that the degree of $\L_{\ll}$ is negative on the fibres of the $\PP^1$-fibration
$\pi_{\vseq,\useq}\colon\Y_{\svseq}\longrightarrow\Y_{\suseq}$.

\tpoint{Lemma} 
\label{lem:induct}
Under the conditions above, the pullback map 

$$\H^{\N+\ell(\useq)}(\Y_{\suseq},f^{*}_{\useq}\L_{\mm}) \stackrel{f_{\useq}^{*}}{\longleftarrow} 
\H^{\N+\ell(\useq)}(\X^{k},\L_{\mm})$$

\np
is ({\em a}) surjective, ({\em b}) zero, or ({\em c}) surjective on the space of $\G$-invariants if the pullback
map \eqref{eqn:top-map-hypoth} has the corresponding property ({\em a}), ({\em b}), or ({\em c}).

\np
Here surjective on the space of $\G$-invariants means (in the case of $\Y_{\svseq}$) that

$$
\H^{\N+\ell(\vseq)}(\Y_{\svseq},f^{*}_{\vseq}\L_{\ll})^{\G} \stackrel{f_{\vseq}^{*}}{\longleftarrow} 
\H^{\N+\ell(\vseq)}(\X^{k},\L_{\ll})^{\G}
$$

\np
is surjective.

\bpf
To reduce notation 
set $\M_{\useq}=\X^{j-1}\times\M_{i_m}\times\X^{k-j}$ and
let $\pi\colon\X^{k}\longrightarrow \M_{\useq}$ be the map $\pi=
(\id_{\X})^{j-1}\times\pi_{i_m}\times(\id_{\X})^{k-j}$. 
The fibre product diagram \eqref{eqn:bigfibredef}
relating $\Y_{\svseq}$, $\Y_{\suseq}$, $\X^{k}$, and $\M_{\useq}$ is

\begin{equation}\label{eqn:productdiagram}
\begin{array}{c}
\xymatrix{
\Y_{\svseq}\ar[r]^(0.5){f_{\vseq}}\ar[d]^{\pi_{\vseq}}\ar @{}[dr]|{\Box} & \X^{\subsmash{k}}\ar[d]^{\pi} \\
\Y_{\suseq}\ar[r]^(0.5){h}\ar @/^10pt/[u]^{\sigma_{\vseq}} & \M_{\subsmash{\useq}} \\
}
\end{array}
\end{equation}

\np
where 
$h=\pi\circ f_{\useq}$ and where we use 
$\pi_{\vseq}$ and $\sigma_{\vseq}$ in place of $\pi_{\vseq,\useq}$ and $\sigma_{\vseq,\useq}$ to reduce notation.

\np
Note that $\L_{\mm}$ is the Demazure reflection of $\L_{\ll}$ with respect to $\pi$.  By \S\ref{sec:demazure}
this means that we have natural isomorphisms 

\begin{equation}\label{eqn:demazure-isom}
\pi_{\vseq*}(f_{\vseq}^{*}\L_{\mm})\cong \R^1\pi_{\vseq*}(f_{\vseq}^{*}\L_{\ll})
\,\,\,\mbox{and}\,\,\,
\pi_{*}\L_{\mm}\cong \R^1\pi_{*}\L_{\ll}
\end{equation}

\np
valid on $\Y_{\suseq}$ and $\M_{\useq}$ respectively.
Diagram \eqref{eqn:productdiagram}, the Leray spectral sequences for $\L_{\ll}$ and $\L_{\mm}$ relative
to $\pi$ and $\pi_{\vseq}$, and the isomorphisms \eqref{eqn:demazure-isom} then give the commutative diagram of 
cohomology groups:

\begin{equation}\label{eqn:bigcohdiag}
\begin{array}{c}
\xymatrix{
\H^{\N+\ell(\vseq)}(\Y_{\svseq},f^{*}_{\vseq}\L_{\ll}) & \H^{\N+\ell(\vseq)}(\X^k,\L_{\ll})
\ar[l]_(0.52){f_{\vseq}^{*}}^{\eqref{eqn:top-map-hypoth}}\\
\H^{\N+\ell(\vseq)-1}(\Y_{\suseq},\R^1\pi_{\vseq*}f^{*}_{\vseq}\L_{\ll})\ar @{=}[u]^{\wr}_{\mbox{\tiny Leray}} & \H^{\N+\ell(\vseq)-1}(\M_{\useq},\R^1\pi_{*}\L_{\ll})\ar[l]_(0.51){h^{*}}\ar @{=}[u]^{\wr}_{\mbox{\tiny Leray}} \\
\H^{\N+\ell(\vseq)-1}(\Y_{\suseq},\pi_{\vseq*}f^{*}_{\vseq}\L_{\mm})\ar @{=}[u]^{\wr}_{\eqref{eqn:demazure-isom}} & \H^{\N+\ell(\vseq)-1}(\M_{\useq},\pi_{*}\L_{\mm})\ar[l]_(0.51){h^{*}}\ar @{=}[u]^{\wr}_{\eqref{eqn:demazure-isom}} \\
\H^{\N+\ell(\vseq)-1}(\Y_{\svseq},f^{*}_{\vseq}\L_{\mm})\ar @{=}[u]^{\wr}_{\mbox{\tiny Leray}} & \H^{\N+\ell(\vseq)-1}(\X^k,\L_{\mm})\ar[l]_(0.52){f_{\vseq}^{*}}\ar @{=}[u]^{\wr}_{\mbox{\tiny Leray}} \\
}
\end{array}
\end{equation}

\np
We conclude that the bottom pullback map $\H^{\N+\ell(\vseq)-1}(\Y_{\svseq},f_{\vseq}^{*}\L_{\mm})\stackrel{f_{\vseq}^{*}}{\longleftarrow}
\H^{\N+\ell(\vseq)-1}(\X^{k},f_{\vseq}^{*}\L_{\mm})$ is surjective, zero, or
surjective on the space of $\G$-invariants if \eqref{eqn:top-map-hypoth} is.

\np
On $\Y_{\svseq}$ we have the exact sequence of bundles:

\begin{equation}\label{eqn:restr-sequence}
0 \longrightarrow f^{*}_{\vseq}\L_{\mm}(-\Y_{\suseq})\longrightarrow f^{*}_{\vseq}\L_{\mm} \longrightarrow f^{*}_{\vseq}\L_{\mm}|_{\Y_{\suseq}} 
\longrightarrow 0,
\end{equation}

\np
where we consider $\Y_{\suseq}$ to be a divisor in $\Y_{\svseq}$ via the section $\sigma_{\vseq}$.
The degree of $f^{*}_{\vseq}\L_{\mm}(-\Y_{\suseq})$ is at least $-1$ on the 
fibres of $\pi_{\vseq}$ so the corresponding Leray spectral sequence gives
$$\H^{\N+\ell(\vseq)}(\Y_{\vseq},f^{*}_{\vseq}\L_{\mm}(-\Y_{\useq}))
=\H^{\N+\ell(\vseq)}\left({\Y_{\suseq},\pi_{_{\vseq}*}(f^{*}_{\vseq}\L_{\mm}(-\Y_{\suseq}))}\right)=0$$
where the second cohomology group above equals zero by reason of dimension: 
$$\N+\ell(\vseq)=\N+\ell(\useq)+1=\dim(\Y_{\suseq})+1.$$
The end of the long exact cohomology sequence associated to \eqref{eqn:restr-sequence} is therefore 

\begin{equation}\label{eqn:surj-seq}
\H^{\N+\ell(\vseq)-1}(\Y_{\vseq},f^{*}_{\vseq}\L_{\mm})\stackrel{\sigma_{\vseq}^{*}}{\longrightarrow}  
\H^{\N+\ell(\vseq)-1}(\Y_{\suseq},f^{*}_{\vseq}\L_{\mm}|_{\Y_{\suseq}})
\longrightarrow 0.
\end{equation}

\np
Since $\ell(\useq)=\ell(\vseq)-1$, $f_{\useq}=f_{\vseq}\circ \sigma_{\vseq}$, and all maps are 
$\G$-equivariant, we conclude that
the pullback map $f_{\useq}^{*}$, being the composite map
$$\H^{\N+\ell(\useq)}(\X^{k},\L_{\mm})\stackrel{f_{\vseq}^{*}}{\longrightarrow} 
\H^{\N+\ell(\useq)}(\Y_{\svseq},f^{*}_{\vseq}\L_{\mm})\stackrel{\sigma_{\vseq}^{*}}{\longrightarrow}
\H^{\N+\ell(\useq)}(\Y_{\suseq},f^{*}_{\vseq}\L_{\mm}|_{\Y_{\suseq}})=
\H^{\N+\ell(\useq)}(\Y_{\suseq},f^{*}_{\useq}\L_{\mm}),$$
is ({\em a}) surjective, ({\em b}) zero,
or ({\em c}) surjective on the space of $\G$-invariants, if the pullback map $f_{\vseq}^{*}$ 
in \eqref{eqn:top-map-hypoth} has the corresponding property ({\em a}), ({\em b}), or ({\em c}).  \epf

\np
{\bf Remark.} In part ({\em c}) of Lemma \ref{lem:induct} we can replace the statement about $\G$-invariants
with a statement about any isotypic component; the proof above goes through without change, 
we will only need the case of $\G$-invariants as part of the proof of Theorem I in \S\ref{sec:ThmI-proof} below.

\bpoint{Proof of Theorem \ref{thm:topsurjectivity} and variation} \label{proof-of-top-dim} 
For the rest of this section, we fix the following notation.
Let $w_1$,\ldots, $w_k$ and $\lambda_1$,\ldots, $\lambda_k$ be as in Theorem \ref{thm:topsurjectivity}.
For each $i=1$,\ldots, $k$ set $v_i:=w_i^{-1}\wo$ and $\lambda_i':=v_i^{-1}\cdot\lambda_i$.
Let $\vv_i$ be a reduced factorization of $v_i$ and let $\vseq=(\vv_1,\ldots, \vv_{k})$.
Finally, set $\ll=(\lambda_1,\ldots,\lambda_k)$ and $\ll'=(\lambda_1',\ldots,\lambda_k')$.

\np
{\em Proof of Theorem \ref{thm:topsurjectivity}.}
Since 
$$\dim(\Y_{\svseq})= \N + \sum_{i=1}^{k} \ell(v_i) =
\N + \sum_{i=1}^{k} (\N - \ell(w_i)) = 
k\N=\dim(\X^{k}),$$
Corollary \ref{cor:map-degree} implies that the degree of $f_{\vseq}\colon\Y_{\svseq}\longrightarrow \X^{k}$ 
is given by the intersection number ${\bigcap_{i=1}^{k}[\Omega_{w_i}]}$.  Therefore
the pullback map $\H^{k\N}(\Y_{\svseq},f_{\vseq}^{*}\L_{\ll'})\stackrel{f_{\vseq}^{*}}{\longleftarrow} \H^{k\N}(\X^{k}, \L_{\ll'})$ is a surjection if ${\bigcap_{i=1}^{k}[\Omega_{w_i}]}=1$
and is zero if ${\bigcap_{i=1}^{k}[\Omega_{w_i}]}=0$. Indeed, if 
${\bigcap_{i=1}^{k}[\Omega_{w_i}]}=1$ then $f_{\vseq}$ is a birational map between the smooth
varieties $\Y_{\svseq}$
and $\X^{k}$ in characteristic zero, 
and so the pullback map $\H^{j}(\Y_{\svseq},f_{\vseq}^{*}\L_{\ll'})\stackrel{f_{\vseq}^{*}}{\longleftarrow} \H^{j}(\X^{k}, \L_{\ll'})$ is an isomorphism in 
all degrees, and in particular is a surjection in degree $j=k\N$. 
On the other hand,
if ${\bigcap_{i=1}^{k}[\Omega_{w_i}]}=0$ then the image $\X_{\vseq}$ of $f_{\vseq}$ is subvariety of $\X^{k}$ of
dimension strictly less than $k\N$ and therefore the pullback map $f^{*}_{\vseq}$ in top cohomology, 
which factors through $\H^{k\N}(\X_{\vseq},\L_{\ll}|_{\X_{\vseq}})=0$, is the zero map.

\np
Consider a sequence 
$$
\vseq=:
\vseq^{0},
\vseq^{1},\ldots,
\vseq^{(k-1)\N}:= \underline{\emptyset}=(\emptyset,\ldots,\emptyset)
$$
of sequences of words which reduces $\vseq$ to the empty sequence, and where at each step $\vseq^{j+1}$ is
obtained by dropping a simple reflection from the right of a single member of $\vseq^{j}$.  Set 
$\underline{\lambda}^{j}=(\vseq^{j})^{-1}\cdot\underline{\lambda}$ 
where (by slight abuse of notation) $\vseq^{j}$ is considered as an element of $\W^{k}$ and the action is componentwise.
Note that $\ll^{0}=\ll'$ and $\ll^{(k-1)\N}=\ll$.
The construction of $\vseq^{j}$ and $\underline{\lambda}^{j}$ implies that the degree of 
$\L_{\underline{\lambda}^{j}}$ is negative on the fibres of the $\PP^1$-fibration 
$\pi_{\vseq^{j},\vseq^{j+1}}\colon \Y_{\svseq^{j}}\longrightarrow\Y_{\svseq^{j+1}}$.

\np
Applying Lemma \ref{lem:induct} to the pairs $(\vseq^{j},\vseq^{j+1})$ for $j=0$,\ldots, $(k-1)\N-1$ we conclude
that 

\begin{equation}\label{eqn:inductbottom}
\H^\N(\Y_{\underline{\emptyset}}, f^{*}_{\underline{\emptyset}}\L_{\ll}) 
\stackrel{f_{\underline{\emptyset}}^{*}}{\longleftarrow} \H^\N(\X^{k}, \L_{\ll})
\end{equation}

\np
is surjective if ${\bigcap_{i=1}^{k}[\Omega_{w_i}]}=1$ and zero if ${\bigcap_{i=1}^{k}[\Omega_{w_i}]}=0$.
By construction $f_{\underline{\emptyset}}\colon \Y_{\underline{\emptyset}}=\X \longrightarrow \X^{k}$ is the diagonal embedding
of $\X$ into $\X^{k}$ and the pullback map \eqref{eqn:inductbottom}
is the cup product map. This proves Theorem \ref{thm:topsurjectivity} and completes the
proof of Theorem III. \epf

\np
We record a statement that will be used in the proof of Theorem I below.

\tpoint{Proposition}\label{prop:variation}
If the pullback map 
$$\H^{k\N}(\Y_{\svseq},f_{\vseq}^{*}\L_{\ll'})\stackrel{f_{\vseq}^{*}}{\longleftarrow} 
\H^{k\N}(\X^{k}, \L_{\ll'})$$ 
is surjective on the space of $\G$-invariants then the cup product map 
\eqref{eqn:symmetriccup} is surjective.  

\np
\bpf We repeat the inductive reduction in the proof of Theorem \ref{thm:topsurjectivity}  above with part
({\em c}) of Lemma \ref{lem:induct} in place of parts ({\em a}) and ({\em b}). As a result we conclude that 
the cup product map \eqref{eqn:symmetriccup} is surjective on the space of $\G$-invariants.  Since
$\H^{\N}(\X,\K_{\X})$ is the trivial $\G$-module we conclude that 
\eqref{eqn:symmetriccup} is surjective. \epf

\section{Proof of Theorem I and corollaries}\label{sec:ThmI-proof}

\np
In this section we use Theorem III and Proposition \ref{prop:variation} to prove Theorem I.
The proof that \eqref{eqn:liningup} is necessary for the surjectivity of the cup product map 
appears in \S\ref{sec:neces-cond} and the proof that \eqref{eqn:liningup} is 
sufficient appears in \S\ref{sec:suff-cond}.

\bpoint{Proof that $\invset_{w}=\sqcup_{i=1}^{k} \invset_{w_i}$ is a necessary condition for surjectivity}
\label{sec:neces-cond}
We assume the  notation of \S\ref{sec:setup}, and set $\mu_i=w_i\cdot \lambda_i$ for $i=1\ldots k$, and 
$\mu=w\cdot \lambda$.  By assumption the weights $\mu_1$,\ldots, $\mu_k$,  and $\mu$ are dominant. 
By the Borel-Weil-Bott theorem each $\H^{\ell(w_i)}(\X,\L_{\lambda_i})=\V_{\smu_i}^{*}$ 
and $\H^{d}(\X,\L_{\lambda})=\V_{\smu}^{*}$.

\np
Since $w_i^{-1}\mu_i = w_i^{-1}\cdot\mu_i - w_i^{-1}\cdot 0$ and $w^{-1}\mu_i = w^{-1}\cdot\mu_i - w^{-1}\cdot 0$, 
we have

$$
\sum_{i=1}^{k} w_i^{-1}\mu_i - w^{-1}\mu 
= \left( \sum_{i=1}^{k} w_i^{-1}\cdot\mu_i - w^{-1}\cdot\mu \right)
- \left( \sum_{i=1}^{k} w_i^{-1}\cdot0 - w^{-1}\cdot0 \right).
$$

\np
Furthermore $\sum_{i=1}^{k} w_i^{-1}\cdot\mu_i - w^{-1}\cdot\mu = \sum \lambda_i - \lambda=0$  and so the equation
above becomes

\begin{equation}\label{eqn:rholinesup}
\sum_{i=1}^{k} w_i^{-1}\mu_i - w^{-1}\mu 
=
- \left( \sum_{i=1}^{k} w_i^{-1}\cdot0 - w^{-1}\cdot0 \right).
\end{equation}

\np
If the cup product map 
$\H^{\ell(w_1)}(\X,\L_{\lambda_1})\otimes\cdots\otimes \H^{\ell(w_k)}(\X,\L_{\lambda_k})
\stackrel{\cup}{\longrightarrow}
\H^d(\X,\L_{\lambda})$
is surjective, then (after dualizing)
$\V_{\smu}$ must be a component of the tensor product $\V_{\smu_1}\otimes\cdots\otimes\V_{\smu_k}$ and
by Theorem III({\em b}), the intersection 
${\cap_{i=1}^{k}[\Omega_{w_i}]\cdot[\X_{w}]\neq 0}$ in $\H^{*}(\X,\ZZ)$; we may therefore apply Theorem
\ref{thm:hammer}.

\np
By Theorem \ref{thm:hammer}({\em a}) the left hand side of \eqref{eqn:rholinesup} belongs to $\Span_{\ZZ\geq0}\delpos$
and by part ({\em c}) of the same theorem the right hand side belongs to $\Span_{\ZZ\leq 0}\delpos$.
We conclude that both sides are zero and so by Theorem \ref{thm:hammer}({\em d}) that 
$\invset_{w}=\sqcup_{i=1}^{k} \invset_{w_i}$. \epf

\np
{\bf Remark.} In the first half of the argument above the hypothesis that the cup product map is surjective 
was used, along with Theorem III, to
conclude that $\V_{\smu}$ is a component of $\V_{\smu_1}\otimes\cdots\otimes\V_{\smu_k}$ 
and ${\cap_{i=1}^{k}[\Omega_{w_i}]\cdot[\X_{w}]\neq 0}$.  If, on the other hand, we assume the latter two conditions
then the second half of the argument still applies to give $\invset_{w}=\sqcup_{i=1}^{k} \invset_{w_i}$. 
We will use this observation in Corollary \ref{cor:iscomponent} below.

\bpoint{Setup for the proof of sufficiency}
For convenience, we collect some of the consequences of condition \eqref{eqn:liningup} in its symmetric form
which have effectively appeared in previous arguments, and which we will use in the proof of sufficiency.

\tpoint{Proposition }
\label{prop:consequences}
Suppose that $w_1$, \ldots, $w_k$ are elements of the Weyl group such that 
$\delpos = \bigsqcup_{i=1}^{k} \invset_{w_i}.$

\np
{\em Combinatorial Consequences:}

\begin{enumerate}
\item $\sum_{i=1}^{k} w_i^{-1}\cdot 0 = -2\rho$.
\item Suppose that $\lambda_1$,\ldots, $\lambda_{k}$ are weights such that $\sum \lambda_i=-2\rho$, and set
$\mu_i=w_i\cdot\lambda_i$ for $i=1$,\ldots, $k$.  Then $\sum_{i=1}^{k} w_i^{-1}\mu_i=0$.
\end{enumerate}

\np
{\em Geometric Consequences:}
For each  $i=1$,\ldots, $k$, 
let $v_i=w_i^{-1}\wo$ and let $\vv_i$ be a word which is a reduced factorization of
$v_i$.  We set $\vseq=(\vv_1,\ldots, \vv_{k})$ and construct as usual the variety
$\Y_{\svseq}$ and the map $f_{\vseq}\colon\Y_{\svseq}\longrightarrow \X^{k}$.

\np
Then

\begin{enumerate}
\setcounter{enumi}{2}
\item $\deg(f_{\vseq})\neq 0$.
\item The weight of  the relative canonical bundle $\K_{\Y_{\svseq}/\K_{\X^{k}}}$ at $p_{\vseq}$ is zero.
\end{enumerate}

\np
\bpf
Part ({\em a}) is immediate from the condition $\delpos=\sqcup_{i=1}^{k}\invset_{w_i}$ and formula \eqref{eqn:rootsum}.
Part ({\em b}) reverses the argument used to arrive at \eqref{eqn:rholinesup} in \S\ref{sec:neces-cond}:

$$\sum_{i=1}^{k} \mu^{-1}\mu_i = 
\left(\sum_{i=1}^{k} w_i^{-1}\cdot\mu_i\right) -\left(\sum_{i=1}^{k} w_i^{-1}\cdot 0\right) 
=\sum_{i=1}^{k} \lambda_i - (-2\rho) = 0.  $$

\np
Part ({\em c}) is Corollary \ref{cor:map-degree} combined with Lemma \ref{lem:nonzero}.
Part ({\em d}) is the symmetric version of the computation in the proof of Theorem \ref{thm:hammer}({\em d}): the
weight of  the relative canonical bundle $\K_{\Y_{\svseq}/\X^{k}}$ at $p_{\vseq}$ is $\sum_{i=1}^{k} w_i^{-1}\cdot 0 + 2\rho$ which is zero 
by part ({\em a}). 
\epf

\bpoint{Proof that $\invset_{w}=\sqcup_{i=1}^{k} \invset_{w_i}$ is a sufficient condition for surjectivity}
\label{sec:suff-cond}
Consider the symmetric version of the problem as in \S\ref{sec:nonsymmetric-symmetric}.  It suffices to 
show the surjectivity of a cup product map 

\begin{equation}\label{eqn:symmetric-cupprod}
\H^{\ell(w_1)}(\X,\L_{\lambda_1})\otimes\cdots\otimes \H^{\ell(w_{k})}(\X,\L_{\lambda_{k}})
\stackrel{\cup}{\longrightarrow}
\H^{\N}(\X,\Ksh_{\X})
\end{equation}

\np
where $w_1$, \ldots, $w_{k}$ are elements of the Weyl group such that 
$\sum \ell(w_i)=\N$; $\lambda_1$, \ldots, $\lambda_{k}$  are weights such that 
$w_i\cdot \lambda_i\in\Lambda^{+}$ for $i=1$,\ldots,$k$ and $\sum \lambda_i=-2\rho$. 
After this reduction condition \eqref{eqn:liningup} becomes 
$
\delpos = \sqcup_{i=1}^{k} \invset_{w_i}.
$
We recall the notation from \S\ref{proof-of-top-dim}: $v_i:=w_i^{-1}\wo$, 
$\lambda_i':=v_i^{-1}\cdot\lambda_i$,
$\vv_i$ is a reduced factorization of $v_i$, $\vseq=(\vv_1,\ldots,\vv_{k})$, 
and $\ll'=(\lambda_1',\ldots,\lambda_{k}')$.

\np
By Proposition \ref{prop:variation} to show the surjectivity of \eqref{eqn:symmetric-cupprod} it is enough
to show that the pullback map 

\begin{equation}\label{eqn:substep-pullback}
\H^{k\N}(\Y_{\svseq},f_{\vseq}^{*}\L_{\ll'})^{\G}\stackrel{f_{\vseq}^{*}}{\longleftarrow} \H^{k\N}(\X^{k}, \L_{\ll'})^{\G}
\end{equation}

\np
on the space of $\G$-invariants is surjective.  We will show that both spaces of $\G$-invariants are 
one-dimensional, and that the induced map is an isomorphism. Note that by Proposition \ref{prop:consequences}({\em c}) 
$\deg(f_{\vseq})\neq 0$ and so $f_{\vseq}$ is surjective.

\medskip

\np
The pullback map on top cohomology is Serre dual to the trace map:

$$
\H^{0}(\Y_{\svseq}, (f_{\vseq}^{*}\L_{\ll'})^{*}\otimes \K_{\Y_{\vseq}}) =
\H^{0}\left({\Y_{\svseq}, f_{\vseq}^{*}(\L_{\ll'}^{*}\otimes \K_{\X^{k}})\otimes \K_{\Y_{\vseq}/\X^{k}}}\right) 
\stackrel{\Tr_{f_{\vseq}}}{\longrightarrow}
\H^{0}(\X^{k},\L_{\ll'}^{*}\otimes \K_{\X^{k}}).
$$

\np
Let $s\in\H^{0}(\Y_{\svseq},\K_{\Y_{\svseq}/\X^{k}})^{\G}$ be the nonzero $\G$-invariant
section giving the map $f_{\vseq}^{*}\K_{\X^{k}}\longrightarrow \K_{\Y_{\svseq}}$ induced by $f_{\vseq}$, and let
$\E$ be the divisor of $s=0$.
The composition 

\noindent
\begin{minipage}{\textwidth}
\noindent
\begin{eqnarray*}
\H^{0}(\X^{k},\L_{\ll'}^{*}\otimes \K_{\X^{k}}) & \stackrel{f_{\vseq}^{*}}{\longrightarrow} & 
\H^{0}\left({\Y_{\svseq}, f_{\vseq}^{*}(\L_{\ll'}^{*}\otimes \K_{\X^{k}})}\right) \\
& \stackrel{\cdot\E}{\longrightarrow} & 
\H^{0}\left({\Y_{\svseq}, f_{\vseq}^{*}(\L_{\ll'}^{*}\otimes \K_{\X^{k}})\otimes \K_{\Y_{\vseq}/\X^{k}}}\right)  
 \stackrel{\Tr_{f_{\vseq}}}{\longrightarrow}  
\H^{0}(\X^{k},\L_{\ll'}^{*}\otimes \K_{\X^{k}}) \\
\end{eqnarray*}
\end{minipage}

\np
of pullback, multiplication by $\E$,  and the trace map is
multiplication by $\deg(f_{\vseq})$, which is nonzero.    This gives us the inequality 

\begin{equation}\label{eqn:surj-inequality}
\dim
\H^{0}(\X^{k},\L_{\ll'}^{*}\otimes \K_{\X^{k}})^{\G} \leq 
\dim 
\H^{0}\left({\Y_{\svseq}, f_{\vseq}^{*}(\L_{\ll'}^{*}\otimes \K_{\X^{k}})\otimes \K_{\Y_{\vseq}/\X^{k}}}\right)^{\G}
\end{equation}

\np
and shows that in order to prove that the trace map induces an isomorphism on $\G$-invariants it
is sufficient to prove that we have equality of dimensions in \eqref{eqn:surj-inequality}.

\np
Set $\mu_i=w_i\cdot\lambda_i=\wo\cdot\lambda_i'$ for $i=1$,\ldots, $k$. By the Borel-Weil-Bott Theorem we have 
$\H^{k\N}(\X^{k},\L_{\ll'}) = \V_{\smu_1}^{*}\otimes\cdots\otimes\V_{\smu_k}^{*}$ and so (by Serre duality) 
$\H^{0}(\X^{k}, \L_{\ll'}^{*}\otimes\K_{\X^{k}})=\V_{\smu_1}\otimes\cdots\otimes\V_{\smu_k}$.
Now set $\nu_i=-\wo\mu_i$ 
for $i=1$,\ldots, $k$ so that $\V_{\snu_i}=\V_{\smu_i}^{*}$ and let $\nn=(\nu_1,\ldots, \nu_k)$.   
By the calculation 
$$\SS(\lambda_i')=-\lambda_i'-2\rho =-\wo\cdot\mu_i -2\rho = -(\wo\mu_i-2\rho)-2\rho=-\wo\mu_i$$
in each coordinate factor (as in \S\ref{sec:serredual}) we conclude that $\L_{\ll'}^{*}\otimes\K_{\X^{k}}=\L_{\nn}$.

\np
The weight of $\L_{\nn}$ at $q:=f_{\vseq}(p_{\vseq})=(v_1,\ldots,v_k)$ is 

$$-\sum_{i=1}^{k} v_i\nu_i = \sum_{i=1}^{k} (w_i^{-1}\wo)(\wo\mu_i)
=\sum_{i=1}^{k} w_i^{-1}\mu_i
\stackrel{\mbox{\tiny\ref{prop:consequences}({\em b})}}{=}0.$$

\np
Since by Proposition \ref{prop:consequences}({\em d}) the weight of  $\K_{\Y_{\svseq}/\X^{k}}$ at $p_{\vseq}$ is zero, 
the weight of $f_{\vseq}^{*}\L_{\nn}\otimes \K_{\Y_{\svseq}/\X^{k}}$ at $p_{\vseq}$ in $\Y_{\svseq}$ is also 
zero and hence 
$\dim \H^{0}(\Y_{\svseq},f_{\svseq}^{*}\L_{\nn}\otimes\K_{\Y_{\svseq}/\X^{k}})^{\G}\leq 1$ by Lemma \ref{lem:key}({\em c}).
On the other hand, Lemma \ref{lem:PRVsymmetric} implies
that $(\V_{\smu_1}\otimes\cdots\otimes\V_{\smu_k})^{\G}\neq 0$ so we conclude that 
$\dim \H^{0}(\X^{k},\L_{\nn})^{\G}\geq 1$.  This gives us 

$$ 1\leq \dim \H^{0}(\X^{k},\L_{\nn})^{\G} \leq 
\dim \H^{0}({\Y_{\svseq},f_{\vseq}^{*}\L_{\nn}\otimes\K_{\Y_{\svseq}/\X^{k}}})^{\G}\leq 1.$$
Therefore the inequality in \eqref{eqn:surj-inequality} is an equality, and the cup product map 
in \eqref{eqn:symmetric-cupprod} is surjective.  
\epf

\bpoint{Corollaries of Theorem I and its proof}
\label{sec:ThmIcor}

\tpoint{Corollary}\label{cor:H0xHk} The cup product map 
$\H^{0}(\X,\L_{\lambda_1})\otimes\H^{d}(\X,\L_{\lambda_2}) \longrightarrow \H^{d}(\X,\L_{\lambda_1}\otimes\L_{\lambda_2})$
is surjective whenever both sides are nonzero.

\np
\bpf
If $w_2$ is the element of the Weyl group so that $w_2\cdot\lambda_2\geq 0$ then the conditions that $\lambda_1$ is 
dominant and that $\L_{\lambda_1+\lambda_2}$ has cohomology in the same degree $d$ as $\L_{\lambda_2}$ imply that 
$w_2\cdot(\lambda_1+\lambda_2)\geq 0$, and so the corollary follows from Theorem I and the obvious
statement that $\invset_{w_2}=\invset_{w_2}\sqcup \invset_{e}$.
\epf

\tpoint{Corollary} {\em (Compatibility with Leray Spectral Sequence).}\label{cor:compat}
Suppose that $\lambda_1$, $\lambda_2$, and
$\lambda=\lambda_1+\lambda_2$ are regular weights 
and that the cup product map 

$$\H^{d_1}(\X,\L_{\lambda_1})\otimes\H^{d_2}(\X,\L_{\lambda_2})\stackrel{\cup}{\longrightarrow}
\H^{d}(\X,\L_{\lambda})$$

\np
is nonzero.   Let $\P$ be any parabolic subgroup of $\G$ containing $\B$, 
and $\pi\colon\X\longrightarrow \M:=\G/\P$ be the corresponding projection.  
Then the cup product map on $\X$ factors as a composition

\np
\hspace{-0.4cm}
$
\xymatrix{
\H^{d_1}(\X,\L_{\lambda_1})\otimes\H^{d_2}(\X,\L_{\lambda_2})\ar@{=}[d]\ar[rr]^{\bigcup_{\X}} & & \H^{d}(\X,\L_{\lambda})\ar@{=}[d] \\
\H^{d_1-i}(\M,\R^{i}_{\pi*}\L_{\lambda_1})\otimes 
\H^{d_2-j}(\M,\R^{j}_{\pi*}\L_{\lambda_2}) \ar[r]^(.55){\bigcup_{\M}} &
\H^{d-i-j}(\M,\R^{i}_{\pi*}\L_{\lambda_1}\otimes\R^{j}_{\pi*}\L_{\lambda_2})\ar[r]^(.6){\bigcup_{\pi}} & 
\H^{d-i-j}(\M,\R^{i+j}_{\pi*}\L_{\lambda})
}
$

\medskip

\np
of the cup product on $\M$ 
followed by the map induced on cohomology by the relative cup product map
$\R^{i}_{\pi*}\L_{\lambda_1}\otimes\R^{j}_{\pi*}\L_{\lambda_2}\longrightarrow \R^{i+j}_{\pi*}\L_{\lambda}$ on
the fibres of $\pi$.   A similar statement holds for the cup product of an arbitrary number of factors.

\np
\bpf
The factorization statement amounts to a numerical condition on the cohomology degrees of the line bundles on
the fibres of $\pi$ ensuring that the cup product map is computed by the map on $\E_2$-terms 
of the Leray spectral sequence.  This numerical condition 
is immediately implied by \eqref{eqn:liningup}.  
We explain this in more detail below.

\np
Set $\L=\L_{\lambda_1}\bt\L_{\lambda_2}$ on $\X\times\X$ and 
consider the following factorization of the diagonal map $\delta_{\X}\colon\X\hookrightarrow\X\times\X$:

\begin{equation}\label{eqn:factordiagram}
\begin{array}{c}
\xymatrix{
\X\ar[d]^{\pi} \ar@{^{(}->}[r]^(0.4){s} & \X\times_{\M}\X \ar@{^{(}->}[r]^(0.55){t}\ar[d]^{\psi} \ar@{}[dr]|\Box & \X\times\X\ar[d]^{\pi\times\pi} \\
\M\ar@{=}[r] & \M \ar@{^{(}->}[r]^{\delta_{\M}} & \M\times \M 
}
\end{array}.
\end{equation}

\np
The cup product map then factors as

\begin{equation}\label{eqn:factorizationsequence}
\,\,\,\,
\H^{d}(\X,\L_{\lambda})
\stackrel{s^{*}}{\longleftarrow} 
\H^{d}(\X\times_{\M}\X,t^{*}\L) 
\stackrel{t^{*}}{\longleftarrow}
\H^{d_1}(\X,\L_{\lambda_1})\otimes\H^{d_2}(\X,\L_{\lambda_2})=\H^{d}(\X\times\X,\L)
\end{equation}

\np
and we claim that \eqref{eqn:factorizationsequence} induces the factorization claimed above.

\np
By the Borel-Weil-Bott theorem applied to the fibres of $\pi$, for each of the line bundles  $\L_{\lambda_1}$,
$\L_{\lambda_2}$, and $\L_{\lambda}$ 
there is precisely one degree for which the higher direct image sheaf is nonzero.
Suppose $i$ is the degree such that $\R^{i}_{\pi*}\L_{\lambda_1}\neq 0$, 
$j$ is the degree such that $\R^{j}_{\pi*}\L_{\lambda_2}\neq 0$, and 
$k$ is the degree such that $\R^{k}_{\pi*}\L_{\lambda}\neq 0$.  
The Leray spectral sequence for the cohomology of these
bundles degenerates at the $\E_2$ term and we have the isomorphisms  
$\H^{d_1}(\X,\L_{\lambda_1})=\H^{d_1-i}(\M,\R^{i}_{\pi*}\L_{\lambda_1})$,
$\H^{d_2}(\X,\L_{\lambda_2})=\H^{d_2-j}(\M,\R^{j}_{\pi*}\L_{\lambda_2})$, and 
$\H^{d}(\X,\L_{\lambda_1})=\H^{d-k}(\M,\R^{k}_{\pi*}\L_{\lambda})$.

\np
Since $\R^{i+j}_{\pi\times\pi*}\L = \R^{i}_{\pi*}\L_{\lambda_1}\bt\R^{j}_{\pi*}\L_{\lambda_2}$ is a vector
bundle on $\M\times\M$, the theorem on cohomology and base change gives us 
$\R^{i+j}_{\psi*}t^{*}\L=\delta_{\M}^{*}(\R^{i}_{\pi*}\L_{\lambda_1}\bt\R^{j}_{\pi*}\L_{\lambda_2})=
\R^{i}_{\pi*}\L_{\lambda_1}\otimes\R^{j}_{\pi*}\L_{\lambda_2}$ on $\M$ and therefore we have 
$\H^{d}(\X\times_{\M}\X,t^{*}\L)=\H^{d-i-j}(\M, \R^{i}_{\pi*}\L_{\lambda_1}\otimes\R^{j}_{\pi*}\L_{\lambda_2})$.
The Leray spectral sequences for $\L$ and $t^{*}\L$ with respect to $\psi$ and $\pi\times\pi$ also degenerate
at the $\E_2$-terms and have nonzero terms in the same degree. The discussion in \S\ref{sec:E2compute-begin} 
implies that the map on $\E_2$-terms 
computes the pullback map $t^{*}$.  Therefore $t^{*}$ in \eqref{eqn:factorizationsequence} is equal to the map

$$
\H^{d-i-j}(\M, \R^{i}_{\pi*}\L_{\lambda_1}\otimes\R^{j}_{\pi*}\L_{\lambda_2}) \stackrel{\delta_{\M}^{*}}{\longleftarrow}
\H^{d_1-i}(\M,\R^{i}_{\pi*}\L_{\lambda_1})\otimes\H^{d_2-j}(\M,\R^{j}_{\pi*}\L_{\lambda_2})
$$

\np
which shows that $t^{*}$ is the first part of the factorization claimed.  

\np
We now study $s^{*}$.   The map $s$ includes $\X$ as the relative diagonal of $\X\times_{\M}\X$ over $\M$. It follows
that $s^{*}$ induces the relative cup product map 
on the higher direct image sheaves of $t^{*}\L$ and $\L_{\lambda}$.
Therefore the map associated to $s^{*}$ on the $\E_2$-terms of the Leray spectral sequences for $t^{*}\L$ 
and $\L_{\lambda}$ is given by the relative cup product map

$$
\H^{d-i-j}(\M,\R^{i+j}\L_{\lambda})=
\H^{d-i-j}(\M,\R^{i+j}(\L_{\lambda_1}\otimes\L_{\lambda_2}))
\stackrel{\bigcup_{\pi}}{\longleftarrow}\H^{d-i-j}(\M,
\R^{i}_{\pi*}\L_{\lambda_1}\otimes\R^{j}_{\pi*}\L_{\lambda_2}).$$

\np
All that is needed to demonstrate the factorization claimed is to demonstrate the condition 
$k=i+j$ which ensures the map on the associated graded pieces in the $\E_2$-terms 
agrees with the global map on the cohomology groups (c.f.\ \S\ref{sec:E2compute-begin}).

\np
Suppose that $w_1$, $w_2$, and $w$ are
the elements of the Weyl group such that $w_1\cdot\lambda_1$, $w_2\cdot\lambda_2$, and 
$w\cdot\lambda$ are dominant.  Then $k=\#(\invset_{w}\cap-\del_{\P})$, $i=\#(\invset_{w_1}\cap-\del_{\P})$, and
$j=\#(\invset_{w_2}\cap-\del_{\P})$ where the symbol $\#$ indicates the cardinality of a set.  
The condition $k=i+j$ guaranteeing the factorization thus amounts to the condition

\begin{equation}\label{eqn:parabolicliningup}
\#(\invset_{w}\cap-\del_{\P})= \#(\invset_{w_1}\cap -\del_{\P})+\#(\invset_{w_2}\cap -\del_{\P}).
\end{equation}

\np
Since the original cup product map was assumed surjective we must have $\invset_{w}=\invset_{w_1}\sqcup\invset_{w_2}$
by Theorem I; this immediately implies that \eqref{eqn:parabolicliningup} holds.
\epf

\tpoint{Corollary} \label{cor:mult-one-again}
Suppose that $w_1$,\ldots, $w_k$ are elements of the Weyl group such that 
$\delpos=\sqcup_{i=1}^{k} \invset_{w_i}$, and that $\mu_1$,\ldots, $\mu_k$ are dominant weights satisfying the
condition $\sum_{i=1}^{k} w_i^{-1}\mu_i=0$.  Then 
$\dim (\V_{\smu_1}\otimes\cdots\otimes\V_{\smu_k})^{\G}=1$.

\np
\bpf 
Set $\lambda_i=w_i^{-1}\cdot\mu_i$ for $i=1$,\ldots, $k$. Then $\sum \lambda_i=-2\rho$ and we have a cup
product problem as in \eqref{eqn:symmetric-cupprod}.  
As part of the proof of Theorem I in \S\ref{sec:suff-cond} it was established that
$\dim (\V_{\smu_1}\otimes\cdots\otimes\V_{\smu_k})^{\G}=1$.  Alternatively, the corollary is simply 
Theorem \ref{thm:hammer}({\em b}) applied in symmetric form, with Lemma \ref{lem:nonzero} used to ensure that
the hypotheses of the theorem are satisfied.
\epf

\tpoint{Corollary} \label{cor:iscomponent}
Suppose that we have a cup product map 
$$
\H^{\ell(w_1)}(\X,\L_{\lambda_1})\otimes\cdots\otimes \H^{\ell(w_k)}(\X,\L_{\lambda_k})
\stackrel{\cup}{\longrightarrow}
\H^d(\X,\L_{\lambda})
$$

\np
and as above Weyl group elements $w_1$,\ldots, $w_k$, and $w$ such that $\mu_i:=w_i\cdot\lambda_i$, 
$i=1$,\ldots, $k$, and $\mu:=w\cdot\mu$ are dominant weights.  Then if 
${\cap_{i=1}^{k}[\Omega_{w_i}]\cdot[\X_{w}]\neq 0}$ the cup product map is surjective if and only if
$\V_{\mu}$ is a component of $\V_{\mu_1}\otimes\cdots\otimes\V_{\mu_k}$.

\np
\bpf If $\V_{\mu}$ is not a component of the tensor product the map is clearly not surjective.  Conversely,
if $\V_{\mu}$ is a component, the assumption on the intersection number and the argument in \S\ref{sec:neces-cond}
for the necessity of condition \eqref{eqn:liningup} show that $\invset_{w}=\sqcup_{i=1}^{k} \invset_{w_i}$, 
and therefore we conclude that the map is surjective by the sufficiency of condition \eqref{eqn:liningup}.  \epf

\np
The following example illustrates Corollary \ref{cor:iscomponent} and provides an example which shows that 
condition \eqref{eqn:liningup} is not necessary in order to have a cup product problem for which both
sides are nonzero.

\tpoint{Example} \label{ex:nolineup}
Let $\G=\Sl_{6}$ and $w_1=w_2=s_2s_4s_3$. For any integers $a_i,b_i\geq 0$ ($i=1,2$) 
set $\mu_i=(0,a_i,0,b_i,0)$ and $\lambda_i=w_i^{-1}\cdot\mu_i=(a_i+1,b_i+1,-4-a_i-b_i,a_i+1,b_i+1)$.
(The weights are written in terms of the fundamental weights of $\Sl_{6}$.)
Finally,
let $w=s_1s_3s_5s_2s_4s_3$ and set $\mu=w\cdot(\lambda_1+\lambda_2)=(0,a_1+a_2+1,0,b_1+b_2+1,0)\in\Lambda^{+}$.
We therefore get a cup product problem:

$$\H^{3}(\X,\L_{\lambda_1})\otimes\H^{3}(\X,\L_{\lambda_2})\stackrel{\bigcup}{\longrightarrow}\H^{6}(\X,\L_{\lambda_1+\lambda_2}).$$

\np
This cup product cannot be surjective by Theorem I since $\invset_{w_1}=\invset_{w_2}$; alternatively the map 
cannot be surjective since $\V_{\mu}$ is clearly not a component of $\V_{\mu_1}\otimes\V_{\mu_2}$.  
The intersection number $([\Omega_{w_1}]\cap[\Omega_{w_2}])\cdot\X_{w}$ is two.

\tpoint{Corollary} \label{cor:unionofinvsets}
If $\delpos=\sqcup_{i=1}^{k} \invset_{w_i}$ then for any subset $\I\subseteq\{1\ldots, k\}$
there is an element $w$ of the Weyl group such that $\invset_{w}=\sqcup_{i\in \I} \invset_{w_i}$.  

\np
\bpf
Let $\lambda_i=w_i^{-1}\cdot 0$ so that we get a cup product problem as in \eqref{eqn:symmetric-cupprod}. (Here
each $\H^{\ell(w_i)}(\X,\L_{\lambda_i})$ is the trivial $\G$-module).  By Theorem I and the assumption
on $w_1$,\ldots, $w_k$ this cup product is surjective.  It can be factored by first taking the cup
product of any subset $\I\subseteq\{1,\ldots, k\}$ of the factors and
the resulting cup product problem must also be nonzero since the larger problem is.  Hence by Theorem I there
is a $w\in \W$ with $w\cdot(\sum_{i\in \I} \lambda_i)\in\Lambda^{+}$ 
and such that $\invset_{w}=\sqcup_{i\in \I}\invset_{w_i}$.
\epf

\bpoint{Comments}

\np
{\bf 1.}
Corollary \ref{cor:unionofinvsets} 
can also be proved independently of any of the constructions in this paper by using a similar
argument in nilpotent cohomology.  
We are grateful to Olivier Mathieu for pointing this out to us.

\np
{\bf 2.} Using the result of Corollary \ref{cor:unionofinvsets} and induction, to prove Theorem I it is 
sufficient to prove it in the case $k=2$ of the cup product of two cohomology groups into a third.  
We have chosen to develop the description of the varieties $\Y_{\svseq}$ for arbitrary $k$ partly since this 
is the natural generality of the construction, partly because it makes no difference in our proofs, but also
because some of the applications (e.g., the multiplicity bounds) do not follow by induction.
Note that by the methods of this paper, even to prove the case $k=2$ of the cup product 
it would be necessary to consider the case of the cup product of three factors into $\H^{\N}(\X,\K_{\X})$, and hence
we would need the construction of $\Y_{\svseq}$ for three factors.

\np
{\bf 3.} As Example \ref{ex:nolineup} shows, the natural numerical condition $\ell(w_1)+\ell(w_2)=\ell(w)$
does not imply condition \eqref{eqn:liningup} even if there is a nontrivial cup product problem corresponding to
$w_1$, $w_2$, and $w$.  On the other hand, Condition \eqref{eqn:parabolicliningup} imposes further necessary 
numerical conditions for \eqref{eqn:liningup}.  Namely, 

\begin{equation}\label{eqn:minrepsadd}
\ell(w_1^{\P})+\ell(w_2^{\P}) = \ell(w^{\P}) \,\,\mbox{for every parabolic subgroup $\P\supseteq\B$ of $\G$},
\end{equation}
where $w_1^{\P}$, $w_2^{\P}$, $w^{\P}$ denote the minimal length representatives in $w_1\W_{\P}$, $w_2\W_{\P}$, 
$w\W_{\P}$.
In the case when $\G=\Sl_{n+1}$ one can show that condition \eqref{eqn:minrepsadd} is sufficient for 
\eqref{eqn:liningup}.
The simple inductive argument relies on the fact 
that if $\G=\Sl_{n+1}$ it is possible to assign a parabolic $\P_{\!\alpha}\supset\B$ to every root 
$\alpha\in\delpos$ in such a way that $-\alpha$ is a root of $\P_{\!\alpha}$ but not a root of any proper
parabolic subgroup of $\P_{\!\alpha}$ containing $\B$.
We do not know if \eqref{eqn:minrepsadd} is sufficient to imply \eqref{eqn:liningup} for general $\G$.

\np
{\bf 4.} 
Corollary \ref{cor:compat} establishes the following Factorization Property: any nonzero cup product map on $\X$
factors as a cup product on $\G/\P$ and fibres of $\pi\colon\X\longrightarrow\G/\P$ for all $\P\supset\B$.  
We know of no apriori reason why this should hold.   
The Factorization Property is equivalent to 
\eqref{eqn:parabolicliningup} holding for all $\P\supset\B$ which is equivalent to \eqref{eqn:minrepsadd}.
Hence, in the case $\G=\Sl_{n+1}$ the Factorization Property is equivalent to \eqref{eqn:liningup}.

\section{Cohomological components and proof of Theorem II}
\label{sec:cohomologicalcomponents}

\bpoint{Conditions on components of tensor products}
\label{sec:componentdefs}

\np
We begin by introducing two relevant conditions.  We also recall the notion of generalized PRV component
from \S\ref{sec:PRVdiscuss} for convenience.

\tpoint{Definitions} \label{def:componentdefs}
Suppose that $\mu_1$,\ldots, $\mu_k$, and $\mu$ are dominant weights.

\begin{enumerate}
\item We say that $\V_{\smu}$ is a {\em generalized PRV component} 
of $\V_{\smu_1}\otimes\cdots\otimes\V_{\smu_k}$ if there
exist $w_1$,\ldots, $w_k$, and $w$ in $\W$ such that $w^{-1}\mu=\sum_{i=1}^{k} w_i^{-1}\mu_i$.

\medskip
\item We say that $\V_{\smu}$ is a component of {\em stable multiplicity one} of
$\V_{\smu_1}\otimes\cdots\otimes\V_{\smu_k}$ if we have 
$\dim (\V_{\!m\mu_1}\otimes\cdots\otimes\V_{\!m\mu_k}\otimes\V_{\!m\mu}^{*})^{\G}=1$ for all $m\gg 0$.

\medskip
\item We say that $\V_{\smu}$ is a {\em cohomological component} of $\V_{\smu_1}\otimes\cdots\otimes\V_{\smu_k}$ if there
exist $w_1$,\ldots, $w_k$, and $w$ in $\W$ such that $w^{-1}\mu=\sum_{i=1}^{k} w_i^{-1}\mu_i$ and such that
$\invset_{w}=\sqcup_{i=1}^{k} \invset_{w_i}$.
\end{enumerate}

\np
Under the hypothesis that $\invset_{w}=\sqcup_{i=1}^{k} \invset_{w_i}$ the condition 
$w^{-1}\mu=\sum_{i=1}^{k} w_i^{-1}\mu_i$ is equivalent to the condition 
$w^{-1}\cdot\mu=\sum_{i=1}^{k} w_i^{-1}\cdot\mu_i$.  Therefore 
by Theorem I condition ({\em c}) is equivalent to having a surjective cup product map

$$
\H^{\ell(w_1)}(\X,\L_{w_1^{-1}\cdot\mu_1})\otimes\cdots\otimes
\H^{\ell(w_k)}(\X,\L_{w_k^{-1}\cdot\mu_k})\stackrel{\cup}{\longrightarrow}
\H^{\ell(w)}(\X,\L_{w^{-1}\cdot\mu})
$$

\np
which, after dualizing gives an injective map 

$$ \V_{\smu} \longrightarrow \V_{\smu_1}\otimes\cdots\otimes\V_{\smu_k}. $$

\np
In other words, we obtain a  
construction of $\V_{\smu}$ as a 
component of $\V_{\smu_1}\otimes\cdots\otimes\V_{\smu_k}$ realized through the cohomology of $\X$.

\np
Note that the conditions in \ref{def:componentdefs} 
are {\em homogeneous}: If $\V_{\smu}$ is a generalized PRV component, component of stable multiplicity one, or 
cohomological component of $\V_{\smu_1}\otimes\cdots\otimes\V_{\smu_k}$
then the same is true of $\V_{m\mu}$ as a component
of $\V_{m\mu_1}\otimes\cdots\otimes\V_{m\mu_k}$ for all $m\geq 1$.  This follows immediately from the definitions.

\np
\bpoint{Proof of Theorem II({\em a}) and restatement of II({\em b})}

\np
{\em Proof of Theorem II({\em a}).}
Every cohomological component has multiplicity one by Theorem \ref{thm:hammer}({\em b}) (the condition on 
nonzero intersection holds by the nonsymmetric version of Lemma \ref{lem:nonzero}).  
By homogeneity we conclude that homological components are of stable multiplicity one.
From Definition \ref{def:componentdefs}({\em a},{\em c}) it is clear that every cohomological component is 
a generalized PRV component.  Thus every cohomological component is a generalized PRV component of stable multiplicity
one.  \epf

\np
For the proof of part ({\em b}) it will be more convenient to work with the symmetric form of the problem.  
Applying the symmetrization procedure from \S\ref{sec:nonsymmetric-symmetric} (and replacing $k+1$ by $k$) 
we obtain the following reformulation of Theorem II({\em b}).

\tpoint{Proposition} \label{prop:IIb}
Let $\mu_1$,\ldots, $\mu_k$ be dominant weights such that 
$\dim (\V_{\!m\mu_1}\otimes\cdots\otimes\V_{\!m\mu_k})^{\G}=1$ for $m\gg 0$ 
and suppose that we have elements $w_1$,\ldots, $w_k$ such that $\sum w_i^{-1}\mu_i=0$. 
Then in either of the following two cases:

\medskip
\begin{itemize}
\item[({\em i})] at least one of $\mu_1$,\ldots, $\mu_k$ is strictly dominant, 
\item[({\em ii})] $\G$ is a classical simple group or product of classical simple groups,
\end{itemize}

\np
there exist $\overline{w}_1,\ldots, \overline{w}_k\in\W$ such that 

\begin{equation}\label{eqn:cohomologicalconditions}
 \sum_{i=1}^{k} \overline{w}_i^{-1}\mu_{i} = 0\,\,\mbox{and}\,\, 
\delpos = \bigsqcup_{i=1}^{k} \invset_{\overline{w}_i}.
\end{equation}

\np
The proof of Proposition \ref{prop:IIb} will be given in 
\S\ref{sec:proof-cohomological}  after some preliminary reduction steps.

\np
For the rest of this section 
we assume that we have fixed dominant weights 
$\mu_1$,\ldots, $\mu_k$ and Weyl group elements 
$w_1$,\ldots, $w_k$ satisfying the conditions of Proposition \ref{prop:IIb}.

\bpoint{Outline of the proof of Proposition \ref{prop:IIb}} 

\np
For $i=1$,\ldots, $k$
let $\P_i$  be the parabolic subgroup of $\G$ such that $\L_{\mu_i}$ is the pullback to
$\X$ of an ample line bundle $\L_{\tilde{\mu}_i}$ on $\G/\P_i$.  
Set $\M= \G/\P_1\times\cdots\times\G/\P_k$ and 
$\L=\L_{\tilde{\mu}_1}\bt\cdots\bt\L_{\tilde{\mu}_k}$.
The condition that $\dim (\V_{\!m\mu_1}\otimes\cdots\otimes\V_{\!m\mu_k})^{\G}=1$ for all $m\gg 1$ implies that
the GIT quotient $\M\quot\G$ with respect to $\L$ is a point. 

\np
If $w_1$,\ldots, $w_k$  are elements such that $\sum_{i=1}^{k} w_i^{-1}\mu_i=0$ 
then by 
Lemma \ref{lem:PRV-point},
the point $q=(w_1^{-1},\ldots,w_k^{-1})$ is a semi-stable point of $\M$ with a closed orbit.
Let $\H\subseteq\G$ be the stabilizer subgroup of $q$, and $\Nsh_{q}$ be the normal space to the orbit at $q$.
By the Luna slice theorem and the fact that the GIT quotient 
$\M\quot\G$ is a point we conclude that $\Sym^{\cdot}(\Nsh_{q})^{\H}$ is one-dimensional.

\np
The explicit combinatorial formula for the weights appearing in $\Nsh_{q}$ shows that a necessary condition 
for a solution of \eqref{eqn:cohomologicalconditions} to exist is that there is $v\in\W$ such that the 
weights of $v\Nsh_{q}$ are contained in $\Delta^{-}$.
In Proposition \ref{prop:reductionI} below we formulate a condition which together with the necessary
condition above guarantees the existence of a solution of \eqref{eqn:cohomologicalconditions}.
Together these two conditions are equivalent to the existence of a parabolic subalgebra $\gp$ 
with reductive part $\Lie(\H)$ such that the weights of $\Nsh_{q}$ are contained in $\gp$.

\np
Finally,
we use the restriction that $\Sym^{\cdot}(\Nsh_{q})^{\H}$ is one-dimensional to show the existence of such a 
parabolic subalgebra when $\G$ is a classical group, or for any semisimple group $\G$ under a genericity condition.

\np
\bpoint{Stabilizer subgroup of a semi-stable $\T$-fixed point}

\np
Let $\P_i$ be the parabolic with roots 
$\Delta_{\P_i}=\left\{{ \alpha \in\Delta \st \kappa(\alpha,\mu_i)\geq 0\rule{0cm}{0.38cm}}\right\}$,  
and let $\M_i=\G/\P_i$.
The stabilizer subgroup of the
point $w_i^{-1}$ in $\M_i$ is $w_i^{-1}\P_iw_i$, whose roots are 

\begin{equation}\label{eqn:parabolic-roots}
\Delta_{w_i^{-1}\P_i w_i}=
\left\{{ \alpha \in\Delta \st \kappa(w_i\alpha,\mu_i)\geq 0\rule{0cm}{0.40cm}}\right\}
=
\left\{{ \alpha \in\Delta \st \kappa(\alpha,w_i^{-1}\mu_i)\geq 0\rule{0cm}{0.40cm}}\right\}.
\end{equation}

\np
Let $\M=\M_1\times\cdots\times\M_k$ and  let 
$q$ be the point $q=(w_1^{-1},\ldots,w_{k}^{-1})$ of $\M$. We set $\H=\bigcap_{i=1}^{k} w_i^{-1}\P w_i$ to be the
stabilizer subgroup of $q$.  The condition $\sum_{i=1}^{k} w_i^{-1}\mu_i=0$ in combination with
\eqref{eqn:parabolic-roots} shows that the roots of $\H$ are given by 

\begin{equation}\label{eqn:stab-roots-II}
\Delta_{\H}=
\left\{{ \alpha \in\Delta \st \kappa(\alpha,w_i^{-1}\mu_i)= 0\,\mbox{for $i=1$,\ldots, $k$}\rule{0cm}{0.40cm}}\right\}.
\end{equation}

\np
We conclude from \eqref{eqn:stab-roots-II} that $\H$ is a reductive subgroup of $\G$.
Noting that $\T\subseteq\H$, the following lemma is another immediate consequence of \eqref{eqn:stab-roots-II}.

\tpoint{Lemma} \label{lem:HisTcondition}
We have $\H=\T$ if and only if 
the span of $\{w_i^{-1}\mu_i\}_{i=1}^{k}$ 
intersects the interior of some Weyl chamber. 
This happens, for instance, if any one of the weights $\mu_i$ is strictly dominant.

\bpoint{Torus action at fixed points of $\M$ and combinatorial deductions}\label{sec:thmIIIsetup}

\np
Let $\W_i=\left\{w\in\W \st w\mu_i=\mu_i\rule{0cm}{0.0cm}\right\} \subseteq \W$
be the stabilizer subgroup of $\mu_i$; this is the Weyl group of $\P_i$.
We will need the formula for the formal character of the tangent space of $\M_i$ at a torus fixed point.
Because of the way that the inverses of group elements enter into our formulas we make the following convention:
For any element $w$ of $\W$ and any $i$ we let $w_{s(i)}$ and $w_{l(i)}$ be respectively the shortest and longest
elements in the coset $\W_i w$.
Recall also that  for $\invset\subseteq\Delta$, $\langle\invset\rangle$ denotes the formal character 
$\sum_{\alpha\in\invset} e^{\alpha}$.

\np
With this convention, if $w_i$ is any element of $\W$, the formal character of the tangent space of $\M_i$ 
at the torus fixed point corresponding to the coset $w_i^{-1}\W_i$ is 

$$\Ch(\T_{w_i^{-1}}\M_i) = \langle  \invset_{w_{i,s(i)}} \rangle  +  \langle -\invset_{w_{i,l(i)}}^{\c} \rangle = 
\left\langle \{\alpha\in\Delta \st \kappa(\alpha,w_i^{-1}\mu_i)<0 \rule{0cm}{0.4cm}\} \right\rangle .$$

\label{sec:setup-theoremIII}

\np
The formal character of the tangent space of $\M$ at $q$ is therefore

\begin{equation}\label{eqn:TMmultiplicities}
\,\Ch(\T_{q}\M)  =  
\sum_{i=1}^{k} 
\left( \langle {\invset_{w_{i,s(i)}} \rangle + \langle -\invset_{w_i,l(i)}^{\c}} \rangle \right)
=  
\sum_{i=1}^{k} 
\left\langle \{\alpha\in\Delta \st \kappa(\alpha,w_i^{-1}\mu_i)<0\rule{0cm}{0.4cm}  \} \right\rangle .
\end{equation}

\np
Note that the multiplicity of each root $\alpha$ in the equations above is the number of $i$ for which 
$\kappa(\alpha,w_i^{-1}\mu_i)<0$.

\np
If $\alpha\not\in\Delta_{\H}$ then there is some $i$ for which $\kappa(\alpha,w_i^{-1}\mu_i)\neq 0$ and hence, by 
the condition $\sum_{i=1}^{k}w_i^{-1}\mu_i=0$, there is some $i$ for which 
$\kappa(\alpha,w_i^{-1}\mu_i)< 0$, i.e., $\alpha$ must appear as a weight in $\T_{q}\M$.  By looking at the positive
roots of $\T_{q}\M$ we therefore conclude that 

\begin{equation}\label{eqn:unionequation}
(\delpos\setminus\delpos_{\H}) = \bigcup \invset_{w_{i,s(i)}}.
\end{equation}

\np
Let $\O_{q}$ be the $\G$-orbit of $q$ in $\M$.  Since $\H$ is the stabilizer of $q$, the formal character
of the tangent space $\T_{q}\O_{q}$ is

\begin{equation}\label{eqn:TOweights}
\Ch(\T_{q}\O_{q}) = \langle\delpos\setminus \delpos_{\H}\rangle + \langle\delneg\setminus\delneg_{\H}\rangle.
\end{equation}

\np
If $\Nsh_{q}=\T_{q}\M/\T_{q}\O_{q}$ is the normal space to the orbit at $q$, then 
the union in \eqref{eqn:unionequation} is disjoint if and only if the formal character of $\Nsh_{q}$ contains
no positive root.

\np
Let $\Msh$ be the subspace of $\ggg$ spanned by the root spaces corresponding to the roots appearing in $\Nsh_{q}$.  
Comparing the multiplicities in \eqref{eqn:TMmultiplicities} and \eqref{eqn:TOweights} we conclude that
the roots of $\Msh$ are 

\begin{equation}\label{eqn:Mshweights}
\Delta_{\Msh}=\left\{ \alpha\in\Delta \st \kappa(\alpha,w_i^{-1}\mu_i)<0\,\mbox{for at least two $i\in\{1,\ldots,k\}$}
\rule{0cm}{0.4cm}\right\}.
\end{equation}
Let $\gs=\Lie(\H)$; equations \eqref{eqn:Mshweights} and \eqref{eqn:stab-roots-II} show that $\Msh$
is an $\gs$-submodule of $\ggg$.

\np
The point $q$ is not the only torus fixed point in its orbit;  for any $v\in \W$ we can act on the left
to get the torus fixed point $vq=(v w_1^{-1},\ldots, v w_{k}^{-1})$.  The weights of the normal
space $\Nsh_{vq}$ to the $\G$-orbit at $vq$ are the result of acting on the weights of $\Nsh_{q}$ by
$v$ and are hence the roots appearing in $\Ch v\Msh$.

\np
Repeating the previous arguments with the new point $vq$ and the new 
stabilizer group $v\H v^{-1}=\Stab(vq)$, gives the following result.

\tpoint{Lemma} \label{lem:disjointcondition}
For any $v\in\W$ we have 
$$
(\delpos\setminus \delpos_{v\H v^{-1}}) = \bigsqcup_{i=1}^{k} \invset_{(w_iv^{-1})_{s(i)}}
$$
if and only if $v\Msh\subseteq \gb^{-}$.  

\bpoint{Reduction to the existence of $\gp_{\Msh}$}

\tpoint{Proposition}  \label{prop:reductionI}
Suppose that there exists $v\in\W$ satisfying the conditions

\begin{enumerate}
\item[({\em i})] $v\Msh \subseteq\gb^{-}$, 

\medskip
\item[({\em ii})]  there is an element $w\in\W$ such that $\invset_{w}=\delpos_{v\H v^{-1}}$.
\end{enumerate}

\np
Then there exist  $\overline{w}_1$,\ldots, $\overline{w}_k\in \W$ such that 
$ \sum_{i=1}^{k} \overline{w}_i^{-1}\mu_{i} = 0$  and 
$\delpos = \bigsqcup_{i=1}^{k} \invset_{\overline{w}_i}.$

\np
\bpf
By condition ({\em i}) and Lemma \ref{lem:disjointcondition} we have 
$(\delpos\setminus \delpos_{v\H v^{-1}}) = \sqcup_{i=1}^{k} \invset_{(w_iv^{-1})_{s(i)}}$. 
Set $\widetilde{w}_{k+1}=w$, $\mu_{k+1}=0$, and $\widetilde{w}_i=(w_iv^{-1})_{s(i)}$ for $i=1$,\ldots, $k$.  
Conditions ({\em i}) and ({\em ii}) above and the original assumption about $w_1$,\ldots, $w_k$ imply

\begin{equation}\label{eqn:cup-reduction}
 \sum_{i=1}^{k+1} \widetilde{w}_i^{-1}\mu_{i} = 0\,\,\mbox{and}\,\, 
\delpos = \bigsqcup_{i=1}^{k+1} \invset_{\widetilde{w}_i}.
\end{equation}

\np
Equation \eqref{eqn:cup-reduction} and Theorem I show that there is a surjective cup product map 

$$
\H^{\ell(\tilde{w}_1)}(\X,\L_{\tilde{w}_1^{-1}\cdot\mu_1})\otimes\cdots\otimes
\H^{\ell(\tilde{w}_{k+1})}(\X,\L_{\tilde{w}_{k+1}^{-1}\cdot\mu_{k+1}})\stackrel{\cup}{\longrightarrow}
\H^{\N}(\X,\Ksh_{\X}).
$$

\np
Since $\H^{\ell(\tilde{w}_{k+1})}(\X,\L_{\tilde{w}_{k+1}^{-1}\cdot\mu_{k+1}})$ is the trivial module, if
we factor the map above by cupping the $k$-th and $(k+1)$-st factors together first, we obtain a surjective cup
product map onto $\H^{\N}(\X,\K_{\X})$ 
only involving the modules $\V_{\mu_1}^{*}$,\ldots, $\V^{*}_{\mu_k}$. By invoking
Theorem I again we conclude that there are $\overline{w}_1$,\ldots, $\overline{w}_k$ such that

\begin{equation}\label{eqn:ThmIII-condition}
 \sum_{i=1}^{k} \overline{w}_i^{-1}\mu_{i} = 0\,\,\mbox{and}\,\, 
\delpos = \bigsqcup_{i=1}^{k} \invset_{\overline{w}_i},
\end{equation}

\np
proving Proposition \ref{prop:reductionI}. \epf

\np
{\bf Remark.} If there do exist $\overline{w}_1$,\ldots, $\overline{w}_k$ satisfying the conclusion of
Proposition \ref{prop:IIb} it is not hard to show that there must exist $v\in\W$ so that ({\em i}) of
\ref{prop:reductionI} holds. 
As a consequence of our method of proof
we see aposteriori that there must be a $v$ so that both ({\em i}) and ({\em ii}) hold when $\G$ is a classical group
or under a genericity condition.
We do not know if condition ({\em ii})  is necessary in general.

\np
It is useful to rephrase the conditions of Proposition \ref{prop:reductionI} 
in terms of the existence of a particular parabolic subalgebra  $\gp_{\Msh}$.

\tpoint{Lemma} \label{lem:parabolic}
Let $\gs=\Lie(\H)$. 
Suppose that there exists a parabolic subalgebra $\gp_{\Msh}$ with reductive part $\gs$ such
that $\Msh\subseteq\gp_{\Msh}$.  Then conditions ({\em i}) and ({\em ii}) of Proposition \ref{prop:reductionI} hold.

\np
\bpf
Let $\gp_{\Msh}$ be such a parabolic subalgebra.  Acting by an element $v\in\W$ we can conjugate $\gp_{\Msh}$ so that 
$\gb^{-}\subseteq v\gp_{\Msh}$.
This implies that $v\Msh\subseteq \gb^{-}$.  
Since $v\gs$ is the radical of a parabolic subalgebra containing $\gb^{-}$, 
if $w$ is the longest element of the Weyl group of $v\gs$ then
$\invset_{w}=\Delta^{+}_{v\gs}=\Delta^{+}_{v\H v^{-1}}$. \epf

\np
{\bf Remark.}
If there exists $v\in\W$ such that condition ({\em ii}) of Proposition \ref{prop:reductionI} holds then 
one can show that $\gp:=\gb^{-}+v\gs$ is a parabolic subalgebra of $\ggg$. 
If condition ({\em i}) also holds for this $v$
then $\gp_{\Msh}:=v^{-1}\gp$ is a parabolic subalgebra satisfying the conditions of Lemma \ref{lem:parabolic}.  
Therefore the existence of the parabolic $\gp_{\Msh}$ is equivalent to the conditions in 
Proposition \ref{prop:reductionI}.  Since we will not need this direction of the equivalence we omit the justification
of the first assertion.

\bpoint{GIT consequences of the stable multiplicity one condition}

\np
Let $\L$ be the line bundle on $\M$ whose pullback to $\X^{k}$ is $\L_{\mu_1}\bt\cdots\bt\L_{\mu_k}$. 
Then $\L$ is a $\G$-equivariant ample line bundle on $\M$.  By the stable multiplicity one condition we have
$\dim(\M,\L^{m})^{\G}=1$ for all $m\gg 1$, and so the GIT quotient $\M\quot\G$ is a point.  

\np
The weight of  $\L$ at $q$ is $\sum_{i=1}^{k}w_i^{-1}\mu_i=0$.  By Lemma \ref{lem:PRV-point} this means
that $q$ is a semi-stable point with a closed orbit.  By the Luna slice theorem,
\cite[Th\'eor\`em du Slice \'Etale,pg. 97]{lu1},
$\spec(\Sym(\Nsh_{q}^{*})^{\H})$ and the image of $q$ 
in the GIT quotient $\M\quot\G$ have a common \etale neighbourhood.  Hence 
$\dim(\Nsh_{q}/\H) = \dim(\M\quot\G)=0$, i.e., 
$\dim \Sym^{\cdot}(\Nsh_{q}^{*})^{\H}=1$. 
Passing to the level of Lie algebras and dualizing we obtain $\dim \Sym^{\cdot}(\Nsh_{q})^{\gs}=1$.

\np
Since $\Msh$ is isomorphic to an $\gs$-submodule of  $\Nsh_{q}$ 
we arrive at the following consequence of the stable multiplicity one condition:

\np
\tpoint{Lemma}  \label{lem:GITconsequences}
Under the hypotheses of Proposition \ref{prop:IIb} and with the notation of \S\ref{sec:thmIIIsetup},
we have $\dim \Sym^{\cdot}(\Msh)^{\gs}=1$, i.e., 
$\Sym^{\cdot}(\Msh)^{\gs}$ consists of just the constants.

\bpoint{Proof of Proposition \ref{prop:IIb}}
\label{sec:proof-cohomological}

\np
By Proposition \ref{prop:reductionI} and Lemma \ref{lem:parabolic}, to prove Proposition \ref{prop:IIb} it is enough
to show the existence of the parabolic subalgebra  $\gp_{\Msh}$.  By Lemma \ref{lem:GITconsequences} we may assume that $\dim \Sym^{\cdot}(\Msh)^{\gs}=1$.

\np
{\em Proof of \ref{prop:IIb}(i)---} If any one of the weights $\mu_1$,\ldots, $\mu_k$ 
is strictly dominant, or more generally,
if the span of $\{w_i^{-1}\mu_i\}_{i=1}^{k}$ 
intersects the interior of some Weyl chamber then by 
Lemma \ref{lem:HisTcondition} $\H=\T$ and so $\gs=\Lie(\T)=\tt$ and $\delpos_{\tt}=\emptyset$. 
The condition that $\dim(\Sym^{\cdot}(\Msh)^{\tt})=1$ is then equivalent 
to the condition that no non-trivial non-negative combination of weights of $\Msh$ 
is zero. 
Hence by Farkas's lemma the weights of $\Msh$ all lie strictly on
one side of a hyperplane and the cone dual to the cone they span is open.  We may therefore pick a weight in the 
interior of the dual cone which is not on any hyperplane of the Weyl chambers.  The roots lying on the positive side
of this hyperplane give the parabolic subalgebra $\gp_{\Msh}$. \epf

\np
{\em Proof of \ref{prop:IIb}(ii)---} 
Equation \eqref{eqn:stab-roots-II} shows 
that the roots of $\gs$ are given by the vanishing of linear forms and hence
$\gs$ is the reductive part of a parabolic subalgebra.
Let $\gt$ be the center of $\gs$.  For any 
$\nu \in \gt^* \setminus\{0\}$  set

$$\ggg^\nu = \left\{{x \in \ggg \, \st \, [t, x] = \nu(t) x \,\,\mbox{for all}\, t \in \gt\rule{0cm}{0.39cm}}\right\}.$$

\np
Following Kostant, \cite{K} we call $\nu \in \gt^* \backslash \{0\}$ an $\gt$--root if $\ggg^{\nu}\neq 0$.  
Let $\Rsh$ be the set of $\gt$--roots of $\ggg$ and $\Ssh$ the
subset of those $\gt$--roots appearing in $\Msh$, 
so that $\Msh=\opp_{\nu\in\Ssh} \ggg^{\nu}$.

\np
A subset $\Rsh'$ of $\Rsh$ is called {\em saturated} if 
whenever $\nu \in \Rsh'$ and $r \nu \in \Rsh$ for some $r \in \QQ_+$ then $r \nu \in \Rsh'$ as well.  
It follows from \eqref{eqn:Mshweights} that $\Ssh$ is a saturated subset of $\Rsh$.

\np
In the appendix we establish the following result (Theorem \ref{thm:classical}).

\np
{\em Theorem ---}
Let $\ggg$ be a classical Lie algebra, $\gs$ be a subalgebra which is the reductive part of a parabolic subalgebra
of $\ggg$, $\Ssh$ be a saturated subset of the $\gt$--roots $\Rsh$, and $\Msh = \oplus_{\nu \in \Ssh} \ggg^\nu$. 
If $\dim(\Sym^\cdot(\Msh))^{\gs} = 1$, then there exists a parabolic subalgebra $\gp_{\Msh}\subseteq\ggg$ 
with reductive part $\gs$ such that $\Msh \subseteq \gp_{\Msh}$.

\np
Thus when $\G$ is a simple classical group or a product of simple classical groups, 
the above theorem along with the previous reductions establish Proposition \ref{prop:IIb}. \epf

\bpoint{Comments on the existence of $\gp_{\Msh}$} We want to prove the implication 

\begin{equation}\label{eqn:condition}
\mbox{generalized PRV of stable multiplicity one $\implies$ cohomological}
\end{equation}

\np
for components of a tensor product.  Our 
strategy has been to reduce the proof of this implication to the existence of the parabolic $\gp_{\Msh}$.  
As is clear from the proof of Proposition \ref{prop:IIb} we can establish
\eqref{eqn:condition} in a few more cases than stated in Theorem II({\em b}).  Namely we can replace ({\em i}) by the
condition that the span of $\{w_i^{-1}\mu_i\}$ intersects the interior of some Weyl chamber, and in ({\em ii})
we can allow $\ggg$ to be a semisimple algebra whose simple ideals are classical or of type $\G_2$
since Theorem \ref{thm:classical} extends to this case.

\np
We also have a proof of the existence of $\gp_{\Msh}$ in the simply-laced case under a hypothesis 
somewhat complementary to the condition in ({\em i}): the semisimple part of each 
$w_i^{-1}\P_iw_i$ be exactly $\H$;  we omit this argument here.

\np
Nonetheless, we have not been able to establish the implication \eqref{eqn:condition}
in full generality.  There are two points where our strategy can be improved.  First, it is not clear to us
that condition ({\em ii}) in Proposition \ref{prop:reductionI} is necessary, although we believe that it is, i.e.,
that proving the implication \eqref{eqn:condition} is equivalent to showing the existence of $\gp_{\Msh}$.
Second, it is not clear to us how to fully utilize the condition that $\Msh$ comes from the data 
in the hypotheses of Proposition \ref{prop:IIb}.  One of the implications of this data 
is that $\dim\Sym^{\cdot}(\Msh)^{\gs}=1$.  In the classical cases (and $\G_2$) or under the genericity condition 
our arguments show that for any $\gs$-submodule $\Msh$ of $\ggg$ satisfying the condition 
$\dim\Sym^{\cdot}(\Msh)^{\gs}=1$ there exists a parabolic subalgebra $\gp_{\Msh}$ as in Lemma \ref{lem:parabolic}.
For groups of the other exceptional types we can find 
explicit examples of $\gs$-submodules $\Msh\subseteq \ggg$ such that $\dim\Sym^{\cdot}(\Msh)^{\gs}=1$ 
but which are not contained in a parabolic $\gp_{\Msh}$ with reductive part $\gs$ (see the Appendix). 
However we have not been able to turn these
examples into counterexamples to \eqref{eqn:condition}
since we are not able to construct data satisfying the hypotheses in Proposition \ref{prop:IIb} which produce
these $\Msh$'s.  Perhaps a more careful analysis of the requirement that $\Msh$ come from this particular
kind of representation-theoretic problem would be sufficient to prove \eqref{eqn:condition} in general.

\section{Littlewood-Richardson cone and bounds on multiplicities}

\np
In this section we relate cohomological components to the boundary of codimension $n$ of the Littlewood-Richardson 
cone, where $n$ is the rank of $\G$.
We extend the argument used in Lemma \ref{lem:key}({\em c}) to give a bound on the multiplicities
of components generalizing that of Theorem \ref{thm:hammer}({\em b}).   
We show how this bound imposes conditions on the order of growth of multiplicities on the boundary of the cone.

\np
We are working over an algebraically closed field of characteristic zero.  For notational convenience 
we will assume that the field is $\CC$.
It is advisable to look over \S\ref{sec:LR-cone} before reading this section.

\bpoint{Cohomological components are of codimension $n$}
\label{sec:cohomologicalcodimn}

\tpoint{Proposition} \label{prop:cohomologicalvertex}
Let $\mu_1$,\ldots, $\mu_k$, and $\mu$ be dominant weights.
If $\V_{\mu}$ is a cohomological component of $\V_{\smu_1}\otimes\cdots\otimes\V_{\smu_k}$ then
$\mu$ is a vertex of the Littlewood-Richardson polytope $\LRP(\mu_1,\ldots, \mu_k)$.

\np
\bpf
Recall that for any polytope $\mathcal{P}$ 
in a vector space $\E$ 
every point $p\in\mathcal{P}$ lies on the relative interior of a unique face.  The dimension of this face
is the same as the dimension of the subspace of $\E$ spanned by the set
$\left\{\epsilon\in\E \st p\pm \epsilon\in\mathcal{P}\right\}.$

\np
Let $w_1$,\ldots, $w_k$, and $w$ be the elements from \S\ref{def:componentdefs}({\em c}) expressing $\V_{\smu}$
as a cohomological component of $\V_{\smu_1}\otimes\cdots\otimes\V_{\smu_k}$.
If $\epsilon$ is such that $(\mu_1,\ldots, \mu_k,\mu\pm\epsilon)$ is in 
$\LRP(\mu_1,\ldots,\mu_k)$ 
we must have that $\sum_{i=1}^{k} w_i^{-1}\mu_i-w^{-1}(\mu\pm\epsilon)=\mp w^{-1}\epsilon$ 
belongs to $\Span_{\QQ\geq0} \delpos$.
This is only possible if $w^{-1}\epsilon=0$, i.e., $\epsilon=0$ and so $\mu$ is a vertex of 
$\LRP(\mu_1,\ldots, \mu_k)$.
\epf

\np
\tpoint{Lemma} 
\label{lem:Iempty}
With the notation of \S\ref{sec:LR-cone},
if $\I=\emptyset$ then $w_1$,\ldots, $w_k$, and $w$ satisfy \eqref{eqn:Iconditions} with respect
to $\I$ if and only if $\invset_{w}=\sqcup_{i=1}^{k} \invset_{w_i}$ and $\cap_{i=1}^{k}[\Omega_{w_i}]\cdot[\X_{w}]=1$.

\bpf
If $\I=\emptyset$ then
\eqref{eqn:Iconditions}({\em iii}) becomes $\sum_{i=1}^{k}w_i^{-1}\cdot 0 - w^{-1}\cdot0=0$. 
Hence by \eqref{eqn:Iconditions}({\em i})
and Theorem \ref{thm:hammer}({\em d}) we conclude that $\invset_{w}=\sqcup_{i=1}^{k} \invset_{w_i}$.
Conversely, if $w_1$,\ldots, $w_k$, and $w$ are such that 
$\invset_{w}=\sqcup_{i=1}^{k} \invset_{w_i}$ and $\cap_{i=1}^{k}[\Omega_{w_i}]\cdot[\X_{w}]=1$ then
\eqref{eqn:Iconditions}({\em i},{\em iii}) hold
and \eqref{eqn:Iconditions}({\em ii}) is automatic since $\P_{\I}=\B$ and so $\W_{\P_{\I}}=\{e\}$.   
\epf

\np
Theorem IV now follows immediately by combining Theorem \ref{thm:Ressayre}, Lemma \ref{lem:Iempty},  and Theorem I.

\np
\tpoint{Lemma} \label{lem:unproved}
Suppose $w_1$,\ldots, $w_k$, and $w$ are such that $\invset_{w}=\sqcup_{i=1}^{k}\invset_{w_i}$ and 
$\cap_{i=1}^{k}[\Omega_{w_i}]\cdot[\X_{w}]=1$. Then for any set $\I$ of simple roots, if we set 
$\overline{w}_i$ to be the shortest element in $w_i\W_{\P_{\I}}$ for $i=1$,\ldots, $k$, and $\overline{w}$ 
to be the shortest element in $w\W_{\P_{\I}}$ then 
$\overline{w}_1$,\ldots, $\overline{w}_k$ and $\overline{w}$ satisfy \eqref{eqn:Iconditions} with respect to $\I$.

\np
\bpf
The only part of \eqref{eqn:Iconditions} which is not immediate is that 
$\cap_{i=1}^{k} [\Omega_{\overline{w}_i}]\cdot[\X_{\overline{w}}]=1$. 
Since we will only use this lemma to derive a picture of the Littlewood-Richardson cone around a 
cohomological component in Corollary \ref{cor:nearvertex}({\em b}), we only sketch the argument.

\np
Let $\G'$ be the Levi component of $\P_{\I}$, $\B'=\G'\cap \B$ and $\X'=\G'/\B'$.  Write
$w_i=\overline{w}_i u_i$ for $i=1$,\ldots, $k$, and  
$w=\overline{w}u$ with each $u_i$, respectively $u$ in $\W_{\P_{\I}}$.
Let $\pi\colon\X\longrightarrow\G/\P_{\I}$ be the projection.
The general fibre
of $\pi|_{\Omega_{w_i}}\colon \Omega_{w_i}\longrightarrow \pi(\Omega_{w_i})$ has cohomology class
$[\Omega'_{u_i}]$ in $\H^{*}(\X',\ZZ)$ and similarly the general fibre of 
of $\pi|_{\X_{w}}\colon \X_{w}\longrightarrow \pi(\X_{w})$ has cohomology class $[\X'_{u}]$.
Furthermore, the intersection numbers 
$\cap_{i=1}^{k} [\Omega_{\overline{w}_i}]\cdot[\X_{\overline{w}}]$ and
$\cap_{i=1}^{k} [\pi(\Omega_{w_i})]\cdot[\pi(\X_{w})]$ are equal,
the last intersection being computed in $\H^{*}(\G/\P_{\I},\ZZ)$.
Using Kleiman's transversality theorem one can factor the intersection number on $\X$ as a product
of an intersection number on $\G/\P_{\I}$ and an intersection number on the fibre:

\begin{equation}\label{eqn:factorintersection}
\bigcap_{i=1}^{k}[\Omega_{w_i}]\cdot[\X_{w}] = 
\left({\bigcap_{i=1}^{k} [\Omega_{\overline{w}_i}]\cdot [\X_{\overline{w}}]}\right)
\left({\bigcap_{i=1}^{k} [\Omega'_{u_i}]\cdot [\X'_{u}]}\right).
\end{equation}

\np
By hypothesis the left side of 
\eqref{eqn:factorintersection} is equal to one.  Therefore both of the intersection numbers
on the right side are equal to one, and in particular
$\cap_{i=1}^{k} [\Omega_{\overline{w}_i}]\cdot[\X_{\overline{w}}]=1$. \epf

\np
Note that the evaluation of the condition in Theorem \ref{thm:Ressayre}({\em a}) does not 
depend on the representatives chosen in the cosets $w_i\W_{\P_{\I}}$.
Thus by Lemma \ref{lem:unproved} if one considers the conditions imposed 
on $(\mu_1,\ldots, \mu_k,\mu)\in\LRC(k)$ obtained by writing
$\sum_{i=1}^{k} w_i^{-1}\mu_i-w^{-1}\mu$ in terms of the simple roots and requiring that any subset
of size $r$ of the coordinates vanish, one obtains a codimension $n-r$ face of $\LRC(k)$. 
Each of these faces contain the codimension $n$ face where all coefficients are zero, consisting of cohomological
components.

\bpoint{Bounds on multiplicities}
Recall that $\mult(\V_{\smu},\V)$ denotes the multiplicity of $\V_{\smu}$ in $\V$.

\tpoint{Definition} For any multiset $\Ssh$ with support in $\delpos$ (i.e., $\Ssh$ is a set of positive roots
with multiplicities) we define the partition function $\kostP_{\!\Ssh}$ by setting 
$\kostP_{\!\Ssh}(\gamma)$ to be the coefficient of $e^{\gamma}$ in the series 
$\prod_{\alpha\in\delpos} \frac{1}{(1-e^{\alpha})^{m_{\Ssh}(\alpha)}}$, where $m_{\Ssh}(\alpha)$ is the multiplicity
of $\alpha$ in $\Ssh$.  

\np
If $\Ssh=\delpos$ then $\kostP_{\!\Ssh}$ is the usual Kostant partition function.
For any $\Ssh$, $\kostP_{\!\Ssh}(0)=1$, and $\kostP_{\!\Ssh}(\gamma)=0$ if $\gamma\not\in\Span_{\ZZ\geq0}\Ssh$ 
(in particular, $\kostP_{\!\Ssh}(\gamma)=0$ if $\gamma \not \in \Span_{\ZZ \geq 0} \delpos$). 

\np
By convention the multiplicities in a multiset are always non-negative.   If $\Ssh'$ is a multiset
with support in $\delpos$, then the set difference $\Ssh'\setminus\delpos$ is the multiset obtained by reducing
each of the nonzero multiplicities in $\Ssh'$ by one.  The operation of union with multiplicity consists
of adding the multiplicities and is denoted by $\Umult$.

\tpoint{Theorem} \label{thm:multbound}
Let $w_1$,\ldots, $w_k$, and $w$ be such that $\cap_{i=1}^{k}[\Omega_{w_i}]\cdot[\X_{w}]\neq 0$
and set

$$\Ssh= \left(\left({\Umult\!\!{\vphantom{\bigcup}}_{i=1}^{k} 
\invset_{w_i}^{\c}}\right)\Umult \invset_{w}\right)\setminus\delpos.$$

\np
Then for any dominant weights $\mu_1$,\ldots, $\mu_k$, and $\mu$ we have the bound

$$\mult(\V_{\mu},\V_{\mu_1}\otimes\cdots\otimes\V_{\mu_k}) \leq 
\kostP_{\!\Ssh}\!\left(\sum_{i=1}^{k}w_i^{-1}\mu_i - w^{-1}\mu\right).$$

\np
\bpf
As in the proof of Theorem \ref{thm:hammer} 
let $v_i=w_i^{-1}\wo$ for $i=1,\ldots, k$, $v_{k+1}=w^{-1}$, 
let $\vv_i$ be a reduced word with product $v_i$ for $i=1$,\ldots $k+1$, 
and set $\vseq=(\vv_1,\ldots, \vv_{k+1})$.
By Corollary \ref{cor:map-degree} the degree of 
$f_{\vseq}\colon\Y_{\svseq}\longrightarrow \X^{k+1}$ is given by the intersection number 
$\cap_{i=1}^{k+1} [\Omega_{\wo v_i^{-1}}]= \cap_{i=1}^{k}[\Omega_{w_i}]\cdot[\X_{w}]$ which is nonzero
by assumption, i.e., $f_{\vseq}$ is surjective.

\np
Let $\lambda_i=-\wo\mu_i$ for $i=1\ldots, k$ and $\lambda_{k+1}=\mu$.  
Set $\L$ to be the line bundle $\L_{\lambda_1}\bt\cdots\bt\L_{\lambda_{k+1}}$ on $\X^{k+1}$,
so that the multiplicity of $\V_{\mu}$ in $\V_{\mu_1}\otimes\cdots\otimes\V_{\mu_k}$ is
equal to $\dim \H^{0}(\X^{k+1},\L)^{\G}$.
Since $f_{\vseq}$ is surjective, pullback induces an inclusion 
$\H^{0}(\Y_{\svseq},f_{\vseq}^{*}\L)^{\G} \stackrel{f_{\vseq}^{*}}{\hookleftarrow} \H^{0}(\X^{k+1},\L)^{\G}$.

\np
As in the proof of Lemma \ref{lem:key}, let 
$\Z:=\Z_{\vv_1}\times\cdots\times\Z_{\vv_k}$ be the fibre of $\fo$ over $e\in\X$.  The pushforward 
${\fo}_{*}(f_{\vseq}^{*}\L)$ is a vector bundle on $\X$ whose fibre over $e$ is $\H^{0}(\Z,f_{\vseq}^{*}{\L}|_{\Z})$.
By the equivalence between $\G$-invariant sections of a $\G$-bundle on $\X$ and $\B$-invariant
vectors of the fibre over $e$, restriction to $\Z$ induces an isomorphism 

$$\H^{0}(\Z,f_{\vseq}^{*}{\L}|_{\Z})^{\B}\stackrel{\sim}{\longleftarrow}\H^{0}(\Y_{\svseq},f_{\vseq}^{*}\L)^{\G}.$$

\np
Let
$p=p_{\vseq}\in\Z$ be the maximum point of $\Y_{\svseq}$.  Either by the theorem of 
Bia\l ynicki-Birula \cite{bial} or by following through the construction in \S\ref{sec:BSDH}
one checks that there is a $\T$-stable
open affine subset $\U\subset\Z$ containing $p$ and isomorphic to affine $k\N$-space.  In fact

\begin{equation}\label{eqn:Vdescription}
\U=(\B v_1\B)/\B \times (\B v_2\B)/\B \times \cdots\times (\B v_{k+1}\B)/\B 
\end{equation}

\np
and so $\U$ is $\B$-stable.  Since $\U$ is open in the irreducible variety $\Z$, restriction
gives an inclusion $\H^{0}(\U,f_{\vseq}^{*}\L|_{\U})^{\B}\hookleftarrow\H^{0}(\Z,f_{\vseq}^{*}{\L}|_{\Z})^{\B}$.
Thus to bound the $\mult(\V_{\smu},\V_{\mu_1}\otimes\cdots\otimes\V_{\mu_k})$ 
it suffices to bound $\dim \H^{0}(\U,f_{\vseq}^{*}\L|_{\U})^{\B}$. 
Set $\gamma=\sum_{i=1}^{k} w_i^{-1}\mu_i-w^{-1}\mu$;  we will show that 
$\dim \H^{0}(\U,f_{\vseq}^{*}\L|_{\U})^{\B}\leq\kostP_{\!\Ssh}(\gamma)$.

\np
Let $\L_{\U}=(f_{\vseq}^{*}\L)|_{\U}$.
Since $\U$ is isomorphic to affine space, $\L_{\U}$ is (non-equivariantly) trivial on $\U$. Let 
$\so$ be a section of $\L_{\U}$ which is nowhere vanishing. The torus $\T$ takes $\so$ to another nowhere vanishing
section which must therefore be a multiple of $\so$, i.e., $\T$ acts on $\so$ via a weight.  This must be the
same as the weight of the action of  $\L_{\U}$ at $p$, and so $\T$ acts on $\so$ with weight $\gamma$.
Let $\B^{+}$ be the unipotent radical of $\B$.  By the same reasoning,
$\B^{+}$ must take $\so$ to a multiple of itself. Since $\B^{+}$ has only the trivial 
one-dimensional representation $\so$ must be fixed by $\B^{+}$.

\np
Every section $s\in\H^{0}(\U,\L_{\U})$ can be written as $s=\so h$ for some function $h\in\H^{0}(\U,\Osh_{\U})$.
The section $s$ is $\B$-invariant if and only if $h$ is $\B^{+}$-invariant  and $h$ is an eigenfunction of $\T$
on which $\T$ acts via $-\gamma$.
For any weight $\delta$, let $\H^{0}(\U,\Osh_{\U})_{\delta}$ denote the space of 
eigenfunctions of $\T$ on which $\T$ acts via $\delta$.  
Let $\gn$ be the Lie algebra of $\B^{+}$ (i.e, the nilpotent radical of $\gb$); $\gn$ acts on 
$\H^{0}(\U,\Osh_{\U})$ via derivations.
By the correspondence above between 
sections of $\L_{\U}$ and functions on $\U$ we have 
$\dim \H^0(\U,\L_{\U})^{\B}=\dim \H^{0}(\U,\Osh_{\U})^{\gn}_{-\gamma}$. 

\np
Let $\U_i=(\B v_i\B)/\B$.
The formal character of the tangent space $\T_{p_{v_i}}\U_i=\T_{p_{v_i}}\Z_{v_i}$ 
is $\langle\invset_{v_i^{-1}}\rangle$ (see \eqref{eqn:Z-tgtweights}) and hence the affine 
space $\U_i$ is equal to $\spec(\CC[z_{i,-\alpha}]_{\alpha\in\invset_{v_i^{-1}}})$ 
where each $z_{i,-\alpha}$ is an independent variable on which $\T$ acts via the weight $-\alpha$.
Since the roots in $\invset_{v_i^{-1}}$ are distinct, the variables
$z_{i,-\alpha}$ are uniquely determined up to scalar.

\np
Let $\R=\CC[z_{i,-\alpha}]_{i=1,\ldots, k+1, \alpha\in\invset_{v_i^{-1}}}$.
By \eqref{eqn:Vdescription} and the description of $\U_i$ above 
$\R=\H^{0}(\U,\Osh_{\U})$ (or, equivalently, $\U=\spec(\R)$).  Let

$$
\Ssh'=
{\displaystyle\Umult\!\!{\displaystyle\vphantom{\displaystyle\bigcup}}_{i=1}^{k+1} 
\invset_{v_i^{-1}}}
=
\left({\displaystyle\Umult\!\!{\displaystyle\vphantom{\displaystyle\bigcup}}_{i=1}^{k} 
\invset_{w_i}^{\c}}\right)\Umult \invset_{w} 
$$

\np
so that $\Ssh=\Ssh'\setminus\delpos$.  The number of monomials in the variables $\{z_{i,-\alpha}\}$ 
of weight $-\gamma$ is $\kostP_{\!\Ssh'}(\gamma)$, i.e., 
$\dim\H^{0}(\U,\Osh_{\U})_{-\gamma}=\kostP_{\!\Ssh'}(\gamma)$.  
By analyzing the condition of $\gn$-invariance
we will obtain the bound  $\dim\H^{0}(\U,\Osh_{\U})_{-\gamma}^{\gn}\leq\kostP_{\!\Ssh}(\gamma)$.

\np
For each $\beta\in\delpos$ let $\die_{\beta}$ be the vector field giving the action of a nonzero element of
$\ggg^{\beta}\subseteq\gn$ 
on $\U$.  
Each $\die_{\beta}$ is a first-order differential operator with polynomial coefficients.  
Separating the coefficients into homogeneous pieces we may write

$$
\die_{\beta} = \die_{\beta}^{0} + \die_{\beta}^{1} + \die_{\beta}^{2}+ \cdots
$$

\np
where $\die_{\beta}^{j}$ denotes the component of $\die_{\beta}$ whose coefficients are homogeneous of degree $j$. 
In particular, 
$$
\die_{\beta}^{0} = \sum_{i\colon\beta\in\invset_{v_i}^{-1}} c_{(i,\beta)} \die_{z_{i,-\beta}} 
$$
where $c_{(i,\beta)}\in\CC$ are nonzero constants.
Each $\die_{\beta}$ is an operator of $\Delta$-degree $\beta$, i.e., 
$$\die_{\beta}\left({\H^{0}(\U,\Osh_{\U})_{\delta}}\right)\subseteq\H^{0}(\U,\Osh_{\U})_{\delta+\beta}$$
for each weight $\delta$. 
Thus, setting $\die_{1}=\sum_{\beta\in\delpos}\die_{\beta}$ we obtain

\begin{eqnarray*}
\H^{0}(\U,\Osh_{\U})_{-\gamma}^{\gn} & =
\bigcap_{\beta\in\delpos} \ker\left({
\H^{0}(\U,\Osh_{\U})_{-\gamma} \stackrel{\die_{\beta}}{\longrightarrow} \H^{0}(\U,\Osh_{\U})_{-\gamma+\beta}}\right) \\
& =
\ker\left({
\H^{0}(\U,\Osh_{\U})_{-\gamma} \stackrel{\die_1}{\longrightarrow} \opp_{\beta\in\delpos} \H^{0}(\U,\Osh_{\U})_{-\gamma+\beta}
}\right).
\end{eqnarray*}

\np
In order to estimate the size of the kernel we deform the differential operator.  
For any $s\in\CC$ and $\beta\in\delpos$ let 
$$
\die_{\beta,s} = \die_{\beta}^{0} + s \die_{\beta}^{1} + s^2 \die_{\beta}^{2}+ \cdots.
$$

\np
Each $\die_{\beta,s}$ is again a homogeneous operator
of $\Delta$-degree $\beta$.  Set $\die_{s}=\sum_{\beta\in\delpos} \die_{\beta,s}$ and consider the maps

\begin{equation}\label{eqn:die-s}
\H^{0}(\U,\Osh_{\U})_{-\gamma} \stackrel{\die_s}{\longrightarrow} \opp_{\beta\in\delpos} \H^{0}(\U,\Osh_{\U})_{-\gamma+\beta}.
\end{equation}

\np
If $s\neq 0$ then by simultaneously scaling the variables $z_{i,-\alpha}$ we can transform $\die_{s}$ 
into a nonzero multiple of $\die_{1}$, i.e., for all
$s\neq 0$ the dimension of the kernel of $\die_{s}$ is the same.  Specializing to $s=0$ can only cause the kernel
to increase and therefore $\dim\H^{0}(\U,\Osh_{\U})_{-\gamma}^{\gn} = \dim \ker(\die_{1})\leq \dim\ker(\die_0)$.

\np
Each $\die_{\beta,0}=\die_{\beta}^{0}$ is a constant coefficient differential operator. Therefore
if $\R_1$ is the degree one piece of 
$\H^{0}(\U,\Osh_{\U})$
and $\R_1^{\circ}$ the subspace annihilated
by $\{\die_{\beta,0}\}_{\beta\in\delpos}$ then the subring of $\H^{0}(\U,\Osh_{\U})$ consisting 
of polynomials annihilated by $\{\die_{\beta,0}\}_{\beta\in\delpos}$ is $\Sym^{\cdot}(\R_1^{\circ})$.

\np
For any weight $\delta$ let $\R_{1,\delta}$ and $\R_{1,\delta}^{\circ}$ denote the subspaces of 
weight $\delta$ of $\R_1$ and $\R_1^{\circ}$ respectively.
Since $\die_{\beta,0}$ acts only on the variables of weight $-\beta$, if $\R_{1,-\beta}\neq 0$ then 
$\R_{1,-\beta}^{\circ}$ is of codimension one in $\R_{1,-\beta}$.
I.e., the effect of each $\die_{\beta,0}$ is to ``reduce the number of variables of weight $-\beta$ by one''.
Therefore $\dim\ker\die_{0}=\dim \Sym^{\cdot}(\R_1^{\circ})_{-\gamma}=\kostP_{\!\Ssh}(\gamma)$.  \epf

\np
\bpoint{Corollaries of Theorem \ref{thm:multbound}}
Note that if $\invset_w = \sqcup_{i=1}^k \invset_{w_i}$ then $\Ssh=(k-1)\delpos$, i.e, the
multiset in which each positive root has multiplicity $k-1$.  Applying Theorem \ref{thm:multbound} yields

\begin{equation}\label{eqn:cohbound}
\mult(\V_{\smu},\V_{\smu_1}\otimes\cdots\otimes\V_{\smu_k})\leq\kostP_{\!(k-1)\delpos}
\left({ \sum_{i=1}^{k}w_i^{-1}\mu_i - w^{-1}\mu }\right).
\end{equation}

\np
The most interesting case is when $k=2$, i.e. $\Ssh = \Delta^+$.   Taking $w_1=w_2=w=e$ 
in \eqref{eqn:cohbound}
we obtain

\tpoint{Corollary} \label{cor:Steinberg} If $\mu_1, \mu_2$, and $\mu$ are dominant weights, then 

\begin{equation}\label{eqn:Steinberg}
\mult(\V_{\smu}, \V_{\smu_1} \otimes \V_{\smu_2}) \leq \kostP (\mu_1 + \mu_2 - \mu),
\end{equation}

\np
where $\kostP = \kostP_{\!\delpos}$ is the usual Kostant partition function.

\np
{\bf Remark.}
The right hand side of \eqref{eqn:Steinberg} is the leading term of Steinberg's multiplicity formula \cite{St}. 
Since Steinberg's formula is a signed double summation over the Weyl group, this bound 
does not seem to be obvious. 

\np
{\bf Remark.}
If $\mu$ is sufficiently close to $\mu_1+\mu_2$ then the Brauer-Klimyk multiplicity formula shows
that $\mult(\V_{\smu},\V_{\smu_1}\otimes\V_{\smu_2})=\kostP(\mu_1+\mu_2-\mu)$, i.e., 
for $\mu$ sufficiently close to $\mu_1+\mu_2$, the upper bound in Corollary \ref{cor:Steinberg} is exact.

\np
For any set $\I$ of simple roots, we define $\Delta_{\I}$ to be the roots of the Levi component
of $\P_{\I}$.  Equivalently, $\Delta_{\I}$ is the subset of $\Delta$ consisting of those roots in $\Span_{\ZZ}\I$.

\tpoint{Corollary} \label{cor:boundarybound}
Suppose that $(\mu_1,\ldots,\mu_k,\mu)\in\LRC(k)$ is on a face of the type from Theorem \ref{thm:Ressayre}({\em a}).
Let $\I$ be the corresponding set of simple roots, 
and $w_1$,\ldots, $w_k$, and $w$ be the corresponding elements satisfying conditions \eqref{eqn:Iconditions} with
respect to $\I$.
Set $\gamma=\sum_{i=1}^{k} w_i^{-1}\mu_k-w^{-1}\mu$.  Then 
$$\mult(\V_{\mu},\V_{\mu_1}\otimes\cdots\otimes\V_{\mu_k}) \leq 
\kostP_{\!(k-1)\delpos_{\I}}\!(\gamma).$$

\np
\bpf 
For any set $\I$ of simple roots and any element $v\in\W$ 
the shortest element in $v\W_{\P_{\I}}$ is the unique element $v'\in v\W_{\P_{\I}}$ such that
$\invset_{v'}\cap \delpos_{\I}=\emptyset$.
Applying this observation to the elements $w_1$,\ldots, $w_k$ and $w$, we see that 
$$\Ssh= \left(\left({\Umult\!\!{\vphantom{\bigcup}}_{i=1}^{k} 
\invset_{w_i}^{\c}}\right)\Umult \invset_{w}\right)\setminus\delpos
=(k-1)\delpos_{\I}\stackrel{\mbox{\scriptsize mult}}{\bigsqcup}
\left\{
\mbox{
\begin{minipage}{0.35\textwidth}
roots outside of $\delpos_{\I}$
counted with some multiplicities
\end{minipage}}
\right\}.$$

\np By Theorem \ref{thm:multbound} we know that $\kostP_{\!\Ssh}(\gamma)$ bounds 
$\mult(\V_{\smu},\V_{\smu_1}\otimes\cdots\otimes\V_{\smu_k})$.  
To write a root outside of $\delpos_{\I}$ as a combination of simple roots it is necessary to use at least one
simple root $\alpha_j$ not contained in $\I$.
By the condition of Theorem
\ref{thm:Ressayre}({\em a}) if we write $\gamma$ as a combination of simple roots the coefficient of every simple root 
$\alpha_j$ outside of $\I$ is zero.  Therefore the roots outside of $\delpos_{\I}$ do not contribute to
$\kostP_{\!\Ssh}(\gamma)$ 
and so $\kostP_{\!\Ssh}(\gamma)=\kostP_{\!(k-1)\delpos_{\I}}(\gamma)$.
\epf

\tpoint{Corollary}\label{cor:nearvertex}

\begin{enumerate}

\item Suppose that $\F$ is a face of $\LRC(2)$ of codimension $n-1$ which intersects the locus of strictly 
dominant weights.
Then for every triple $(\mu_1,\mu_2,\mu)$ of weights in $\F$, 
$\mult(\V_{\smu},\V_{\smu_1}\otimes\V_{\smu_2})\leq 1$.
After scaling (i.e., modulo the saturation problem) $\mult(\V_{\smu},\V_{\smu_1}\otimes\V_{\smu_2})=1$.

\medskip
\item If $\G$ is classical, then for every codimension $n$ face $\F'$ of cohomological components, 
there is a face $\F$ of codimension 
$\lfloor\frac{n}{2}\rfloor$ containing $\F'$ such that for every triple 
$(\mu_1,\mu_2,\mu)$ of weights in $\F$, $\mult(\V_{\smu}, \V_{\smu_1}\otimes\V_{\smu_2})\leq 1$
(and after scaling, $\mult(\V_{\smu},\V_{\smu_1}\otimes\V_{\smu_2})=1$). 
\end{enumerate}

\bpf
If $\F$ is a codimension $n-1$ face which intersects the locus of strictly dominant weights then 
by Theorem \ref{thm:Ressayre}({\em b}) $\F$ is of the type in Theorem 
\ref{thm:Ressayre}({\em a}) with $\I$ a single simple root.  Since the partition function $\kostP_{\!\delpos_{\I}}$
is the constant function $1$ for every $\gamma$ in the span of $\I$, part ({\em a}) follows from Corollary
\ref{cor:boundarybound}.

\np
For ({\em b}), pick
$\I$ to be the largest set of simple roots such that the subdiagram of the Dynkin diagram of $\G$ corresponding
to $\I$ consists of disconnected $\A_1$ pieces.  
Then the partition function $\kostP_{\!\delpos_{\I}}$ is also $1$ for every $\gamma$ in the span of $\delpos_{\I}$.

\np
Let $w_1$, $w_2$, and $w$ be the elements defining $\F'$, and $\overline{w}_1$, $\overline{w}_2$, and $\overline{w}$ the
elements of shortest length in the corresponding left cosets of $\W_{\P_{\I}}$.  
By Lemma \ref{lem:unproved}, $\overline{w}_1$, $\overline{w}_2$, and $\overline{w}$ 
satisfy \eqref{eqn:Iconditions} with respect to $\I$.  
Let $\F$ be the face defined by this data.  Corollary  \ref{cor:boundarybound} then gives the 
multiplicity bound $1$. \epf

\np
{\bf Remark.}
In Corollary \ref{cor:nearvertex} one can do slightly better when $\G$ is of type $\D_{n}$: 
there is a face $\F$ of codimension $\lfloor\frac{n-1}{2}\rfloor$ with
this property.  For the exceptional Lie algebras, the proof above guarantees a face of codimension $1$, $2$, 
$3$, $4$, and $4$, if $\G$ is isomorphic to $\G_2$, $\F_4$, $\E_6$, $\E_7$, and $\E_8$ respectively.

\np
{\bf Remark.}
After Corollary \ref{cor:Steinberg} we noted that near the highest weight $\mu_1+\mu_2$ the bound given there is
exact.  If $\invset_{w}=\invset_{w_1}\sqcup\invset_{w_2}$ and $([\Omega_{w_1}]\cap[\Omega_{w_2}])\cdot [\X_{w}]=1$ then 
experimental evidence suggests that the bound in \eqref{eqn:cohbound} is also exact when $\mu$ is sufficiently
close to the cohomological face determined by $w_1$, $w_2$, and $w$.   In other words, that the multiplicities
near a cohomological component are given by the Kostant partition function applied to the sum 
$w_1^{-1}\mu_1+w_2^{-1}\mu_2-w^{-1}\mu$ measuring the distance between $\mu$ and that component.
Corollary \ref{cor:nearvertex} above establishes some consequences consistent with this description, but without
assuming it.   We illustrate this behaviour with the following example.

\psset{unit=0.8,framesep=1pt}
\ifthenelse{\boolean{showpicture}}
{
\noindent
\begin{centering}
\begin{tabular}{c}
\begin{pspicture}(-1.000000,0.732051)(11.000000,13.258330)
\rput(5.000000,0.4){\scriptsize {$\V_{6,5}\otimes \V_{9,7}$}}
\rput(5,12.5){\scriptsize {$\Sl_3$}}
\SpecialCoor
\pscustom[fillstyle=solid,fillcolor=superlightgray,linecolor=gray]{
\wtmvto{3}{18}
\wtlnto{0}{15}
\wtlnto{0}{9}
\wtlnto{4}{1}
\wtlnto{6}{0}
\wtlnto{18}{0}
\wtlnto{20}{2}
\wtlnto{15}{12}
\wtlnto{3}{18}
}
\Agrid{-1.000000}{0.732051}{11.000000}{12.258330}
\SpecialCoor
\psset{hatchcolor=gray}
\pscustom[linecolor=white,linewidth=2.0pt,fillstyle=hlines]{%
\wtmvto{3}{18}  
\widewtlnto{1}{16}
\widewtlnto{3}{12}
\widewtlnto{7}{10}
\widewtlnto{9}{12}
\widewtlnto{7}{16}
\widewtlnto{3}{18}
}
\bigwhitewtcirc{3}{18}
\pscustom[linecolor=white,linewidth=2.0pt,fillstyle=vlines]{%
\wtmvto{4}{1}
\widewtlnto{6}{0}
\widewtlnto{6}{0}
\widewtlnto{7}{1}
\widewtlnto{6}{3}
\widewtlnto{4}{4}
\widewtlnto{3}{3}
\widewtlnto{4}{1}
}
\bigwhitewtcirc{4}{1}
\pscustom[linecolor=white,linewidth=2.0pt,fillstyle=hlines]{%
\wtmvto{20}{2} 
\widewtlnto{18}{0}
\widewtlnto{14}{2}
\widewtlnto{12}{6}
\widewtlnto{14}{8}
\widewtlnto{18}{6}
\widewtlnto{20}{2}
}
\bigwhitewtcirc{20}{2}
\pscustom[linecolor=white,linewidth=2.0pt,fillstyle=vlines]{%
\wtmvto{15}{12} 
\widewtlnto{17}{8}
\widewtlnto{15}{6}
\widewtlnto{11}{8}
\widewtlnto{9}{12}
\widewtlnto{11}{14}
\widewtlnto{15}{12}
}
\bigwhitewtcirc{15}{12}
\pscustom[linecolor=white,linewidth=2.0pt]{%
\wtmvto{16}{7} 
\widewtlnto{14}{8}
\widewtlnto{13}{7}
}

\littlewhitewtcirc{5}{17}
\littlewhitewtcirc{7}{16}
\littlewhitewtcirc{9}{15}
\littlewhitewtcirc{11}{14}
\littlewhitewtcirc{13}{13}
\littlewhitewtcirc{16}{10}
\littlewhitewtcirc{17}{8}
\littlewhitewtcirc{18}{6}
\littlewhitewtcirc{19}{4}
\littlewhitewtcirc{19}{1}
\littlewhitewtcirc{18}{0}
\littlewhitewtcirc{15}{0}
\littlewhitewtcirc{12}{0}
\littlewhitewtcirc{9}{0}
\littlewhitewtcirc{6}{0}
\littlewhitewtcirc{3}{3}
\littlewhitewtcirc{2}{5}
\littlewhitewtcirc{1}{7}
\littlewhitewtcirc{0}{9}
\littlewhitewtcirc{0}{12}
\littlewhitewtcirc{0}{15}
\littlewhitewtcirc{1}{16}
\littlewhitewtcirc{2}{17}
\psclip{%
\pscustom[linestyle=none]{
\wtmvto{3}{18}
\wtlnto{0}{15}
\wtlnto{0}{9}
\wtlnto{4}{1}
\wtlnto{6}{0}
\wtlnto{18}{0}
\wtlnto{20}{2}
\wtlnto{15}{12}
\wtlnto{3}{18}
}
}
\pscircle[fillstyle=solid,fillcolor=superlightgray,linecolor=superlightgray](1.500000,11.258330){0.2}
\pscircle[fillstyle=solid,fillcolor=superlightgray,linecolor=superlightgray](2.500000,11.258330){0.2}
\pscircle[fillstyle=solid,fillcolor=superlightgray,linecolor=superlightgray](3.500000,11.258330){0.2}
\pscircle[fillstyle=solid,fillcolor=superlightgray,linecolor=superlightgray](4.500000,11.258330){0.2}
\pscircle[fillstyle=solid,fillcolor=superlightgray,linecolor=superlightgray](5.500000,11.258330){0.2}
\pscircle[fillstyle=solid,fillcolor=superlightgray,linecolor=superlightgray](6.500000,11.258330){0.2}
\pscircle[fillstyle=solid,fillcolor=superlightgray,linecolor=superlightgray](7.500000,11.258330){0.2}
\pscircle[fillstyle=solid,fillcolor=superlightgray,linecolor=superlightgray](1.500000,2.598076){0.2}
\pscircle[fillstyle=solid,fillcolor=superlightgray,linecolor=superlightgray](2.000000,1.732051){0.2}
\pscircle[fillstyle=solid,fillcolor=superlightgray,linecolor=superlightgray](0.000000,5.196152){0.2}
\pscircle[fillstyle=solid,fillcolor=superlightgray,linecolor=superlightgray](1.000000,3.464102){0.2}
\pscircle[fillstyle=solid,fillcolor=superlightgray,linecolor=superlightgray](2.500000,2.598076){0.2}
\pscircle[fillstyle=solid,fillcolor=superlightgray,linecolor=superlightgray](0.500000,4.330127){0.2}
\pscircle[fillstyle=solid,fillcolor=superlightgray,linecolor=superlightgray](3.000000,1.732051){0.2}
\pscircle[fillstyle=solid,fillcolor=superlightgray,linecolor=superlightgray](3.000000,3.464102){0.2}
\pscircle[fillstyle=solid,fillcolor=superlightgray,linecolor=superlightgray](2.000000,3.464102){0.2}
\pscircle[fillstyle=solid,fillcolor=superlightgray,linecolor=superlightgray](0.500000,6.062178){0.2}
\pscircle[fillstyle=solid,fillcolor=superlightgray,linecolor=superlightgray](1.500000,4.330127){0.2}
\pscircle[fillstyle=solid,fillcolor=superlightgray,linecolor=superlightgray](3.500000,2.598076){0.2}
\pscircle[fillstyle=solid,fillcolor=superlightgray,linecolor=superlightgray](1.000000,5.196152){0.2}
\pscircle[fillstyle=solid,fillcolor=superlightgray,linecolor=superlightgray](0.000000,6.928203){0.2}
\pscircle[fillstyle=solid,fillcolor=superlightgray,linecolor=superlightgray](4.500000,2.598076){0.2}
\pscircle[fillstyle=solid,fillcolor=superlightgray,linecolor=superlightgray](2.500000,4.330127){0.2}
\pscircle[fillstyle=solid,fillcolor=superlightgray,linecolor=superlightgray](3.000000,5.196152){0.2}
\pscircle[fillstyle=solid,fillcolor=superlightgray,linecolor=superlightgray](2.000000,5.196152){0.2}
\pscircle[fillstyle=solid,fillcolor=superlightgray,linecolor=superlightgray](0.000000,8.660254){0.2}
\pscircle[fillstyle=solid,fillcolor=superlightgray,linecolor=superlightgray](1.500000,6.062178){0.2}
\pscircle[fillstyle=solid,fillcolor=superlightgray,linecolor=superlightgray](3.500000,4.330127){0.2}
\pscircle[fillstyle=solid,fillcolor=superlightgray,linecolor=superlightgray](1.000000,6.928203){0.2}
\pscircle[fillstyle=solid,fillcolor=superlightgray,linecolor=superlightgray](0.500000,7.794229){0.2}
\pscircle[fillstyle=solid,fillcolor=superlightgray,linecolor=superlightgray](4.000000,3.464102){0.2}
\pscircle[fillstyle=solid,fillcolor=superlightgray,linecolor=superlightgray](3.000000,6.928203){0.2}
\pscircle[fillstyle=solid,fillcolor=superlightgray,linecolor=superlightgray](2.500000,6.062178){0.2}
\pscircle[fillstyle=solid,fillcolor=superlightgray,linecolor=superlightgray](6.000000,3.464102){0.2}
\pscircle[fillstyle=solid,fillcolor=superlightgray,linecolor=superlightgray](2.000000,6.928203){0.2}
\pscircle[fillstyle=solid,fillcolor=superlightgray,linecolor=superlightgray](3.500000,6.062178){0.2}
\pscircle[fillstyle=solid,fillcolor=superlightgray,linecolor=superlightgray](1.500000,7.794229){0.2}
\pscircle[fillstyle=solid,fillcolor=superlightgray,linecolor=superlightgray](2.000000,8.660254){0.2}
\pscircle[fillstyle=solid,fillcolor=superlightgray,linecolor=superlightgray](1.000000,8.660254){0.2}
\pscircle[fillstyle=solid,fillcolor=superlightgray,linecolor=superlightgray](4.000000,5.196152){0.2}
\pscircle[fillstyle=solid,fillcolor=superlightgray,linecolor=superlightgray](0.500000,9.526279){0.2}
\pscircle[fillstyle=solid,fillcolor=superlightgray,linecolor=superlightgray](1.500000,9.526279){0.2}
\pscircle[fillstyle=solid,fillcolor=superlightgray,linecolor=superlightgray](5.000000,3.464102){0.2}
\pscircle[fillstyle=solid,fillcolor=superlightgray,linecolor=superlightgray](4.500000,4.330127){0.2}
\pscircle[fillstyle=solid,fillcolor=superlightgray,linecolor=superlightgray](2.500000,7.794229){0.2}
\pscircle[fillstyle=solid,fillcolor=superlightgray,linecolor=superlightgray](7.000000,5.196152){0.2}
\pscircle[fillstyle=solid,fillcolor=superlightgray,linecolor=superlightgray](1.000000,10.392305){0.2}
\pscircle[fillstyle=solid,fillcolor=superlightgray,linecolor=superlightgray](4.500000,7.794229){0.2}
\pscircle[fillstyle=solid,fillcolor=superlightgray,linecolor=superlightgray](5.500000,4.330127){0.2}
\pscircle[fillstyle=solid,fillcolor=superlightgray,linecolor=superlightgray](3.500000,9.526279){0.2}
\pscircle[fillstyle=solid,fillcolor=superlightgray,linecolor=superlightgray](5.000000,5.196152){0.2}
\pscircle[fillstyle=solid,fillcolor=superlightgray,linecolor=superlightgray](5.000000,6.928203){0.2}
\pscircle[fillstyle=solid,fillcolor=superlightgray,linecolor=superlightgray](4.500000,6.062178){0.2}
\pscircle[fillstyle=solid,fillcolor=superlightgray,linecolor=superlightgray](6.500000,6.062178){0.2}
\pscircle[fillstyle=solid,fillcolor=superlightgray,linecolor=superlightgray](4.000000,6.928203){0.2}
\pscircle[fillstyle=solid,fillcolor=superlightgray,linecolor=superlightgray](5.500000,6.062178){0.2}
\pscircle[fillstyle=solid,fillcolor=superlightgray,linecolor=superlightgray](3.500000,7.794229){0.2}
\pscircle[fillstyle=solid,fillcolor=superlightgray,linecolor=superlightgray](3.000000,10.392305){0.2}
\pscircle[fillstyle=solid,fillcolor=superlightgray,linecolor=superlightgray](3.000000,8.660254){0.2}
\pscircle[fillstyle=solid,fillcolor=superlightgray,linecolor=superlightgray](6.000000,5.196152){0.2}
\pscircle[fillstyle=solid,fillcolor=superlightgray,linecolor=superlightgray](2.500000,9.526279){0.2}
\pscircle[fillstyle=solid,fillcolor=superlightgray,linecolor=superlightgray](7.500000,4.330127){0.2}
\pscircle[fillstyle=solid,fillcolor=superlightgray,linecolor=superlightgray](2.000000,10.392305){0.2}
\pscircle[fillstyle=solid,fillcolor=superlightgray,linecolor=superlightgray](6.500000,4.330127){0.2}
\pscircle[fillstyle=solid,fillcolor=superlightgray,linecolor=superlightgray](4.000000,8.660254){0.2}
\pscircle[fillstyle=solid,fillcolor=superlightgray,linecolor=superlightgray](9.500000,7.794229){0.2}
\pscircle[fillstyle=solid,fillcolor=superlightgray,linecolor=superlightgray](6.000000,6.928203){0.2}
\pscircle[fillstyle=solid,fillcolor=superlightgray,linecolor=superlightgray](9.500000,6.062178){0.2}
\pscircle[fillstyle=solid,fillcolor=superlightgray,linecolor=superlightgray](5.500000,7.794229){0.2}
\pscircle[fillstyle=solid,fillcolor=superlightgray,linecolor=superlightgray](8.500000,7.794229){0.2}
\pscircle[fillstyle=solid,fillcolor=superlightgray,linecolor=superlightgray](5.000000,8.660254){0.2}
\pscircle[fillstyle=solid,fillcolor=superlightgray,linecolor=superlightgray](6.000000,10.392305){0.2}
\pscircle[fillstyle=solid,fillcolor=superlightgray,linecolor=superlightgray](4.500000,9.526279){0.2}
\pscircle[fillstyle=solid,fillcolor=superlightgray,linecolor=superlightgray](10.000000,6.928203){0.2}
\pscircle[fillstyle=solid,fillcolor=superlightgray,linecolor=superlightgray](4.000000,10.392305){0.2}
\pscircle[fillstyle=solid,fillcolor=superlightgray,linecolor=superlightgray](6.500000,9.526279){0.2}
\pscircle[fillstyle=solid,fillcolor=superlightgray,linecolor=superlightgray](8.000000,5.196152){0.2}
\pscircle[fillstyle=solid,fillcolor=superlightgray,linecolor=superlightgray](8.000000,8.660254){0.2}
\pscircle[fillstyle=solid,fillcolor=superlightgray,linecolor=superlightgray](7.500000,6.062178){0.2}
\pscircle[fillstyle=solid,fillcolor=superlightgray,linecolor=superlightgray](7.000000,8.660254){0.2}
\pscircle[fillstyle=solid,fillcolor=superlightgray,linecolor=superlightgray](7.000000,6.928203){0.2}
\pscircle[fillstyle=solid,fillcolor=superlightgray,linecolor=superlightgray](9.000000,8.660254){0.2}
\pscircle[fillstyle=solid,fillcolor=superlightgray,linecolor=superlightgray](6.500000,7.794229){0.2}
\pscircle[fillstyle=solid,fillcolor=superlightgray,linecolor=superlightgray](7.500000,7.794229){0.2}
\pscircle[fillstyle=solid,fillcolor=superlightgray,linecolor=superlightgray](6.000000,8.660254){0.2}
\pscircle[fillstyle=solid,fillcolor=superlightgray,linecolor=superlightgray](7.500000,9.526279){0.2}
\pscircle[fillstyle=solid,fillcolor=superlightgray,linecolor=superlightgray](5.500000,9.526279){0.2}
\pscircle[fillstyle=solid,fillcolor=superlightgray,linecolor=superlightgray](8.000000,6.928203){0.2}
\pscircle[fillstyle=solid,fillcolor=superlightgray,linecolor=superlightgray](5.000000,10.392305){0.2}
\pscircle[fillstyle=solid,fillcolor=superlightgray,linecolor=superlightgray](7.000000,10.392305){0.2}
\pscircle[fillstyle=solid,fillcolor=superlightgray,linecolor=superlightgray](9.000000,5.196152){0.2}
\pscircle[fillstyle=solid,fillcolor=superlightgray,linecolor=superlightgray](8.500000,6.062178){0.2}
\pscircle[fillstyle=solid,fillcolor=superlightgray,linecolor=superlightgray](9.000000,6.928203){0.2}
\pscircle[fillstyle=solid,fillcolor=superlightgray,linecolor=superlightgray](8.500000,9.526279){0.2}
\pscircle[fillstyle=solid,fillcolor=superlightgray,linecolor=superlightgray](8.000000,10.392305){0.2}
\endpsclip
\bigwhitewtcirc{15}{12}
\bigwhitewtcirc{20}{2}
\bigwhitewtcirc{4}{1}
\bigwhitewtcirc{3}{18}
\rput(1.500000,11.258330){\scriptsize {$1$}}
\rput(2.500000,11.258330){\scriptsize {$1$}}
\rput(3.500000,11.258330){\scriptsize {$1$}}
\rput(4.500000,11.258330){\scriptsize {$1$}}
\rput(5.500000,11.258330){\scriptsize {$1$}}
\rput(6.500000,11.258330){\scriptsize {$1$}}
\rput(7.500000,11.258330){\scriptsize {$1$}}
\rput(2.000000,3.464102){\scriptsize {$2$}}
\rput(2.000000,1.732051){\scriptsize {$1$}}
\rput(0.000000,5.196152){\scriptsize {$1$}}
\rput(1.500000,2.598076){\scriptsize {$1$}}
\rput(2.500000,2.598076){\scriptsize {$2$}}
\rput(1.000000,3.464102){\scriptsize {$1$}}
\rput(0.500000,4.330127){\scriptsize {$1$}}
\rput(3.000000,1.732051){\scriptsize {$1$}}
\rput(1.000000,5.196152){\scriptsize {$2$}}
\rput(1.500000,4.330127){\scriptsize {$2$}}
\rput(3.500000,2.598076){\scriptsize {$2$}}
\rput(0.500000,6.062178){\scriptsize {$2$}}
\rput(2.500000,4.330127){\scriptsize {$3$}}
\rput(0.000000,6.928203){\scriptsize {$1$}}
\rput(3.000000,3.464102){\scriptsize {$3$}}
\rput(4.000000,3.464102){\scriptsize {$3$}}
\rput(2.000000,5.196152){\scriptsize {$3$}}
\rput(0.500000,7.794229){\scriptsize {$2$}}
\rput(1.500000,6.062178){\scriptsize {$3$}}
\rput(4.500000,2.598076){\scriptsize {$1$}}
\rput(1.000000,6.928203){\scriptsize {$3$}}
\rput(0.000000,8.660254){\scriptsize {$1$}}
\rput(0.500000,9.526279){\scriptsize {$1$}}
\rput(3.500000,4.330127){\scriptsize {$4$}}
\rput(4.000000,5.196152){\scriptsize {$5$}}
\rput(3.000000,5.196152){\scriptsize {$4$}}
\rput(1.000000,8.660254){\scriptsize {$2$}}
\rput(2.500000,6.062178){\scriptsize {$4$}}
\rput(4.500000,4.330127){\scriptsize {$4$}}
\rput(2.000000,6.928203){\scriptsize {$4$}}
\rput(1.500000,7.794229){\scriptsize {$3$}}
\rput(5.000000,3.464102){\scriptsize {$2$}}
\rput(5.000000,6.928203){\scriptsize {$6$}}
\rput(3.500000,6.062178){\scriptsize {$5$}}
\rput(2.500000,9.526279){\scriptsize {$3$}}
\rput(3.000000,6.928203){\scriptsize {$5$}}
\rput(6.500000,4.330127){\scriptsize {$2$}}
\rput(2.500000,7.794229){\scriptsize {$4$}}
\rput(3.000000,8.660254){\scriptsize {$4$}}
\rput(2.000000,8.660254){\scriptsize {$3$}}
\rput(4.500000,7.794229){\scriptsize {$5$}}
\rput(1.500000,9.526279){\scriptsize {$2$}}
\rput(3.500000,7.794229){\scriptsize {$5$}}
\rput(6.000000,3.464102){\scriptsize {$1$}}
\rput(6.000000,5.196152){\scriptsize {$4$}}
\rput(1.000000,10.392305){\scriptsize {$1$}}
\rput(4.000000,6.928203){\scriptsize {$6$}}
\rput(5.500000,4.330127){\scriptsize {$3$}}
\rput(5.500000,6.062178){\scriptsize {$5$}}
\rput(5.000000,5.196152){\scriptsize {$5$}}
\rput(4.500000,6.062178){\scriptsize {$6$}}
\rput(2.000000,10.392305){\scriptsize {$2$}}
\rput(9.000000,6.928203){\scriptsize {$2$}}
\rput(4.000000,8.660254){\scriptsize {$4$}}
\rput(8.500000,6.062178){\scriptsize {$2$}}
\rput(3.500000,9.526279){\scriptsize {$3$}}
\rput(7.500000,7.794229){\scriptsize {$3$}}
\rput(3.000000,10.392305){\scriptsize {$2$}}
\rput(9.000000,5.196152){\scriptsize {$1$}}
\rput(7.500000,4.330127){\scriptsize {$1$}}
\rput(9.500000,6.062178){\scriptsize {$1$}}
\rput(7.000000,5.196152){\scriptsize {$3$}}
\rput(5.000000,10.392305){\scriptsize {$2$}}
\rput(6.500000,6.062178){\scriptsize {$4$}}
\rput(7.000000,8.660254){\scriptsize {$3$}}
\rput(6.000000,6.928203){\scriptsize {$5$}}
\rput(5.500000,9.526279){\scriptsize {$3$}}
\rput(5.500000,7.794229){\scriptsize {$5$}}
\rput(8.500000,7.794229){\scriptsize {$2$}}
\rput(5.000000,8.660254){\scriptsize {$4$}}
\rput(6.000000,8.660254){\scriptsize {$4$}}
\rput(4.500000,9.526279){\scriptsize {$3$}}
\rput(6.500000,9.526279){\scriptsize {$3$}}
\rput(4.000000,10.392305){\scriptsize {$2$}}
\rput(6.500000,7.794229){\scriptsize {$4$}}
\rput(8.000000,5.196152){\scriptsize {$2$}}
\rput(6.000000,10.392305){\scriptsize {$2$}}
\rput(7.500000,6.062178){\scriptsize {$3$}}
\rput(7.000000,6.928203){\scriptsize {$4$}}
\rput(8.000000,6.928203){\scriptsize {$3$}}
\rput(8.000000,8.660254){\scriptsize {$2$}}
\rput(7.500000,9.526279){\scriptsize {$2$}}
\rput(7.000000,10.392305){\scriptsize {$2$}}
\rput(10.000000,6.928203){\scriptsize {$1$}}
\rput(9.500000,7.794229){\scriptsize {$1$}}
\rput(9.000000,8.660254){\scriptsize {$1$}}
\rput(8.500000,9.526279){\scriptsize {$1$}}
\rput(8.000000,10.392305){\scriptsize {$1$}}
\biggerwtsq{2}{11}
\biggerwtsq{10}{4}
\end{pspicture}\\
\\
Figure 2 \\
\end{tabular}\\
\end{centering}
}
{}

\np
In the figure above we have outlined hexagons at each cohomological component 
where the multiplicities behave as proscribed by the Kostant partition function.

\np
\bpoint{Order of growth of partition functions}
Suppose that we are given a matrix $\A$ with integer entries.  The {\em vector partition function}
corresponding to $\A$ counts, for each $\gamma$ in the span of the columns of $\A$, the number of vectors $v$
with non-negative integer entries such that $\A v=\gamma$. 

\np
Both 
the Kostant partition function and its extension to multisets above are examples of vector partition functions.
For the partition function $\kostP_{\!\Ssh}$ one takes
the matrix $\A$ whose rows are indexed by the simple roots and whose columns are indexed
by the elements of $\Ssh$ (repeated according to their multiplicities).  The column corresponding to an entry
of $\Ssh$ is its expression in terms of the simple roots.

\np
A result of Blakley \cite{Bl} (see also \cite[Theorem 1]{Stu}) shows that for any vector partition function
one can divide the space spanned by the 
columns of $\A$ into chambers so that the partition function is quasi-poly\-nomial on each chamber of 
degree at most $\#(\mbox{columns of $\A$})-\rank \A$.
When $\Ssh=(k-1)\delpos_{\I}$ this structure theorem for partition functions combined with 
Corollary \ref{cor:boundarybound} gives the following result.

\tpoint{Corollary} \label{cor:boundarygrowth}
Suppose that $(\mu_1,\ldots, \mu_k,\mu)\in\LRC(k)$ is on a face of the type from Theorem 
\ref{thm:Ressayre}({\em a}).  Then the function 

\begin{equation}\label{eqn:hilbfunc}
m\mapsto \mult(\V_{\!m\mu},\V_{\!m\mu_1}\otimes\cdots\otimes\V_{\!m\mu_k})
\end{equation}

\np 
is bounded above by a polynomial in $m$ of degree $(k-1)\#\delpos_{\I}-\#\I$.

\np
{\bf Remark.} The function 
\eqref{eqn:hilbfunc} 
can also be viewed as the Hilbert function of an appropriate GIT quotient of $\X^{k+1}$ (or some product of 
$\G/\P$'s if not every weight is strictly dominant).  Corollary \ref{cor:boundarygrowth} can therefore also be
viewed as a bound on the dimension of the GIT quotients when the line bundle giving the GIT quotient comes from
a face of $\LRC(k)$.

\np
{\bf Remark.} Heckman established \cite[(3.17)]{he} that the asymptotic growth of the multiplicity function 
$$(\mu,\mu_1,\ldots,\mu_k)\mapsto\mult(\V_{\smu},\V_{\smu_1}\otimes\cdots\otimes \V_{\smu_k}),$$ 
when viewed as a function on $\Lambda^{k+1}_{+}$, is polynomial of degree $(k-1)\#\delpos-n$.  
Equation \eqref{eqn:cohbound} provides a family of bounds (one for each $w_1$,\ldots, $w_k$, and $w$ 
satisfying \eqref{eqn:liningup}) each of which has the same order of growth as the actual multiplicity function.
It would be interesting to compare the leading terms given by \eqref{eqn:cohbound} with the ones in \cite{he}.

\renewcommand{\thesection}{\Alph{section}}
\setcounter{section}{0}
\renewenvironment{equation}{\medskip\noindent\refstepcounter{subsection}\makebox[0pt][l]{({\bf\thesubsection})}\begin{minipage}[b]{\textwidth}$$}{$$\end{minipage}\medskip\noindent}
\renewcommand{\tpoint}[1]{\np\vspace{3mm}\par\refstepcounter{subsection}\noindent{\em #1 {\rm(}{\em \thesubsection}{\rm)} ---} }

\section{Appendix on invariants of $\Sym^{\cdot} \ggg$} 
\label{sec:appendix}

\renewcommand{\labelenumi}{{\bf\oldstylenums{\arabic{enumi}}.}}
\newcommand{\old}{\oldstylenums}

\np
Let $\gs$ be a subalgebra of $\ggg$ which is the reductive part of a parabolic subalgebra of $\ggg$ and
let $\gt$ be the center of $\gs$.  
Following Kostant, \cite{K}, for any $\nu \in \gt^* \setminus\{0\}$  we set

$$\ggg^\nu = \left\{{x \in \ggg \, \st \, [t, x] = \nu(t) x \,\,\mbox{for all}\, t \in \gt\rule{0cm}{0.39cm}}\right\}.$$

\np
An element $\nu \in \gt^* \setminus\{0\}$ is called an $\gt$--root if $\ggg^{\nu}\neq 0$. 
The properties of $\gt$--roots and the corresponding $\gt$--root decomposition of $\ggg$
inherit most of the properties of the usual roots and the usual root decomposition of $\ggg$, as established by
Kostant. We refer the reader to \cite{K} for a complete account of the theory of $\gt$--roots.

\np
Let $\Rsh$ denote the set of $\gt$--roots.
A subset $\Ssh$ of $\Rsh$ is called {\em saturated} if 
whenever $\nu \in \Ssh$ and $r \nu \in \Rsh$ for some $r \in \QQ_+$ then $r \nu \in \Ssh$ as well. 
With this notation we prove:

\tpoint{Theorem}\label{thm:classical} 
Let $\ggg$ be a simple classical Lie algebra, $\gs$ be a subalgebra which is the reductive part of a parabolic
subalgebra of $\ggg$, $\gt$ be the center of $\gs$, 
$\Ssh$ be a saturated subset of the $\gt$-roots of $\ggg$, 
and $\Msh = \oplus_{\nu \in \Ssh} \ggg^\nu$. If 
$\dim(\Sym^\cdot(\Msh))^{\gs} = 1$, then there exists a
parabolic subalgebra $\gp_{\Msh}$ of $\ggg$ with reductive part $\gs$ such that $\Msh \subset \gp_{\Msh}$.

\np
Theorem \ref{thm:classical} also holds if $\ggg=\G_2$; 
we will give examples to show that it does not hold for the other exceptional algebras.

\np
First we describe the parabolic subalgebras and the corresponding sets $\Rsh$ for the classical Lie algebras. 
For convenience of notation we will work with the reductive Lie algebra $gl_n$ instead of $sl_n$.
A partition $\cP = \{\I_1, \ldots, \I_k\}$ of $\{1, \ldots, n\}$ is {\it linearly ordered} if the set of its parts is linearly ordered.
We write $\cP(i)$ for the part of $\cP$ which contains $i$. The inequalities
$\cP(i) \prec \cP(j)$ and $\cP(i) \preceq \cP(j)$ are taken in the linear order of the parts of $\cP$.
For the standard basis $\{ \vep_1, \ldots, \vep_n\}$ of $\gh^*$ we denote the dual basis of $\gh$ by
$\{h_1, \ldots, h_n\}$. 
A linear order on the set $\{ \pm \delta_1, \ldots, \pm \delta_k\}$ is {\it compatible with multiplication by $-1$} if 
$\pm \delta_i \prec \pm \delta_j$ implies $\mp \delta_j \prec \mp \delta_i$.
To simplify notation we adopt the convention that
$\B_1$, respectively $\C_1$, is a subalgebra of $\ggg = \B_n$, respectively $\ggg = \C_n$, isomorphic to 
$\A_1$ and whose roots are short, respectively long roots, of $\ggg$. The subalgebras $\D_2=\A_1\oplus\A_1$ and 
$\D_3=\A_3$ of $\D_n$ have similar meaning.  

\np
Let $\ggg$ be of type $\X_{n}=\A_n$, $\B_n$, $\C_n$, or $\D_n$ and let $\gs$ be a subalgebra of $\ggg$ which 
is the reductive part of a parabolic subalgebra of $\ggg$.  Every simple ideal of $\gs$ is isomorphic to $\A_r$
or $\X_r$ for some $r$.  Furthermore, if $\ggg$ is not of type $\A_n$, $\gs$ has at most 
one simple ideal of type $\X_r$. 
For $\ggg$ of type $\X_n=\B_n$, $\C_n$, or $\D_n$ the parabolic subalgebras of $\ggg$
are split into two types depending on whether their reductive parts contain (Type II) or do not contain (Type I)
a simple ideal of type $\X_r$ (including $\B_1$, $\C_1$, $\D_2$, or $\D_3$).

\np
In the description of the combinatorics of the simple classical Lie algebras below,
the formulas for their parabolic subalgebras $\gp$ containing a fixed reductive part $\gs$ look very uniform
(e.g. \old{11}).
In some instances this is misleading since the formulas do not explicitly indicate the subalgebra $\gs$ which, 
however, is an integral part of the structure of $\gp$.

\np
{\bf Roots and parabolic subsets of the classical Lie algebras.} 
We now list the combinatorial descriptions of the parabolic subalgebras and related data in the classical cases.

\noindent
\underline{$\ggg=\A_{n}$}

\begin{enumerate}
\item The roots of $\ggg$ are: $ \Delta = \{ \vep_i - \vep_j \, | \, 1 \leq i \neq j \leq n \} $

\item
Parabolic subalgebras of $\ggg$ are in one-to-one correspondence with: \\

\begin{centering}
linearly ordered partitions $\cP = (\I_1, \ldots, \I_k)$ of $\{1,\ldots,n\}$. \\
\end{centering}
\end{enumerate}

\np
Given a linearly ordered partition $\cP$,

\begin{enumerate}
\setcounter{enumi}{2}
\item The roots of $\gp_\cP$ are $\{ \vep_i - \vep_j \, | \, i \neq j, \cP(i) \preceq \cP(j) \}$  
\item The roots of $\gs_\cP$ are $\{ \vep_i - \vep_j \, | \, i \neq j, \cP(i) = \cP(j) \}$
\item $\gs_\cP = \oplus_i \, \gs_\cP^i$, where $\gs_\cP^i \cong gl_{|I_i|}$; 
\item The Cartan subalgebra of $\gs_\cP^i$ is spanned by $\{h_j\}_{j \in \I_i}$  
\item The roots of $\gs_\cP^i$ are $\{\vep_j - \vep_l \, | \, j \neq l \in \I_i\}$.   
\item $\gt_\cP$ has a basis $\{t_1, \ldots, t_k\}$ with $t_i = \frac{1}{|I_i|} \sum_{j \in I_i} h_j$ 
\end{enumerate}

If $\{\delta_1, \ldots, \delta_k\}$ is the basis of $\gt^*$ dual to $\{t_1, \ldots, t_k\}$ then

\begin{enumerate}
\setcounter{enumi}{8}
\item $\Rsh = \{ \delta_i - \delta_j \, | \, 1 \leq i \neq j \leq k\}.$
\item For $\nu = \delta_i - \delta_j \in \Rsh$, $\gs_\cP$--module $\ggg^\nu\cong \V_i \otimes \V_j^*$, 
where $\V_i$ and $\V_j^*$ are the natural $\gs_\cP^i$-module and the dual of the natural $\gs_\cP^j$-module respectively,  all other factors of $\gs_{\cP}$ acting trivially.

\medskip
\item
The parabolic subalgebras of $\ggg$ whose reductive part is $\gs_\cP$ are in a bijection with the ordered partitions $\cQ$ of
$\{1, \ldots, n\}$ whose parts are the same as the parts of $\cP$ or, equivalently, with linear orders on the set $\{\delta_1, \ldots,
\delta_k\}$. 
\end{enumerate}

\bigskip
\noindent
\underline{$\ggg=\B_{n}$}

\begin{enumerate}
\item The roots of $\ggg$ are: 
$\Delta = \{ \pm \vep_i \pm \vep_j, \pm \vep_i  \, | \, 1 \leq i \neq j \leq n \}.$

\item
Parabolic subalgebras of $\ggg$ are in one-to-one correspondence with: \\

\begin{centering}
{\bf Type I:} 
\begin{minipage}[t]{0.75\textwidth}
pairs $(\cP, \sigma)$, where $\cP = (\I_1, \ldots, \I_k)$ is a linearly ordered partition 
of $\{1, \ldots, n\}$ and $\sigma\colon \{1, \ldots, n \} \to \{\pm 1\}$ is a choice of signs. 
\end{minipage} \\

\bigskip
{\bf Type II:} 
\begin{minipage}[t]{0.75\textwidth}
pairs $(\cP, \sigma)$, where $\cP = (\I_0, \I_1, \ldots, \I_k)$ is a linearly ordered partition 
of $\{1, \ldots, n\}$ with largest element $\I_0$ and 
$\sigma\colon \{1, \ldots, n \} \backslash \I_0 \to \{\pm 1\}$ is a choice of signs. 
\end{minipage}\\
\end{centering}
\end{enumerate}

\np
\underline{In Type I:}

\begin{enumerate}
\setcounter{enumi}{2}
\item The roots of $\gp_{(\cP,\sigma)}$ are
$$\left\{{ \sigma(i) \vep_i - \sigma(j) \vep_j \, \st \, i \neq j, \cP(i) \preceq \cP(j) }\right\} \cup \{\sigma(i) \vep_i + \sigma(j) \vep_j, \sigma(i) \vep_i \, \st \, i \neq j\}$$
\item The roots of $\gs_{(\cP,\sigma)}$ are $\{ \sigma(i)\vep_i - \sigma(j)\vep_j \, | \, i \neq j, \cP(i) = \cP(j) \}$.
\item $\gs_{(\cP,\sigma)} = \oplus_i \, \gs_{(\cP,\sigma)}^i$, where $\gs_{(\cP,\sigma)}^i \cong gl_{|\I_i|}$.
\item The Cartan subalgebra of $\gs_\cP^i$ is spanned by $\{\sigma(j) h_j\}_{j \in \I_i}$ 
\item The roots of $\gs_\cP^i$ are $\{\sigma(j)\vep_j - \sigma(l)\vep_l \, | \, j \neq l \in \I_i\}$.  
\item 
$\gt_{(\cP,\sigma)}$ has a basis $\{t_1, \ldots, t_k\}$ with $t_i = \frac{1}{|\I_i|} \sum_{j \in \I_i} \sigma(j) h_j$.
\end{enumerate}

If $\{\delta_1, \ldots, \delta_k\}$ the basis of $\gt^*$ dual to $\{t_1, \ldots, t_k\}$ then

\begin{enumerate}
\setcounter{enumi}{8}
\item 
$\Rsh = \{ \pm \delta_i \pm \delta_j, \pm \delta_i \, | \, 1 \leq i \neq j \leq k\} \cup \{\pm 2 \delta_i \, | \, |\I_i|>1\}.$
\item For $\nu \in \Rsh$, 

\begin{enumerate}
\item $\ggg^\nu \cong \V_i^\pm \otimes \V_j^\pm$ if $\nu = \pm \delta_i \pm \delta_j$,
\item $\ggg^\nu\cong\V_i^\pm$ if $\nu = \pm \delta_i$, and 
\item $\ggg^\nu\cong\Lambda^2 \V_i^\pm$ if $\nu = \pm 2 \delta_i$, 
\end{enumerate}

\np
where $\V_i^+$ and $\V_i^-$ respectively are the natural $\gs_{(\cP,\sigma)}^i$-module and its dual, and
all other factors of $\gs_{(\cP,\sigma)}$ act trivially.

\medskip
\item
The parabolic subalgebras of $\ggg$ whose reductive part is $\gs_{\cP,\sigma}$ are in a bijection with the pairs $(\cQ, \tau)$ 
such that the parts of $\cQ$ are the same as the parts of $\cP$ and $\sigma_{|\I_i} = \pm \tau_{|\I_i}$ for every part $\I_i$ or, 
equivalently, with linear orders on the set $\{\pm \delta_1, \ldots, \pm \delta_k\}$ compatible with multiplication by $-1$.
\end{enumerate}

\np
\underline{In Type II:}

\begin{enumerate}
\setcounter{enumi}{2}
\item The roots of $\gp_{(\cP,\sigma)}$ are 

\medskip
$
\{ \sigma(i) \vep_i - \sigma(j) \vep_j \, | \,  i \neq j, \cP(i) \preceq \cP(j) \prec \I_0\} 
\cup \{\pm \vep_i \pm \vep_j, \pm \vep_i \, | \, i \neq j \in \I_0\} $ \hfill

\medskip
\hfill $\cup
\{\sigma(i) \vep_i + \sigma(j) \vep_j, \sigma(i) \vep_i \, | \, i \neq j \not \in \I_0\} \cup 
\{\sigma(i) \vep_i \pm \vep_j \, | \, i \not \in \I_0, j \in \I_0\}$ 

\item The roots of $\gs_{(\cP,\sigma)}$ are
$$\{ \sigma(i)\vep_i - \sigma(j)\vep_j \, | \, i \neq j, \cP(i) = \cP(j) \prec \I_0 \} \cup \{\pm \vep_i \pm \vep_j, \pm \vep_i \, | \, i \neq j \in \I_0\}.$$
\item $\gs_{(\cP,\sigma)} = \oplus_i \, \gs_{(\cP,\sigma)}^i$, 
where $\gs_{(\cP,\sigma)}^0 \cong \B_{|\I_0|}$ and $\gs_{(\cP,\sigma)}^i \cong gl_{|\I_i|}$ for $i>0$.
\item The Cartan subalgebra of $\gs_\cP^i$ is spanned by $\{h_j\}_{j \in \I_0}$ for $i=0$ and 
$\{\sigma(j) h_j\}_{j \in \I_i}$ for $i>0$. 
\item The roots of $\gs_\cP^i$ are $\{\pm \vep_j \pm \vep_l, \pm \vep_j \, | \, j \neq l \in \I_0\}$ for $i = 0$ and
$\{\sigma(j)\vep_j - \sigma(l)\vep_l \, | \, j \neq l \in \I_i\}$ for $i>0$.  
\item 
$\gt_{(\cP,\sigma)}$ has a basis $\{t_1, \ldots, t_k\}$ with $t_i = \frac{1}{|\I_i|} \sum_{j \in \I_i} \sigma(j) h_j$.
\end{enumerate}

\np
If $\{\delta_1, \ldots, \delta_k\}$ is the basis of $\gt^*$ dual to $\{t_1, \ldots, t_k\}$ then

\begin{enumerate}
\setcounter{enumi}{8}
\item 
$\Rsh = \{ \pm \delta_i \pm \delta_j, \pm \delta_i \, | \, 1 \leq i \neq j \leq k\} \cup \{\pm 2 \delta_i \, | \, |\I_i|>1\}.$
\item For $\nu \in \Rsh$, 

\begin{enumerate}
\item $\ggg^\nu\cong \V_i^\pm \otimes \V_j^\pm$ if $\nu = \pm \delta_i \pm \delta_j$, 
\item $\ggg^\nu\cong \V_i^\pm \otimes \V_0$ if $\nu = \pm \delta_i$,  and 
\item $\ggg^\nu\cong\Lambda^2 \V_i^\pm$ if $\nu = \pm 2 \delta_i$ 
\end{enumerate}

\np
where
$\V_i^+$ and $\V_i^-$ denote the natural $\gs_{(\cP,\sigma)}^i$-module and its dual respectively for $i>0$,
$\V_0$ denotes the natural $\gs_{(\cP,\sigma)}^0$-module,
and all other factors of $\gs_{(\cP,\sigma)}$ act trivially.

\medskip
\item
The parabolic subalgebras of $\ggg$ whose reductive part is $\gs_{\cP,\sigma}$ are in a bijection with the pairs $(\cQ, \tau)$ 
such that the parts of $\cQ$ are the same as the parts of $\cP$ and $\sigma_{|\I_i} = \pm \tau_{|\I_i}$ for every part $\I_i$ or, 
equivalently, with linear orders on the set $\{\pm \delta_1, \ldots, \pm \delta_k\}$ compatible with multiplication by $-1$. 
\end{enumerate}

\bigskip
\noindent
\underline{$\ggg=\C_{n}$}

\begin{enumerate}
\item The roots of $\ggg$ are: 
$\Delta = \{ \pm \vep_i \pm \vep_j, \pm 2 \vep_i  \, | \, 1 \leq i \neq j \leq n \}.$

\item
Parabolic subalgebras of $\ggg$ are in one-to-one correspondence with: \\

\begin{centering}
{\bf Type I:} 
\begin{minipage}[t]{0.75\textwidth}
pairs $(\cP, \sigma)$, where $\cP = (\I_1, \ldots, \I_k)$ is a linearly ordered partition 
of $\{1, \ldots, n\}$ and $\sigma\colon \{1, \ldots, n \} \to \{\pm 1\}$ is a choice of signs. 
\end{minipage} \\

\bigskip
{\bf Type II:} 
\begin{minipage}[t]{0.75\textwidth}
pairs $(\cP, \sigma)$, where $\cP = (\I_0, \I_1, \ldots, \I_k)$ is a linearly ordered partition 
of $\{1, \ldots, n\}$ with largest element $\I_0$ and 
$\sigma\colon \{1, \ldots, n \} \backslash \I_0 \to \{\pm 1\}$ is a choice of signs. 
\end{minipage}\\
\end{centering}
\end{enumerate}

\np
\underline{In Type I:}

\begin{enumerate}
\setcounter{enumi}{2}
\item The roots of $\gp_{(\cP,\sigma)}$  are
$$\{ \sigma(i) \vep_i - \sigma(j) \vep_j \, | \, i \neq j, \cP(i) \preceq \cP(j) \} \cup \{\sigma(i) \vep_i + \sigma(j) \vep_j, 2 \sigma(i) \vep_i \, | \, i \neq j\}$$
\item The roots of $\gs_{(\cP,\sigma)}$ are $\{ \sigma(i)\vep_i - \sigma(j)\vep_j \, | \, i \neq j, \cP(i) = \cP(j) \}$.
\item $\gs_{(\cP,\sigma)} = \oplus_i \, \gs_{(\cP,\sigma)}^i$, where $\gs_{(\cP,\sigma)}^i \cong gl_{|\I_i|}$.
\item The Cartan subalgebra of $\gs_\cP^i$ is spanned by $\{\sigma(j) h_j\}_{j \in \I_i}$.
\item The roots of $\gs_{(\cP,\sigma)}^i$ are $\{\sigma(j)\vep_j - \sigma(l)\vep_l \, | \, j \neq l \in \I_i\}$.  
\item 
$\gt_{(\cP,\sigma)}$ has a basis $\{t_1, \ldots, t_k\}$ with $t_i = \frac{1}{|\I_i|} \sum_{j \in \I_i} \sigma(j) h_j$.
\end{enumerate}

\noindent
If $\{\delta_1, \ldots, \delta_k\}$ is the basis of $\gt^*$ dual to $\{t_1, \ldots, t_k\}$ then

\begin{enumerate}
\setcounter{enumi}{8}
\item 
$\Rsh = \{ \pm \delta_i \pm \delta_j, \pm 2 \delta_i \, | \, 1 \leq i \neq j \leq k\}.$
\item For $\nu \in \Rsh$, 

\begin{enumerate}
\item $\ggg^\nu\cong \V_i^\pm \otimes \V_j^\pm$ if  $\nu = \pm \delta_i \pm \delta_j$.
\item $\ggg^\nu\cong\Sym^2 \V_i^\pm$ if for $\nu = \pm 2 \delta_i$.
\end{enumerate}

\np
where $\V_i^+$ and $\V_i^-$ are the natural $\gs_{(\cP,\sigma)}^i$-module and its dual, and
all other factors of $\gs_{(\cP,\sigma)}$ act trivially.

\medskip
\item
The parabolic subalgebras of $\ggg$ whose reductive part is $\gs_{\cP,\sigma}$ are in a bijection with the pairs $(\cQ, \tau)$ 
such that the parts of $\cQ$ are the same as the parts of $\cP$ and $\sigma_{|\I_i} = \pm \tau_{|\I_i}$ for every part $\I_i$ or, 
equivalently, with linear orders on the set $\{\pm \delta_1, \ldots, \pm \delta_k\}$ compatible with multiplication by $-1$.
\end{enumerate}

\np
\underline{In Type II:}

\begin{enumerate}
\setcounter{enumi}{2}
\item The roots of $\gp_{(\cP,\sigma)}$ are

\medskip
$\{ \sigma(i) \vep_i - \sigma(j) \vep_j \, | \,  i \neq j, \cP(i) \preceq \cP(j) \prec \I_0\} 
\cup \{\pm \vep_i \pm \vep_j, \pm 2 \vep_i \, | \, i \neq j \in \I_0\} \cup$ 

\medskip
\hfill $
\{\sigma(i) \vep_i + \sigma(j) \vep_j, \sigma(i) 2 \vep_i \, | \, i \neq j \not \in \I_0\} \cup 
\{\sigma(i) \vep_i \pm \vep_j \, | \, i \not \in \I_0, j \in \I_0\}$ 

\medskip
\item The roots of $\gs_{(\cP,\sigma)}$ are

\medskip
\begin{centering}
$\{ \sigma(i)\vep_i - \sigma(j)\vep_j \, | \, i \neq j, \cP(i) = \cP(j) \prec \I_0 \} \cup \{\pm \vep_i \pm \vep_j, \pm 2 \vep_i \, | \, i \neq j \in \I_0\}.$ \\
\end{centering}
\medskip
\item 
$\gs_{(\cP,\sigma)} = \oplus_{i=0}^k \, \gs_{(\cP,\sigma)}^i$, where $\gs_{(\cP,\sigma)}^0 \cong \C_{|\I_0|}$ 
and $\gs_{(\cP,\sigma)}^i \cong gl_{|\I_i|}$ for $i>0$.
\item The Cartan subalgebra of $\gs_{(\cP,\sigma)}^i$ is spanned by $\{h_j\}_{j \in \I_0}$ for $i=0$ and 
$\{\sigma(j) h_j\}_{j \in \I_i}$ for $i>0$. 
\item The roots of $\gs_{(\cP,\sigma)}^i$ are
$\{\pm \vep_j \pm \vep_l, \pm 2 \vep_j \, | \, j \neq l \in \I_0\}$ for $i = 0$ and
$\{\sigma(j)\vep_j - \sigma(l)\vep_l \, | \, j \neq l \in \I_i\}$ for $i>0$.  
\item 
$\gt_{(\cP,\sigma)}$ has a basis $\{t_1, \ldots, t_k\}$ with $t_i = \frac{1}{|\I_i|} \sum_{j \in \I_i} \sigma(j) h_j$ 
\end{enumerate}

\noindent
If $\{\delta_1, \ldots, \delta_k\}$ is the basis of $\gt^*$ dual to $\{t_1, \ldots, t_k\}$ then

\begin{enumerate}
\setcounter{enumi}{8}
\item $\Rsh = \{ \pm \delta_i \pm \delta_j, \pm \delta_i, \pm 2 \delta_i \, | \, 1 \leq i \neq j \leq k\}$.
\item For $\nu \in \Rsh$, 

\begin{enumerate}
\item $\ggg^\nu\cong\V_i^\pm \otimes \V_j^\pm$ if $\nu = \pm \delta_i \pm \delta_j$,
\item $\ggg^\nu\cong\V_i^\pm \otimes \V_0$ if $\nu = \pm \delta_i$,
\item $\ggg^\nu$ is isomorphic to $\Sym^2 \V_i^\pm$ if $\nu = \pm 2 \delta_i$,
\end{enumerate}

\np
where $\V_i^+$ and $\V_i^-$ denote the natural $\gs_{(\cP,\sigma)}^i$-module and its dual 
for $i>0$, 
$\V_0$ is the natural $\gs_{(\cP,\sigma)}^0$-module, and where 
all other factors of $\gs_{(\cP,\sigma)}$ act trivially.

\medskip
\item
The parabolic subalgebras of $\ggg$ whose reductive part is $\gs_{\cP,\sigma}$ are in a bijection 
with the pairs $(\cQ, \tau)$ such that the parts of $\cQ$ are the same as the parts 
of $\cP$ and $\sigma_{|\I_i} = \pm \tau_{|\I_i}$ for every part $\I_i$ or, equivalently, with linear orders on 
the set $\{\pm \delta_1, \ldots, \pm \delta_k\}$ compatible with multiplication by $-1$. 
\end{enumerate}

\bigskip
\noindent
\underline{$\ggg=\D_{n}$}

\begin{enumerate}
\item The roots of $\ggg$ are: 
$\Delta = \{ \pm \vep_i \pm \vep_j,  \, | \, 1 \leq i \neq j \leq n \}.$

\item
Parabolic subalgebras of $\ggg$ are determined by: \\

\begin{centering}
{\bf Type I:} 
\begin{minipage}[t]{0.75\textwidth}
pairs $(\cP, \sigma)$, where $\cP = (\I_1, \ldots, \I_k)$ is a linearly ordered partition 
of $\{1, \ldots, n\}$ and $\sigma: \{1, \ldots, n \} \to \{\pm 1\}$ is a choice of signs. 

\np
Two pairs $(\cP', \sigma')$ and $(\cP'', \sigma'')$
determine the same parabolic subalgebra if and only if $\cP'$ and $\cP''$ are the same ordered 
partitions whose maximal part $\I_0$ contains one element and $\sigma'$ and $\sigma''$ 
coincide on $\{1, \ldots, n\} \backslash \I_0$.
\end{minipage} \\

\bigskip
{\bf Type II:} 
\begin{minipage}[t]{0.75\textwidth}
pairs $(\cP, \sigma)$, where $\cP = (\I_0, \I_1, \ldots, \I_k)$ is a linearly ordered partition 
of $\{1, \ldots, n\}$ with largest element $\I_0$ such that $|\I_0| \geq 2$ and 
$\sigma: \{1, \ldots, n \} \backslash \I_0 \to \{\pm 1\}$ is a choice of signs. 
\end{minipage}\\
\end{centering}
\end{enumerate}

\np
\underline{In Type I:}

\begin{enumerate}
\setcounter{enumi}{2}
\item The roots of $\gp_{(\cP,\sigma)}$ are  
$\{ \sigma(i) \vep_i - \sigma(j) \vep_j \, | \, i \neq j, \cP(i) \preceq \cP(j) \} \cup \{\sigma(i) \vep_i + \sigma(j) \vep_j \, | \, i \neq j\}$ 
\item The roots of $\gs_{(\cP,\sigma)}$ are $\{ \sigma(i)\vep_i - \sigma(j)\vep_j \, | \, i \neq j, \cP(i) = \cP(j) \}$.
\item $\gs_{(\cP,\sigma)} = \oplus_i \, \gs_{(\cP,\sigma)}^i$, where $\gs_{(\cP,\sigma)}^i \cong gl_{|\I_i|}$.
\item The Cartan subalgebra of $\gs_\cP^i$ is spanned by $\{\sigma(j) h_j\}_{j \in \I_i}$ 
\item The roots of $\gs_{(\cP,\sigma)}^i$ are $\{\sigma(j)\vep_j - \sigma(l)\vep_l \, | \, j \neq l \in \I_i\}$.  
\item 
$\gt_{(\cP,\sigma)}$ has a basis $\{t_1, \ldots, t_k\}$ with $t_i = \frac{1}{|\I_i|} \sum_{j \in \I_i} \sigma(j) h_j$.
\end{enumerate}

\noindent
If $\{\delta_1, \ldots, \delta_k\}$ is the basis of $\gt^*$ dual to $\{t_1, \ldots, t_k\}$ then

\begin{enumerate}
\setcounter{enumi}{8}
\item 
$\Rsh = \{ \pm \delta_i \pm \delta_j | \, 1 \leq i \neq j \leq k\} \cup \{ 2 \delta_i \, | \, |\I_i| >1\}.$
\item For $\nu \in \Rsh$, 

\begin{enumerate}
\item $\ggg^\nu\cong \V_i^\pm \otimes \V_j^\pm$ if  $\nu = \pm \delta_i \pm \delta_j$.
\item $\ggg^\nu\cong\Lambda^2 \V_i^\pm$ if $\nu = \pm 2 \delta_i$.
\end{enumerate}

\np
where $\V_i^+$ and $\V_i^-$ are the natural $\gs_{(\cP,\sigma)}^i$-module and its dual, and
all other factors of $\gs_{(\cP,\sigma)}$ act trivially.

\medskip
\item
Every parabolic subalgebra of $\ggg$ whose reductive part is $\gs_{\cP,\sigma}$ corresponds to a pair $(\cQ, \tau)$ 
such that the parts of $\cQ$ are the same as the parts of $\cP$ and $\sigma_{|\I_i} = \pm \tau_{|\I_i}$ for every part $\I_i$ or, 
equivalently, to a linear order on the set $\{\pm \delta_1, \ldots, \pm \delta_k\}$ compatible with multiplication by $-1$. 
Note that this correspondence is not bijective since
two different linear orders may determine the same parabolic subalgebra.
\end{enumerate}

\np
\underline{In Type II:}

\begin{enumerate}
\setcounter{enumi}{2}
\item Roots of $\gp_{(\cP,\sigma)} = $

\medskip
$\left\{ \sigma(i) \vep_i - \sigma(j) \vep_j \, \st  \,  i \neq j, \cP(i) \preceq \cP(j) \prec \I_0\right\} 
\cup \{\pm \vep_i \pm \vep_j \, | \, i \neq j \in \I_0\} $

\medskip
\hfill
$\cup
\{\sigma(i) \vep_i + \sigma(j) \vep_j \, | \, i \neq j \not \in \I_0\} \cup 
\{\sigma(i) \vep_i \pm \vep_j \, | \, i \not \in \I_0, j \in \I_0\}$ 

\medskip
\item Roots of $\gs_{(\cP,\sigma)}= 
\{ \sigma(i)\vep_i - \sigma(j)\vep_j \, | \, i \neq j, \cP(i) = \cP(j) \prec \I_0 \} \cup \{\pm \vep_i \pm \vep_j \, | \, i \neq j \in \I_0\}.$ 

\medskip
\item $\gs_{(\cP,\sigma)} = \oplus_{i=0}^k \, \gs_{(\cP,\sigma)}^i$, where $\gs_{(\cP,\sigma)}^0 \cong D_{|\I_0|}$ and 
$\gs_{(\cP,\sigma)}^i \cong gl_{|\I_i|}$ for $i>0$.
\item 
Cartan subalgebra of $\gs_{(\cP,\sigma)}^i$ is spanned by $\{h_j\}_{j \in \I_0}$ for $i=0$ and by 
$\{\sigma(j) h_j\}_{j \in \I_i}$ for $i>0$; 

\medskip
\item 
roots of $\gs_{(\cP,\sigma)}^i$ are
$\{\pm \vep_j \pm \vep_l \, | \, j \neq l \in \I_0\}$ for $i = 0$ and
$\{\sigma(j)\vep_j - \sigma(l)\vep_l \, | \, j \neq l \in \I_i\}$ for $i>0$.  

\medskip
\item 
$\gt_{(\cP,\sigma)}$ has a basis $\{t_1, \ldots, t_k\}$ with $t_i = \frac{1}{|\I_i|} \sum_{j \in \I_i} \sigma(j) h_j$.
\end{enumerate}

\noindent
If $\{\delta_1, \ldots, \delta_k\}$ is the basis of $\gt^*$ dual to $\{t_1, \ldots, t_k\}$ then

\begin{enumerate}
\setcounter{enumi}{8}
\item 
$\Rsh = \{ \pm \delta_i \pm \delta_j, \pm \delta_i \, | \, 1 \leq i \neq j \leq k\} \cup \{\pm 2 \delta_i \, | \, |\I_k| >1\}.$

\item For $\nu \in \Rsh$, 

\begin{enumerate}
\item $\ggg^\nu\cong \V_i^\pm \otimes \V_j^\pm$ if $\nu = \pm \delta_i \pm \delta_j$,
\item $\ggg^\nu\cong\V_i^\pm \otimes \V_0$ if $\nu = \pm \delta_i$, 
\item $\ggg^\nu\cong\Lambda^2 \V_i^\pm$ if $\nu = \pm 2 \delta_i$,
\end{enumerate}

\np
where $\V_i^+$ and $\V_i^-$ denote the natural $\gs_{(\cP,\sigma)}^i$-module and its dual 
for $i>0$, 
$\V_0$ is the natural  $\gs_{(\cP,\sigma)}^0$-module, and where 
all other factors of $\gs_{(\cP,\sigma)}$ act trivially.

\medskip
\item
The parabolic subalgebras of $\ggg$ whose reductive part is $\gs_{(\cP,\sigma)}$ are in a bijection with the pairs $(\cQ, \tau)$ 
such that the parts of $\cQ$ are the same as the parts of $\cP$ and $\sigma_{|\I_i} = \pm \tau_{|\I_i}$ for every part $\I_i$ or, 
equivalently, with linear orders on the set $\{\pm \delta_1, \ldots, \pm \delta_k\}$. 
\end{enumerate}

\np
{\bf Proof of Theorem \ref{thm:classical}.}
The idea is simple: using $\Ssh$ we will define a partial order on the set $\{\delta_1, \ldots, \delta_k\}$
(respectively on $\{\pm \delta_1, \ldots, \pm \delta_k\}$) and using the fact that $\dim(\Sym^\cdot(\Msh))^\gs = 1$ 
we will prove
that this order can be extended to a linear order (respectively, to a linear order compatible with multiplication by $-1$).
We show this case--by--case.

\noindent
\underline{Case 1: $\ggg = gl_n$.} 

\np
Set

\begin{equation} \label{eq2}
\delta_i \prec \delta_j \quad  {\text { if }} \quad \nu = \delta_i - \delta_j \in \Ssh.
\end{equation}
 We claim that
$\prec$ can be extended to a linear order on the set $\{\delta_1, \ldots, \delta_k\}$. Assume not. Then there is a sequence
$\nu_1 = \delta_{i_1} - \delta_{i_2}, \nu_2 = \delta_{i_2} - \delta_{i_3}, \ldots, \nu_l = \delta_{i_l} - \delta_{i_1}$ of elements
of $\Ssh$. We may assume that this sequence is the shortest among all such sequences, and hence that the $\gs$--modules
$\ggg^{\nu_1}, \ldots, \ggg^{\nu_l}$ are distinct submodules of $\Msh$.  As a consequence, 
$\ggg^{\nu_1} \otimes \ldots \otimes \ggg^{\nu_l}$ is a submodule of $\Sym^\cdot(\Msh)$. On the other hand,

\begin{equation} \label{eq12}
\rule{0.5cm}{0cm}\ggg^{\nu_1} \otimes \ldots \otimes \ggg^{\nu_l} \cong (\V_{i_1} \otimes \V_{i_2}^*) \otimes \ldots \otimes (\V_{i_l} \otimes \V_{i_1}^*) \cong
(\V_{i_1} \otimes \V_{i_1}^*) \otimes \ldots \otimes (\V_{i_l} \otimes \V_{i_l}^*),
\end{equation}
where the lower index of a module shows which component of $\gs_\cP$ acts non--trivially on it. 
The right hand side module in (\ref{eq12}) clearly contains the trivial $\gs_\cP$--module which contradicts the assumption that
$\dim(\Sym^\cdot(\Msh))^\gs = 1$. This contradiction shows that the order (\ref{eq2}) can be extended to a linear order on the 
set $\{\delta_1, \ldots, \delta_k\}$, which completes the proof in this case.

\bigskip
\noindent
\underline{Case 2: $\ggg = \B_n$, Type I.} 

\np
Set

\begin{equation} \label{eq22} 
\begin{array}{lll}
s_i \delta_i \prec s_j \delta_j  {\text { for }} i \neq j & {\text { if }} &  \nu = s_i \delta_i - s_j \delta_j \in \Ssh \\
s_i \delta_i \prec - s_i \delta_i & {\text { if }} & \nu = s_i \delta_i \in \Ssh,
\end{array}
\end{equation}
where $s_i, s_j = \pm 1$. We claim that
$\prec$ can be extended to a linear order on the set $\{\pm \delta_1, \ldots, \pm \delta_k\}$ compatible with multiplication by $-1$. 
Assume not. Then there is a cycle

\begin{equation} \label{eq31}
s_1 \delta_{i_1} \prec s_2 \delta_{i_2} \prec \ldots \prec s_l \delta_{i_l} \prec s_1 \delta_{i_1}
\end{equation}
of elements of $\{\pm\delta_1,\ldots,\pm\delta_k\}$
We may assume that we have chosen this cycle to be of minimal length. It gives rise to 
a sequence $\nu_1, \ldots, \nu_l \in \Ssh$ induced from (\ref{eq22}). More precisely, 
$$
\nu_j = \begin{cases} s_j \delta_{i_j} - s_{j+1} \delta_{i_{j+1}} & {\text{ if }}  \delta_{i_j} \neq \delta_{i_{j+1}}\\
s_j \delta_{i_j} & {\text { if }}  \delta_{i_j} = \delta_{i_{j+1}},
\end{cases}
$$
where $s_{l+1} = s_1$ and $\delta_{i_{l+1}} = \delta_{i_1}$. Notice that the minimality of the cycle (\ref{eq31})
implies that $\delta_{i_j} = \delta_{i_{j+1}}$ for at most two values of $j$.

\np
If $\delta_{i_j} \neq \delta_{i_{j+1}}$ for every $j$, then, as in Case 1,  
$\ggg^{\nu_1} \otimes \ldots \otimes \ggg^{\nu_l}$ is a submodule of $\Sym^\cdot(\Msh)$ and
$$
\ggg^{\nu_1} \otimes \ldots \otimes \ggg^{\nu_l} \cong (\V_{i_1}^{s_1} \otimes \V_{i_2}^{-s_2}) 
\otimes \ldots \otimes (\V_{i_l}^{s_l} \otimes \V_{i_1}^{-s_1}) \cong
(\V_{i_1}^{s_1} \otimes \V_{i_1}^{-s_1}) \otimes \ldots \otimes (\V_{i_l}^{s_l} \otimes \V_{i_l}^{-s_l})
$$
 contains the trivial $\gs_\cP$--module which contradicts the assumption that
$\dim(\Sym^\cdot(\Msh))^\gs = 1$. 

\np
If $\delta_{i_1} = \delta_{i_l}$ (which we may assume in case there is only one $j$ with $\delta_{i_j} = \delta_{i_{j+1}}$),
then $s_l = - s_1$. Furthermore,
$\ggg^{\nu_1} \otimes \ldots \otimes \ggg^{\nu_{l-1}}  \otimes \Sym^2 \ggg^{\nu_l}$ is a submodule of $\Sym^\cdot(\Msh)$ and
\begin{align*}
& \ggg^{\nu_1} \otimes \ldots \otimes \ggg^{\nu_{l-1}}  \otimes \Sym^2 \ggg^{\nu_l}  \cong (\V_{i_1}^{s_1} \otimes \V_{i_2}^{-s_2}) 
\otimes \ldots \otimes (\V_{i_{l-1}}^{s_{l-1}} \otimes \V_{i_1}^{s_1}) \otimes (\Sym^2 \V_{i_1}^{-s_1}) \\ & \cong
((\V_{i_1}^{s_1})^{\otimes 2} \otimes \Sym^2 \V_{i_1}^{-s_1}) \otimes (\V_{i_2}^{s_2} \otimes \V_{i_2}^{-s_2}) 
\otimes \ldots \otimes (\V_{i_{l-1}}^{s_{l-1}} \otimes \V_{i_{l-1}}^{-s_{l-1}})
\end{align*}
 contains the trivial $\gs_\cP$--module which contradicts the assumption that
$\dim(\Sym^\cdot(\Msh))^\gs = 1$. 

\np
Finally, if $\delta_{i_j} = \delta_{i_{j+1}}$ for two values of $j$ we proceed as above. The only difference is that we will
obtain a term analogous to $(\V_{i_1}^{s_1})^{\otimes 2} \otimes \Sym^2 \V_{i_1}^{-s_1}$ twice. The proof is complete in this case.

\bigskip
\noindent
\underline{Case 3: $\ggg = \B_n, \C_n$ or $\D_n$, Type II.} 

\np
This case is analogous to Case 2. The only difference is that if $\nu = s_i \delta_i$
then $\ggg^\nu \cong \V_i^{s_i} \otimes \V_0$ and consequently
$$
\Sym^2 \ggg^{\nu} \cong \Sym^2 \V_i^{s_i} \otimes \Sym^2 \V_0 \oplus \Lambda^2 \V_i^{s_i} \otimes \Lambda^2 \V_0.
$$
The observation that the trivial $\gs_{(\cP, \sigma)}^0$--module is contained in $(\Sym^2 \V_0)^{\gs}$ for $\ggg = \B_n, \D_n$ and in
$\Lambda^2 \V_0$ for $\ggg=\C_n$ concludes the proof.

\noindent
\underline{Case 4: $\ggg = \C_n$ or $\D_n$, Type I.} 

\np
This case is again analogous to Case 2. The only difference is that we need to modify the 
second line of formula (\ref{eq22}) to read

\begin{equation} \label{eq41}
s_i \delta_i \prec - s_i \delta_i \quad {\text { if }} \quad \nu = s_i 2 \delta_i \in \Ssh.
\end{equation}
Consequently, for $\nu$ as in (\ref{eq41}) we have $\ggg^\nu \cong \Sym^2 \V_i^{s_i}$ for $\ggg=\C_n$ and
$\ggg^\nu \cong \Lambda^2 \V_i^{s_i}$ for $\ggg = \D_n$ and proceeding as
in Case 2 we find a non--trivial invariant in the module $\ggg^{\nu_1} \otimes \ldots \otimes \ggg^{\nu_l} \subset (\Sym^\cdot\Msh)^{\gs}$.
\hfill $\square$

\np
{\bf Exceptional Lie algebras.}
We first note that Theorem \ref{thm:classical} holds for $\ggg = \G_2$. 
Indeed, if $\gs$ is a proper subalgebra of $\ggg$ which is not equal to
$\gh$, then all elements of $\Rsh$ are proportional and there is nothing to prove. If $\gs = \gh$, then the spaces $\ggg^\nu$ are just
the root spaces of $\ggg$ which are one dimensional and again the statement is clear.

\np
The proposition, however, fails to hold for the remaining exceptional Lie algebras. Here is a uniform example.
Let $\ggg = \F_4, \E_6, \E_7$, or $\E_8$. Denote the rank of $\ggg$ by $n$. 
Then $\ggg$ has a reductive subalgebra $\gm \oplus \gc$ of rank $n$ such that $\gm \cong \A_1$ and 
$\gc \cong \C_3, \A_5, \D_6$, or $\E_7$ respectively. As an $\gm$--module $\ggg$ decomposes as

\begin{equation} \label{eq55}
\ggg = (\Ad_\gm \otimes \tr_\gc) \oplus (\tr_\gm \otimes \Ad_\gc) \oplus (\V \otimes \U),
\end{equation}
where $\Ad_\gm$ and $\Ad_\gc$ are the adjoint modules of $\gm$ and $\gc$ respectively; $\tr_\gm$ and $\tr_\gc$ ---the respective trivial modules;
$\V$ is the natural $\gm \cong \A_1$--module; and the weights of $\U$ belong to a single $\W$-orbit.

\np
Set $\gs = \gm + \gh$. From the construction of $\gs$ we conclude that $\gt = \gh_\gc$, the Cartan subalgebra of $\gc$. Furthermore,
 (\ref{eq55}) implies $\Rsh = \Delta_\gc \cup \supp \U$, where $\supp \U$ denotes the set of weights 
 of $\U$ and, for $\nu \in \Rsh$ the 
 $\gs = \gm \oplus \gh_\gc$--module $\ggg^\nu$ is given by
 $$
 \ggg^\nu \cong \begin{cases}
 \tr_\gm \otimes \nu & {\text { if }} \nu \in \Delta_\gc\\
 \V \otimes \nu & {\text { if }} \nu \in \supp U.
 \end{cases}
 $$
Pick $\nu_0 \in \supp \U$ and write $\nu_0 = q_1 \nu_1 + \ldots + q_r \nu_r$ with $q_i \in \QQ_+$, $\nu_i \in \Delta_\gc^{+}$. Set
$\Ssh = \{ -\nu_0, \nu_1, \ldots, \nu_r\}$ and $\Msh=\opp_{\nu\in\Ssh} \ggg^{\nu}$.
It is clear from the construction that $\dim\Sym^\cdot(\Msh)^{\gs} = 1$ but no parabolic subalgebra $\gp$ as in 
Theorem \ref{thm:classical} exists.

\end{document}